\documentclass{article}
\makeatletter
\ifdefined\tikzexternalrealjob
  \@namedef{ver@hyperref.sty}{9999/99/99 latex-pipeline: real hyperref skipped during TikZ externalize}

  \providecommand{\url}[1]{\texttt{#1}}

  \providecommand{\hypersetup}[1]{}
  
  \providecommand{\pdfbookmark}[3][0]{}
  \providecommand{\texorpdfstring}[2]{#1}
  \let\latexpipeline@origlatexerror\@latex@error
  \let\@latex@error\@gobbletwo
  \AtBeginDocument{\global\let\@latex@error\latexpipeline@origlatexerror}\fi
\makeatother
\usepackage{tikz}
\usetikzlibrary{external}
\usepackage{amsmath}
\usepackage{amssymb}
\usepackage{graphicx}
\usepackage{natbib}
\usepackage{tikz}
\usepackage{endnotes}
\usepackage{amsthm}
\usepackage{pipeline-authoring}
\usepackage{informs4-compat}
\usepackage[margin=1in]{geometry}

\usepackage{snapshot}

\DeclareMathOperator{\proj}{proj}
\DeclareMathOperator{\epi}{epi}

\DeclareMathOperator{\conv}{conv}

\DeclareMathOperator{\graph}{gr}

\DeclareMathOperator{\vertex}{vert}

\DeclareMathOperator{\hypo}{hyp}

\DeclareMathOperator{\R}{\mathbb{R}}
\DeclareMathOperator{\Z}{\mathbb{Z}}

\DeclareMathOperator{\for}{for}

\DeclareMathOperator{\lnatural}{L^\natural}

\usepackage{comment}
\usepackage{bbm}
\usepackage{rotating}

\usepackage{tikz, tikz-3dplot}
\usepackage{pgfplots}
\usepackage{tkz-fct}

\usetikzlibrary{positioning,calc}
\usepackage{etoolbox}
\usepackage{xcolor}
\usepackage[default]{mycolor}
\usepackage[colorlinks=true,breaklinks=true,bookmarks=true,urlcolor=links,citecolor=links,linkcolor=links,bookmarksopen=false,draft=false]{hyperref}

\usepackage[ruled]{algorithm2e}
\usepackage{amsmath,amsfonts}
\usepackage{caption}
\usepackage{subcaption}

\usepackage{endnotes}
\let\footnote=\endnote
\let\enotesize=\normalsize

\definecolor{links}{RGB}{204,36,29}
\usepackage[colorlinks=true,breaklinks=true,bookmarks=true,urlcolor=links,citecolor=links,linkcolor=links,bookmarksopen=false,draft=false]{hyperref}

\usepackage{xcolor}

\NewDocumentEnvironment{diffblock}{ +O{} }
  {
    \par\medskip
        \if\relax\detokenize{#1}\relax\else
      {\color{red} #1 \par\medskip}
    \fi
          }
  {
    \par\medskip
  }

\usepackage{natbib}
\usepackage{multirow}
\usepackage{booktabs}
 \bibpunct[, ]{(}{)}{,}{a}{}{,} \def\bibfont{\small}   \def\newblock{\ } 
\TheoremsNumberedThrough     \ECRepeatTheorems

\EquationsNumberedThrough     \renewcommand{\theARTICLETOP}{}

\begin{document}

\pipelinetitle{Discreteness to Convexity:\\ An Application to Promotion Planning}
\pipelinetitle{Discreteness to Convexity:\\ Promotion Planning via Simplotope Triangulation}
\pipelineauthor{Taotao He}{Antai School of Management, Shanghai Jiao Tong University}{hetaotao@sjtu.edu.cn}
\pipelineauthor{Mohit Tawarmalani}{Mitch Daniels School of Business, Purdue University}{mtawarma@purdue.edu}

\pipelineabstract{Price promotion optimization is a computationally challenging problem central to supermarket operations, requiring simultaneous pricing decisions across multiple products and periods. This paper introduces a new formulation for price promotion by developing convex hull results for supermodular compositions of univariate functions over a simplotope.  Leveraging this reformulation with Gurobi, we achieve substantial performance gains: instances with up to 125 products, 20 periods, and 5 price levels are solved in an average of 7 minutes, demonstrating the potential to handle even larger instances.  Our exact solution methods extract 25--48\% additional profit from promotion planning relative to state-of-the-art heuristic approaches. Additionally, we extend the polynomially solvable cases from two to multiple price levels and expand our results to allow for multiplicative historical effects. Our core methodological innovation applies to a broad class of nonlinear discrete optimization problems. Specifically, our results convexify a class of nonlinear functions that includes monomials and the widely studied $\lnatural$ function structure.
}

\pipelinekeywords{retail operations, promotion planning, submodularity, discrete convexity} 

\maketitle

\def \X {\mathcal{X}}
\def \z {\boldsymbol{z}}
\def \x {\boldsymbol{x}}
\def \blambda {\boldsymbol{\lambda}}
\def \bdelta {\boldsymbol{\delta}}
\def \p {\boldsymbol{p}}
\def \e {\boldsymbol{e}}
\def \s {\boldsymbol{s}}
\def \A {\boldsymbol{A}}
\def \b {\boldsymbol{b}}
\def \w {\boldsymbol{w}}
\def \revision {\color{blue}}
\def \revision {}

\newcommand{\drop}[1]{}
\newcommand{\q}{q}

\def \bgamma {\boldsymbol{\gamma}}
\def \y {\boldsymbol{y}}
\def \v {\boldsymbol{v}}
\def \bsigma {\boldsymbol{\sigma}}
\def \bmu {\boldsymbol{\mu}}

\def \bA {\boldsymbol{A}}
\def \j {\boldsymbol{j}}
\def\id{\text{id}}
\def \D {{\mathcal{D}}}
 \def \V{\mathcal{V}}
 \def \s {\boldsymbol{s}}
\def \r {\boldsymbol{r}}
\def\mcirc{\mathop{\circ}}
\def \B{\mathcal{B} }
\def \P{\mathcal{P}}
\def \Q{\mathcal{Q}}
\def \pt {p}
\newcommand{\G}{\mathcal{G}}
\def\App{\mathop{\text{App}}}
\def\st{\mathop{\text{s.t.}}}
\def\stair{\mathop{\text{stair}}}
\def\bi{\mathop{\text{bsc}}}
\def\rvsbi{\mathop{\text{rvs-bsc}}}
\def\rvsstair{\mathop{\text{rvs-stair}}}
\def \one{\mathop{\boldsymbol{1}}}
\def \ind {\mathbbm{1}}

\section{Introduction}
The growing interest in discrete nonlinear optimization for operations management is driven by increased data availability and advances in machine learning \citep{mivsic2020data}. We address the problem of price promotion optimization in retail operations~\citep{cohen2021promotion}, which involves determining optimal promotion levels across a finite set of items and periods to maximize profit. This problem accounts for complex historical and cross-item effects. Nonlinearity arises from modeling consumer behavior through data-driven demand functions, while discretization is imposed by fixed price promotion levels. Although multi-product demand models incorporating promotion effects have been developed~\citep{cohen2021promotion}, existing formulations {\revision can produce less profitable promotion plans}, and struggle to  scale to complex real-world scenarios. Despite extensive research in the OR/MS community on nonlinear discrete optimization and its formulations~\citep{watters1967letter,glover1975improved,oral1992linearization,murota1998discrete,adams2012base,sherali2013reformulation,vielma2015mixed,bergman2018discrete}, we identify significant opportunities to improve existing formulations by exploiting a largely overlooked structural property of the resulting binary variables: their restriction to a simplotope. These improvements enable the design of  more  {\revision profitable} and scalable price promotion campaigns and have broader implications for discrete nonlinear optimization.

Several approaches address discrete nonlinear functions, with binarization and polyhedral approximation techniques being particularly relevant to our work. Binarization introduces indicator variables for each discrete value \citep{dash2018binary}, followed by linearization of the resulting nonlinear functions to construct integer programming formulations~\citep{glover1975improved,sherali2013reformulation}. This approach is widely used in promotion planning \citep{cohen2017impact,cohen2020optimizing,cohen2021promotion}. Polyhedral approximation of $\text{L}^\natural$ functions~\citep{murota1998discrete} offers another powerful tool for modeling indivisibilities and has found applications in revenue management \citep{chen2021discrete}.  Our work focuses on strengthening formulations for these problems by exploiting additional mathematical structure in the discrete domain. We demonstrate these improvements in the context of price promotion optimization, but our methodology is broadly applicable across fields including engineering design \citep{thanedar1995survey,boyd2005digital}, service operations \citep{birolini2021airline,liu2021time,lejeune2024drone,liu2024design}, and defense \citep{ahuja2007exact,he2025proactive}.

The central contribution of our work is the exploitation of additional mathematical structure in formulation design in three ways. First, we propose permuting the discretization space when the univariate functions are non-monotonic. Second, we convexify a class of functions, including $L^\natural$ functions,  in the discretization space. Third, we exploit the specialized domain in which the binary variables introduced by binarization reside. This domain, a specially structured polytope expressible as a Cartesian product of simplices, is known as a simplotope, a structure largely overlooked in formulation design. By leveraging these ideas, we derive tighter formulations and broaden the function classes with specialized convexification results. Consequently, our techniques reveal new polynomially solvable cases for price promotion optimization, particularly those involving multiple price levels and demand functions with multiplicative historic effects, and significantly expand the computationally tractable problem instances.
\subsection{Problem Description and Contributions}\label{sec:descriptionAndContributions}
{\revision
Price promotion optimization (PPO) presents a critical lever for retail profitability, yet it remains computationally intractable for large-scale instances. While
the impact of promotions on sales is well-documented~\citep{blattberg1995promotions,anderson2019price}, accurately modeling complex cross-item and
historical effects, such as cannibalization and post-promotion dips, poses significant challenges. Existing integer programming
formulations~\citep{cohen2021promotion} are exact only under restrictive assumptions (e.g., additive separability of historical effects) and suffer from
quadratic growth in auxiliary variables, severely limiting scalability.

This paper introduces a novel formulation framework that exploits the \textit{simplotope} structure inherent in binarization variables. By deriving
tight convex hull descriptions for supermodular compositions of univariate functions, we depart from prior reliance on rectangular domains. This
structural insight yields substantial computational gains: instances with up to 125 products, 20 periods, and 5 price levels are solved in an average of 7 minutes.

Our contributions are threefold:

\begin{enumerate}
    \item \textbf{Scalable and Exact Price Promotion Models.}  We develop exact reformulations for price promotion optimization that accommodate broad
demand models with multiple price levels and multiplicative historical effects. In contrast to prior approximations that suffer from quadratic growth
in auxiliary variables, our formulations are exact and scale linearly with the number of price discretization levels. Numerical results demonstrate
that our exact formulations recover an additional 25\% to 48\% of potential promotion value compared to approximations in~\citet{cohen2021promotion},
particularly under less constrained business rules.

    \item \textbf{Ideal Formulations for Submodular Compositions.} We provide the first ideal linear programming formulations for epigraphs of submodular compositions of univariate functions over the simplotope domain. We demonstrate that this framework encompasses standard formulations for $L^\natural$-convex functions over bounded domains, offering new insights into Lov\'{a}sz-extension based approaches.

    \item \textbf{Tight Relaxations for Nonlinear Discrete Functions.} Our results apply to products of two univariate functions, providing the
tightest possible linear programming relaxation for this class of nonlinear discrete functions. This includes products of powers of two discrete
variables, a common structure in revenue management and inventory optimization applications.
\end{enumerate}

The remainder of this paper is organized as follows. Section~\ref{sec:discrete_pricing} establishes the methodological framework for simplotope based formulations of PPO, and demonstrate its importance in the context of cross-item effects. Section~\ref{section:composite-MIP} presents ideal formulations for supermodular composite effects, including historical effects under log-log demand model. Section~\ref{section:ppo} integrates  these tools to obtain strong and compact formulations of PPO, and studies their theoretical properties. 
Section~\ref{section:computation} presents computational results, with Section~\ref{sec:profit-performance} quantifying the economic value of the
resulting promotion plans. Section~\ref{section:discreteconvexity} discuss connections with $\lnatural$-convex functions.

\paragraph{Notation.} We use $[n]$ to denote $\{1,\ldots,n\}$. For a set $S$, $\conv(S)$ denotes its convex hull, $\proj_{\x} S$ denotes its projection
onto the $\x$-space, and $\vertex(S)$ denotes its set of vertices. For a function $f$, $\graph(f)$, $\epi(f)$, and $\hypo(f)$ denote its graph,
epigraph, and hypograph, respectively.
}
\section{Discrete Pricing Opportunity and Methodological Framework}\label{sec:discrete_pricing}

\begin{diffblock}[
]
This section formally defines the Price Promotion Optimization (PPO) and outlines how its structure maps onto the theoretical framework developed in this paper. We present a three-stage methodology: (1) characterizing the nonlinear structure of PPO, (2) introducing the unary binarization that induces a simplotope domain, and (3) deriving tight convex hull formulations for supermodular compositions over this domain. This approach translates discrete pricing decisions into a mixed 0-1 linear optimization problem with strong integrality properties, enabling efficient solution
via linear programming relaxations under various settings. 

\paragraph{1. Nonlinear Structure of PPO.}
We formulate PPO as a discrete nonlinear program that maximizes total profit across products and periods, subject to discrete price levels and business rules:
\begin{subequations}\label{eq:ppowithrules}
    \begin{alignat}{3}
 & \max_{\x} &\quad& \sum_{i  \in [N]}\sum_{t  \in [T] } (x^i_{t} - c^i_t) \cdot D^i_t(\x) \label{eq:objppo}\\
    &&& x^i_t \in \{p^i_0, p^i_1, \ldots,p^i_{d_i} \}&& \for \, i \in [N] \text{ and } t \in [T] \label{eq:brules1} \\
        &&&\sum_{t \in [T]}\ind \{x^i_t < p^i_{d_i}\} \leq K^i&& \for\, i \in [N]\label{eq:brules2} \\
        &&&\sum_{i \in [N]}\ind \{x^i_t < p^{i}_{d_i}\} \leq C_t&& \for\, t \in [T],\label{eq:brules3}
    \end{alignat}
\end{subequations}
where $\x = (x^i_t)_{i \in [N], t \in [T]}$ denotes the vector of all price decisions. Here, $x^i_t$ is the price of item $i$ in period $t$, and
$c^i_t$ is the unit cost. Constraint \eqref{eq:brules1} restricts $x^i_t$ to a finite set of
allowed price levels, where $p^i_{d_i}$ denotes the regular (list) price.  This discrete structure reflects common retailing practices in two ways. First, retailers frequently offer promotions based on a limited number of discount levels, such as 10\%, 20\%, and 30\%. Second, discrete price ladders align with marketing and behavioral considerations, including the well-documented consumer preference for prices ending in the digit $9$~\citep{anderson2003effects}.  The indicator $\ind\{x^i_t < p^i_{d_i}\}$ equals 1 if item $i$ is discounted in period $t$, and 0 otherwise. Thus,
\eqref{eq:brules2} limits the number of periods in which item $i$ can be discounted to at most $K^i$, while \eqref{eq:brules3} ensures that in any
period $t$, at most $C_t$ items are discounted. Together, \eqref{eq:brules2} and \eqref{eq:brules3} are referred to as the \textit{business rules}.

The demand function $D^i_t(\x)$ \eqref{eq:multdemand} captures the nonlinear interactions such as:
\begin{equation}\label{eq:multdemand}
    D^i_t(\x) = a^i_t \cdot \overbrace{(x^i_t)^{-b^i_0}}^{\text{self effect}} \cdot \overbrace{\prod_{m \in [M_i]}(x^i_{t-m})^{b^i_m}}^{\text{historical effects}} + \sum_{k \in
[N]\setminus \{i\}}\overbrace{\delta^{ik}x^k_t}^{\text{cross-item effect}}.
\end{equation}
The model \eqref{eq:multdemand} will be referred to as the \emph{log-log} demand model.
Substituting this into \eqref{eq:objppo} and expanding, these effects affect the retailer's profit as follows
\begin{equation}\label{eq:multobj}
    \text{Profit} = \sum_{i \in [N]}\sum_{t \in [T]}\biggl[ a^i_t \cdot \overbrace{(x^i_t-c^i_t)\cdot (x^i_t)^{-b^i_0}}^{\text{self effect}} \cdot \overbrace{\prod_{m \in [M_i]}(x^i_{t-m})^{b^i_m}}^{\text{historical effects}} + \sum_{k \in
[N]\setminus \{i\}}\overbrace{ \delta^{ik}(x^i_t-c^i_t)x^k_t}^{\text{cross-item effect}}\biggr].
\end{equation}
where each \emph{profit term} in the above summation is a product of univariate functions in the discrete price variables. Standard linearization techniques, such as McCormick envelopes, yield weak relaxations for
these discrete functions because they treat each price level for an item $i$ in period $t$ as an independent choice, ignoring their mutual exclusivity. This results in significant integrality gaps that hinder scalability. In contrast, our approach exploits the mutual exclusivity of price levels via the simplotope structure, yielding tight relaxations and even convex hulls.

\paragraph{2. Unary Binarization and the Simplotope}
To handle the discrete nature of prices $x^i_t \in \{p^i_0, \ldots, p^i_{d_i}\}$, one could employ a unary binarization scheme. One introduces binary variables
$\lambda^i_{tj}$ such that $\lambda^i_{tj}=1$ if and only if $x^i_t = p^i_j$. The constraints $\lambda^i_{tj} \geq 0$ and $\sum_{j=0}^{d_i}
\lambda^i_{tj} = 1$ define a product of simplices, one for each $(i,t)$.

To exploit the lattice structure inherent in the price levels, we introduce auxiliary binary variables $z^i_{tj}$, where $z^i_{tj}=1$ if and only if $x^i_t \geq p^i_j$. These are derived from $\lambda$
via partial sums: $z^i_{tj} = \sum_{j'=j}^{d_i} \lambda^i_{tj'}$. The $\z$ variables satisfy chain constraints $1
= z^i_{t0} \geq z^i_{t1} \geq \cdots \geq z^i_{t,d_i} \geq 0$. This chain structure, referred to as a simplotope, allows us to
represent the discrete price $x^i_t$ as a linear
function of $\z$:
\[
    x^i_t = \sum_{j=0}^{d_i} (p^i_j - p^i_{j-1}) z^i_{tj}  \quad \text{with } \quad p^i_{-1}:=0.
\]
This lifted representation is pivotal. It enables the representation of the profit terms in \eqref{eq:multobj} as supermodular compositions
over the $z$-variables, which is the key to deriving tight convex hulls.

\paragraph{3. Convex Hulls for Supermodular Compositions}
The core insight is that each univariate function in a profit term can be written as a non-decreasing linear function of auxiliary variables that are partial sums of $\lambda^i_{tj}$ variables (see Examples~\ref{ex:reorder-ordering} and~\ref{ex:cuts}). The order of summation for defining these auxiliary variables is chosen to align with the monotonicity of the specific function component. Being partial sums, they induce a natural lattice structure and the profit term in \eqref{eq:multobj} exhibits \textit{supermodularity} over this lattice. Leveraging this, we derive the convex hull of the hypograph of supermodular compositions, which includes profit terms and monomials, over the simplotope of $\z$ variables (detailed in Theorem~\ref{them:supermodular-model}). This construction yields ideal formulations (LP relaxations with only integer extreme points) for PPO under specific sign conditions enabling exact solution via linear programming. Even when these conditions are not met, the resulting relaxations are significantly tighter than standard approaches. Specifically:
\begin{itemize}
    \item \textbf{Base Case (Two Price Levels):} Proposition~\ref{prop:simple-exact} establishes the convex hull for the base case: a bilinear term involving a binary
   variable and a discrete variable with multiple levels. In the PPO context, this provides the exact convex hull for interactions involving a
binary decision (e.g., promote vs. no promote) and a discrete price level. Example~\ref{ex:mcweak} demonstrates that the traditional rectangular-domain formulations~\citep{cohen2021promotion} fail to be ideal even in
this simple setting.
    \item \textbf{General Supermodular Compositions:} Proposition~\ref{prop:bilinear-sub} extends this construction to the product of two arbitrary univariate functions. For the cross-item effect, this provides the convex hull for multiple price levels. Moreover, it provides the convex hull for the historical effect when $M_i=1$, \textit{i.e.}, when the historical effect depends only on the previous period's price. In Examples~\ref{ex:reorder}-\ref{ex:cuts}, we show that for products of functions and, more specifically, for modeling multiplicative historical effects in \eqref{eq:ppowithrules}, the transformation given in Theorem~\ref{them:supermodular-model} is essential for obtaining an ideal formulation.
    \item \textbf{Global Optimality via LP Relaxation:} Theorem~\ref{thm:ppo-lp} demonstrates that under certain sign conditions on the demand parameters, the PPO problem with the demand model specified in \eqref{eq:multdemand} can be solved exactly via the LP relaxation of our formulation. This result exploits the supermodularity of the objective function~\eqref{eq:multobj} over an appropriately defined lattice.
\end{itemize}
This framework not only provides a tighter relaxation for PPO but also generalizes to other discrete nonlinear optimization problems involving supermodular functions. The following sections detail the core ingredients of this framework, starting with the simplotope structure and its application to cross-item effects.
\end{diffblock}
\subsection{Discrete domains with unary binarization lead to simplotope structure}\label{section:binarization}
A prevalent technique to deal with discreteness requirement on a variable is to represent it using new binary variables. Binarization schemes have been studied in~\citet{glover1975improved,sherali2013reformulation,owen2002value,roy2007binarize}, and \citet{dash2018binary}.  They are useful since imposing integrality on binary variables, rather than on general discrete variables, has some attractive theoretical properties with respect to cutting planes and to smaller branch-and-bound trees, see~\cite{dash2018binary}. {\revision Binarization schemes also introduce a special polyhedral structure. Since our formulation exploits this structure to linearize the nonlinear profit function in PPO, we briefly discuss binarization schemes and related polyhedral results.}

Our formulations in Section~\ref{section:composite-MIP} will be built on a particular binarization, typically referred to as  \emph{unary binarization}~\citep{roy2007binarize}.  Consider a  discrete domain $\X$ in $\R^n$ given as $\prod_{i\in[n]}
\{p_{i0},p_{i1},\ldots,p_{id_i}\}$, where $p_{i0} <\cdots <p_{id_i}$ and $d_i \geq 1$. In the context of PPO, this models the price levels of multiple items across multiple periods. To model the constraint $\x=(x_1,\ldots,x_n)\in\X$, it suffices to treat each variable individually. For $i\in[n]$, a discrete variable $x_i \in \{\pt_{i0}, \pt_{i1},\ldots, \pt_{id_i}\}$ can be binarized as follows:
\begin{equation}\label{eq:unary}
x_i = \pt_{i0} + \sum_{ j \in [d_i]} (\pt_{ij} - \pt_{ij-1}) z_{ij} \quad   1\geq z_{i1} \geq z_{i2} \geq \cdots \geq z_{id_i} \geq 0 \quad \text{ and } \quad   \z_i \in  \{0,1\}^{d_i}.\tag{\textsc{Unary}}
\end{equation}
Note that the binary vectors $\z_i$ are precisely the vertices of a $d_i$-dimensional simplex, which we refer to as the \textit{unary simplex}:
\[
\Delta^{d_i} = \bigl\{\z \in \R^{d_i} \bigm| 1\geq z_{i1} \geq z_{i2} \geq \cdots \geq z_{id_i} \geq 0 \bigr\}.
\]

 Later, in Sections~\ref{section:otherbinarization} and~\ref{section:computation}, we extend our formulations to accommodate alternative binarizations, including the \textit{full binarization}~\citep{sherali2013reformulation}, which represents $x_i$ as:
\begin{equation}\label{eq:fullbinarization}
	x_i = \sum_{j=0}^{d_i}p_{ij} \lambda_{ij} \qquad \lambda_{i} \in \{0,1\}^{d_i+1} \quad \text{ and } \quad \sum_{j = 0}^{d_i} \lambda_{ij} = 1, \tag{\textsc{Full}}
\end{equation}
where note that the binary vector $\blambda_i$ lies in the $d_i$-dimensional standard simplex. The unary and full binarizations are related via a one-to-one affine transformation given as follows. For each $i \in [n]$, define the mapping $\Lambda_i : \z_i \in \R^{d_i} \mapsto \blambda_i \in \R^{d_i+1}$ by
\begin{equation}\label{eq:fromztolambda}
\lambda_{ij} = z_{ij} - z_{ij+1} \quad \text{for } j \in 0 \cup [d_i],
\end{equation}
where we set $z_{i0} = 1$ and $z_{i(d_i+1)} = 0$. The inverse transformation $\Lambda_i^{-1} : \blambda_i \mapsto \z_i$ is given by:
\[
z_{ij} = \sum_{k =j}^{d_i}\lambda_{ik} \quad \text{ for } j \in [d_i].
\]
We note that \cite{dash2018binary} study the behavior of binarizations with respect to the split closure and show that a binarization is optimal in this sense if it belongs to the family of unimodular binarizations. Both binarizations are unimodular, providing theoretical justification for their use in our formulations. 

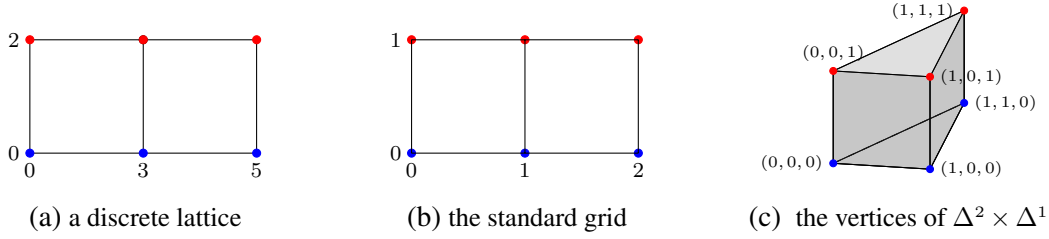
\begin{figure}
\centering
  \tdplotsetmaincoords{70}{370}
\begin{subfigure}[b]{0.3\linewidth}
\centering
\tikzsetnextfilename{promotion-1}
\begin{tikzpicture}[scale=1.5]
\begin{scope}
    \node[below] at (0 ,0) {\scriptsize $0$};
    \node[below] at (1 ,0) {\scriptsize $3$};
    \node[below] at (2 ,0) {\scriptsize $5$};
    
    \node[left] at (0, 0) {\scriptsize $0$};
    \node[left] at (0, 1) {\scriptsize $2$};
\filldraw[black] (1,1) circle (1pt) ; 

    \draw (0, 0) grid (2, 1);
    \filldraw[red] (1,1) circle (1pt) ; 
\filldraw[blue] (0,0) circle (1pt) ; 
\filldraw[red] (0,1) circle (1pt) ; 
\filldraw[red] (2,1) circle (1pt) ; 
\filldraw[blue] (2,0) circle (1pt) ; 
\filldraw[blue] (1,0) circle (1pt) ; 

\end{scope}
\end{tikzpicture}
  \caption{{\small a discrete lattice}}
\end{subfigure}
\begin{subfigure}[b]{0.3\linewidth}
\centering
\tikzsetnextfilename{promotion-2}
\begin{tikzpicture}[scale=1.5]
\begin{scope}
    \node[below] at (0 ,0) {\scriptsize $0$};
    \node[below] at (1 ,0) {\scriptsize $1$};
    
    \node[below] at (2 ,0) {\scriptsize $2$};
    \node[left] at (0, 0) {\scriptsize $0$};
    \node[left] at (0, 1) {\scriptsize $1$};
\filldraw[red] (1,1) circle (1pt) ; 
\filldraw[red] (0,1) circle (1pt) ; 
\filldraw[blue] (0,0) circle (1pt) ; 
\filldraw[red] (2,1) circle (1pt) ; 
\filldraw[blue] (2,0) circle (1pt) ; 
\filldraw[blue] (1,0) circle (1pt) ; 

    \draw (0, 0) grid (2, 1);
\end{scope}
\end{tikzpicture}
  \caption{{\small the standard grid}}
\end{subfigure}
\begin{subfigure}[b]{0.3\linewidth}
\centering
  \tikzsetnextfilename{promotion-3}
\begin{tikzpicture}[scale=1.3, tdplot_main_coords]
        \definecolor{point_color}{rgb}{ 0,0,0 }
    \tikzstyle{point_style} = [fill=point_color]

    \coordinate (q00) at (0, 0, 0);
    \coordinate (q10) at (1, 0, 0);
    \coordinate (q20) at (1, 2, 0);
    \coordinate (q01) at (0, 0, 1);
    \coordinate (q11) at (1, 0, 1);
    \coordinate (q21) at (1, 2, 1);

        \definecolor{edge_color}{RGB}{100,100,100}

      \definecolor{plane_color}{RGB}{249, 6, 33}
    \definecolor{shadow_color}{RGB}{100,100,100}

    \tikzstyle{plane_style} = [fill=plane_color, fill opacity=0.8, draw=plane_color, line width=1 pt, line cap=round, line join=round] 
    \tikzstyle{shadow_style} = [fill=shadow_color, fill opacity=0.2, line width=0.05 pt, line cap=round, line join=round]
\draw[black,dashed] (q00) -- (q20);
\draw[black] (q00) -- (q10);
\draw[black] (q20) -- (q10);
\draw[black] (q00) -- (q01);
\draw[black] (q11) -- (q01);
\draw[black] (q11) -- (q21);
\draw[black] (q11) -- (q10);
\draw[black] (q20) -- (q21);
\draw[black] (q01) -- (q21);

      \draw[shadow_style] (q00) -- (q10) -- (q20) -- (q00) -- cycle;
	  \draw[shadow_style] (q00) -- (q01) -- (q21) -- (q20) -- (q00)-- cycle;
	  \draw[shadow_style] (q00) -- (q01) -- (q11) -- (q10) -- (q00)-- cycle;
	  \draw[shadow_style] (q10) -- (q11) -- (q21) -- (q20) -- (q10)-- cycle;

\node [left] at (q00) {\tiny{$(0,0,0)$}};
\node [right] at (q10) {\tiny{$(1,0,0)$}};
\node [right] at (q20) {\tiny{$(1,1,0)$}};

\node [above] at (q01) {\tiny{$(0,0,1)$}};
\node [left] at (q21) {\tiny{$(1,1,1)$}};
\node [right] at (q11) {\tiny{$(1,0,1)$}};

\filldraw[blue] (q00) circle (1pt) ; 
\filldraw[blue] (q10) circle (1pt); 
\filldraw[blue] (q20) circle (1pt) ; 
\filldraw[red] (q01) circle (1pt) ; 
\filldraw[red] (q11) circle (1pt); 
\filldraw[red] (q21) circle (1pt) ; 
  \end{tikzpicture}
  \caption{{ \small the vertices of $\Delta^2 \times \Delta^1$}}
\end{subfigure}
\caption{ The grid representation of a discrete lattice and vertices of a simplotope}\label{fig:grid-rep}
\end{figure}

\begin{figure}[h]
    \begin{subfigure}[b]{0.5\linewidth}
        \centering
        \includegraphics[page=2, scale=0.6]{gridvisualize.pdf}
        \caption{Function: $\max\{5x+1,10x-14\}\times (y+1)$}
    \end{subfigure}
    \begin{subfigure}[b]{0.5\linewidth}
        \centering
        \includegraphics[page=1, scale=0.6]{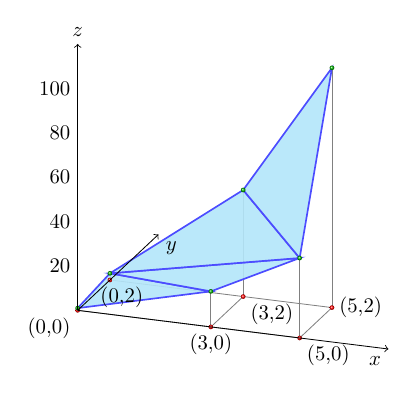}
        \caption{Convex Envelope}
    \end{subfigure}
    \caption{Grid: $\{0,3,5\}\times\{0,2\}$. Function (resp. convex envelope) at $(x,y) = (2,1)$ evaluate to $22$ (resp. $12$).}\label{fig:gridfuncenv}
\end{figure}

For our theoretical developments, we will represent the possible discrete values of $\x$ using the standard grid $\G := \prod_{i=1}^n\{0,\ldots,d_i\}$. To do so, let $j_i$ be the index of the value that variable $x_i$ takes, so that $j_i = k$ when $x_i = p_{ik}$. Thus, each grid point $(j_1,\ldots,j_n)$ corresponds to a unique setting of $\x$. 
This grid structure will also allow us to relate $\x=(x_1, \ldots, x_n)$ to its lifting  $\z = (\z_1,\ldots,\z_n)$. To see this, observe that $\z$ belongs to the following Cartesian product of simplices, referred to as the \textit{unary simplotope}:
\[
\Delta : = \Delta^{d_1} \times \Delta^{d_2} \times \cdots \times 	\Delta^{d_n}.
\] 
The grid-points of $\G$, \textit{a.k.a.} feasible $\x$ vectors, are in one-to-one correspondence with the vertices of $\Delta$, which can be written as:
\[
\vertex(\Delta) = \vertex(\Delta^{d_1})  \times \vertex(\Delta^{d_2}) \times \cdots \times 	 \vertex(\Delta^{d_n}) = \Delta \cap \{0,1\}^{\sum_{i \in [n]}d_i}.
\]
The vertex set of $\Delta^{d_i}$ consists of $d_i+1$ elements, $\v_{i0}, \v_{i1},\ldots,\v_{id_i}$, where 
\begin{equation}\label{eq:verticesDelta}
\v_{i0} = (0, \ldots, 0) \quad \text{ and } \quad \v_{ij} = \v_{ij-1} + \e_{ij} \quad \for  j \in [d_i], \tag{\textsc{Vert}}
\end{equation}
and $\e_{ij}$ is the $j^\text{th}$ principal vector in the space spanned by variables $(z_{i1}, \ldots, z_{id_i})$. In particular, these correspondences are:
\[
\x =(p_{1j_1},\ldots,p_{nj_n})  \Leftrightarrow \j=(j_1,\ldots,j_n) \Leftrightarrow \z = (\v_{1j_1},\ldots,\v_{nj_n}).
\]
See Figure~\ref{fig:grid-rep} for an illustration with $n=2$ and $(x_1,x_2)\in\{0,3,5\}\times\{0,2\}$. Our interest will be in functions defined over this region and their convex envelopes (tightest convex underestimators). See Figure~\ref{fig:gridfuncenv} for an illustration of a function whose envelope consists of four linear inequalities. We will focus on describing this convex envelope in the (lifted) space of $\z$ variables which reside in a simplotope shown in Figure~\ref{fig:grid-rep}(c).

To describe our formulations, we will use the following prevalent concepts and notations. For a set $\x \in S \subseteq \R^{n}$, let $Q$ be a polyhedron with additional $r+q$ continuous variables defined as follows:
\[
Q = \bigl\{(\x,\y,\z) \in \R^n \times \R^r \times \R^q \bigm| A\x + B\y + C\z \leq d \bigr\},
\]
where $A \in \R^{m \times n}$, $B \in \R^{m \times r}$, $C \in \R^{m \times q}$ and $d \in \R^m$.
Now, consider a binary polyhedral set $E$ given by imposing binary conditions on some of the continuous variables in $Q$, where $E = \bigl\{(\x,\y,\z) \in Q \bigm| \z \in \{0,1\}^q \bigr\}$.
We say that $E$ is a mixed-integer linear programming (MIP) formulation of $S$ if the projection of the constraint set of $E$ onto the space of $\x$ variables is $S$. One of the most important properties of an MIP formulation is the strength of its natural linear programming (LP) relaxation obtained by ignoring the integrality condition since a tighter relaxation often indicates a faster convergence of the branch-and-bound algorithm, on which most  MIP solvers are built~\citep{vielma2015mixed}. We say an MIP formulation is \textit{ideal} if it has the strongest LP relaxation, that is, for any vertex $(\x,\y,\z)$ of its LP relaxation $Q$, we have $\z \in \{0,1\}^q$.

\subsection{Cross-item effects: tighter formulations via simplotope}\label{section:example} 
In this subsection, we focus on modeling the cross-item effects in price promotion planning problem. The insights developed here will be useful to introduce a more general class of functions later in Section~\ref{section:composite-MIP}. Under the linear demand function, the cross-item effect on revenue between two items, \textit{e.g.}, item $1$ and $2$, can be modeled as $\beta x_1x_2$, where $\beta$ is non-zero and $x_i$ is the price decision variable for item $i$. If $\beta$ is negative then two items are complementary; otherwise, they are substitutable. To capture these effects, we need to model the following discrete bilinear set:
\begin{equation}\label{eq:example}
	\Bigl\{(\x, \mu) \in \R^2 \times \R \Bigm| \mu = x_1x_2,\ x_i \in \{p_{i0}, \ldots, p_{id_i}\}\; \for \; i = 1,2 \Bigr\}. \tag{\textsc{Bi-Special}}
\end{equation}
In this subsection, we limit our attention to $d_2 = 1$. That is, the retailer decides whether the second product is promoted or not. \citet{cohen2021promotion} provide an MIP formulation for \eqref{eq:example} using McCormick inequalities~\citep{mccormick1976computability}, that is not ideal. Our goal is to show how to derive an ideal formulation.

The formulation of~\cite{cohen2021promotion} models the set~\eqref{eq:example} using the full binarization scheme given in~\eqref{eq:fullbinarization}. 
Then, the cross-item effect $x_1x_2$ is $\sum_{j,k}(p_{1j}p_{2k})\lambda_{1j}\lambda_{2k}$. \cite{cohen2021promotion} linearize $\lambda_{1j}\lambda_{2k}$ by introducing a variable $\gamma_{jk}$ and using McCormick inequalities, and obtain:
\begin{equation}\label{eq:cohen}
\begin{aligned}
\biggl\{(\x,\mu,\blambda,  \bgamma )& \biggm|  \eqref{eq:fullbinarization},\  \mu = \sum_{j = 0}^{d_1} \sum_{k = 0}^{d_2} (p_{1j}p_{2k}) \cdot \gamma_{jk},\  \gamma_{jk} \geq 0,\ \gamma_{jk} \geq \lambda_{1j} + \lambda_{2k}-1, \\
     & \gamma_{jk} \leq \lambda_{1j},\ \gamma_{jk}\leq \lambda_{2k} \, \for j \in 0 \cup [d_1] \text{ and } k \in 0\cup [d_2] \biggr\}.
\end{aligned}\tag{\textsc{Mc-Box}}
\end{equation}
Geometrically, this formulation convexifies $\lambda_{1j}\lambda_{2k}$ over the standard box $[0,1]^2$, ignoring that the variable $\blambda = (\blambda_1, \blambda_2 )$ lies in a simplotope. This results in an \emph{non}-ideal formulation as illustrated next.  
\begin{example}\label{ex:mcweak}
Let $(x_1,x_2) \in \{1,2,4\}\times\{1,2\}$. The set of our interest consists of $6$ price settings in $\R^3$, written compactly as the following matrix,
\renewcommand\arraystretch{1.3}
\setlength\arraycolsep{4pt}
\begin{equation}\label{eq:ex-vertices}
	\begin{pmatrix}
		1 & 2 & 4 & 1 & 2 & 4 \\
		1 & 1 & 1 & 2 & 2 & 2 \\
		1 & 2 & 4 & 2 & 4 & 8
	\end{pmatrix},
\end{equation}
where each column lists the two prices and their product in $\R^3$. We consider the LP relaxation of formulation~\eqref{eq:cohen} obtained by relaxing the binary requirement, and project it into the space of the original variables $(\x,\mu)$\endnote{This projection is obtained using the polyhedral computation package \textsf{Polyhedra.jl}~\citep{legat2023polyhedral} in \textsf{Julia} to project the LP relaxation of formulation~\eqref{eq:cohen} onto the space of the original variables $(\x,\mu)$.} using Fourier-Motzkin elimination. The vertices of the projection, listed as a matrix, are:

\renewcommand\arraystretch{2}
\begin{equation}\label{eq:ex-projected}
\begin{pmatrix}
1 & 1 &  4 & 4 & 3 &         \dfrac{7}{3}  &  \dfrac{7}{3}   &   \dfrac{3}{2}    &  3 & \dfrac{7}{3} & \dfrac{7}{3} & \dfrac{3}{2} \\
1 & 2 & 1 & 2 & \dfrac{3}{2}    & \dfrac{5}{3}     &\dfrac{4}{3}  & \dfrac{3}{2}&  \dfrac{3}{2} & \dfrac{5}{3} & \dfrac{4}{3} & \dfrac{3}{2}\\
 1 & 2 & 4 & 8 & 0&  0  & 0& 0 & 9 &7  & 7 & \dfrac{9}{2}
\end{pmatrix}, 
\end{equation}
and this matrix differs from~\eqref{eq:ex-vertices}. For example, let $\gamma_{1i} = 0$, $\lambda_{11}=0$ and, $\gamma_{ij} = \lambda_{1,i} = \lambda_{2,j} = 0.5$ for $i=2,3$ and $j=1,2$. This results in $x_1=3$, $x_2=1.5$, and $\mu=9$, which satisfies the constraints in \eqref{eq:cohen} and corresponds to the ninth column above. The third row is, however, not an element-wise product of the first two rows signifying that $\mu\ne x_1x_2$ even at extreme points. This means that even without business rules and with substitutable products, the LP relaxation of \eqref{eq:cohen} strictly overestimates the optimal revenue, showing that \eqref{eq:cohen} is not an ideal formulation. \hfill \Halmos

\end{example}

To make the formulation ideal, we will exploit the simplotope structure implicit in the binarization scheme. Using the unary binarization scheme for pricing variables, we obtain
\begin{equation}\label{eq:unary-example}
	x_i = p_{i0} + \sum_{j \in [d_i]} (p_{ij} - p_{ij-1})z_{ij}=:l_i(\z_i) \quad \text{ and } \quad \z_i \in \Delta^{d_i} \cap \{0,1\}^{d_i} \quad \for i = 1,2.
\end{equation}
The cross-item effect, now, is represented as $l_1(\z_1)l_2(\z_2)$. 
Treating each cross-item effect as a function over the vertex set of the simplotope $\Delta^{d_1} \times \Delta^{d_2}$, we develop the following affine under- and over-estimators:
\begin{subequations}\label{eq:generalziedMC}
\begin{alignat}{3}
  \mu &\leq p_{10}p_{20} + p_{10}(p_{21}-p_{20})  z_{21} + \sum_{j \in [d_1]}p_{21}(p_{1j}-p_{1(j-1)})  z_{1j} \label{eq:generalziedMC1}\\
\mu &\leq p_{10}p_{20}	+ \sum_{j \in [d_1]} p_{20}(p_{1j} - p_{1(j-1)})  z_{1j} + p_{1d_1}(p_{21}-p_{20}) z_{21}\label{eq:generalziedMC2} \\
\mu &\geq p_{10}p_{21} + p_{10}(p_{20}-p_{21})(1- z_{21}) + \sum_{j \in [d_1]}p_{20}(p_{1j} - p_{1(j-1)}) z_{1j} \label{eq:generalziedMC3}  \\
 \mu & \geq p_{10}p_{21} + \sum_{j \in [d_1]}p_{21}(p_{1j} - p_{1(j-1)}) z_{1j} + p_{1d_1}(p_{21}-p_{20})(1- z_{21}). \label{eq:generalziedMC4} 
\end{alignat}
\end{subequations}
In the setting of~\eqref{eq:example}, our inequalities generalize the following McCormick inequalities~\citep{mccormick1976computability} of $x_1x_2$ over $[p_{10},p_{11}] \times [p_{20},p_{21}]$:
\[
\begin{aligned}
	\mu &\leq \min \{p_{21}x_1 + p_{10}x_2 - p_{10}p_{21}, p_{20}x_1 + p_{11}x_2 - p_{11}p_{20} \} \\
	\mu &\geq \max\{ p_{20}x_1 + p_{10}x_2 -p_{10}p_{20}, p_{21}x_1 + p_{11}x_2 - p_{11}p_{21}\}.
\end{aligned}
\]
To see this, we take $d_1 = 1$ and additionally use the relation $x_i = p_{i0} + (p_{i1} - p_{i0})z_{i1}$ for $i = 1,2$.

The inequalities in~\eqref{eq:generalziedMC} arise from simplicial subsets of the underlying simplotope. Specifically, they are derived by affinely interpolating the cross-item effect $l_1(\cdot)l_2(\cdot)$ over the vertices of these simplicial subsets. This geometric intuition will be useful to treat more general cases in Section~\ref{section:supermodular}. For concreteness, we illustrate the geometry of~\eqref{eq:generalziedMC1} and~\eqref{eq:generalziedMC2} on a three dimensional example. Inequalities~\eqref{eq:generalziedMC3} and~\eqref{eq:generalziedMC4} can be derived similarly with a different simplicial partition. 
\begin{example}\label{ex:fourineq}
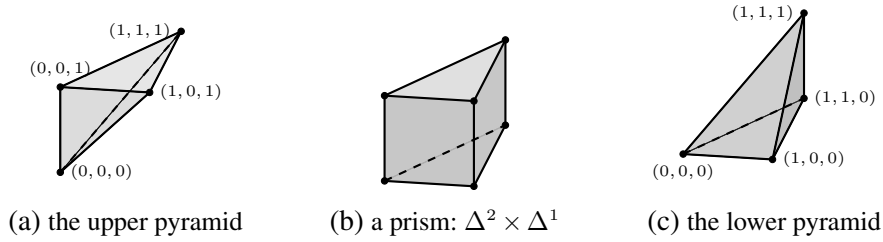
\begin{figure}[h]
	\centering
\begin{subfigure}[b]{0.3\linewidth}
\centering
  \tdplotsetmaincoords{70}{370}
\tikzsetnextfilename{promotion-4}
\begin{tikzpicture}[scale=1.2, tdplot_main_coords]
        \definecolor{point_color}{rgb}{ 0,0,0 }
    \tikzstyle{point_style} = [fill=point_color]

    \coordinate (q00) at (0, 0, 0);
    \coordinate (q01) at (0, 0, 1);
    \coordinate (q11) at (1, 0, 1);
    \coordinate (q21) at (1, 2, 1);

\node [right] at (q00) {\tiny{$(0,0,0)$}};
\node [above] at (q01) {\tiny{$(0,0,1)$}};
\node [left] at (q21) {\tiny{$(1,1,1)$}};
\node [right] at (q11) {\tiny{$(1,0,1)$}};

        \definecolor{edge_color}{RGB}{100,100,100}

      \definecolor{plane_color}{RGB}{249, 6, 33}
    \definecolor{shadow_color}{RGB}{0,64,152}

    \tikzstyle{plane_style} = [fill=plane_color, fill opacity=0.8, draw=plane_color, line width=0.1 pt, line cap=round, line join=round] 
    \tikzstyle{shadow_style} = [fill=gray, fill opacity=0.1, line width=0.005 pt, line cap=round, line join=round]

 	  \draw[shadow_style] (q00) -- (q01) -- (q21) -- (q00)-- cycle;
	  
	  \draw[shadow_style] (q00) -- (q01) -- (q11) -- (q00)-- cycle;
	       \draw[shadow_style] (q00) -- (q11) -- (q21) -- (q00) -- cycle;
	  \draw[shadow_style] (q01) -- (q11) -- (q21) -- (q01)-- cycle;
    	  \draw[shadow_style] (q00) -- (q01) -- (q11) -- (q00)-- cycle;

\draw[black, thick] (q00) -- (q01);
\draw[black, thick] (q11) -- (q01);
\draw[black, thick] (q11) -- (q21);
\draw[black, thick] (q01) -- (q21);
\draw[black, thick] (q00) -- (q11);
\draw[black, thick,dashed] (q00) -- (q21);
    \filldraw[black] (q00) circle (1pt) ; 
\filldraw[black] (q01) circle (1pt) ; 
\filldraw[black] (q11) circle (1pt); 
\filldraw[black] (q21) circle (1pt) ; 
  \end{tikzpicture}
    \caption{\small the upper pyramid}\label{fig:staircase-1}
  \end{subfigure}
  \begin{subfigure}[b]{0.2\linewidth}
  \centering
    \tdplotsetmaincoords{70}{370}
  \tikzsetnextfilename{promotion-5}
\begin{tikzpicture}[scale=1.2, color={lightgray},tdplot_main_coords]
    
    \coordinate (q00) at (0, 0, 0);
    \coordinate (q10) at (1, 0, 0);
    \coordinate (q20) at (1, 2, 0);
    \coordinate (q01) at (0, 0, 1);
    \coordinate (q11) at (1, 0, 1);
    \coordinate (q21) at (1, 2, 1);

        \definecolor{edge_color}{RGB}{100,100,100}

      \definecolor{plane_color}{RGB}{249, 6, 33}
    \definecolor{shadow_color}{RGB}{100,100,100}

    \tikzstyle{plane_style} = [fill=plane_color, fill opacity=0.8, draw=plane_color, line width=0.1 pt, line cap=round, line join=round] 
    \tikzstyle{shadow_style} = [fill=shadow_color, fill opacity=0.2, line width=0.05 pt, line cap=round, line join=round]

      \draw[shadow_style] (q00) -- (q10) -- (q20) -- (q00) -- cycle;
	  \draw[shadow_style] (q00) -- (q01) -- (q21) -- (q20) -- (q00)-- cycle;
	  \draw[shadow_style] (q00) -- (q01) -- (q11) -- (q10) -- (q00)-- cycle;
	  \draw[shadow_style] (q10) -- (q11) -- (q21) -- (q20) -- (q10)-- cycle;
    
    \draw[black,  thick,dashed] (q00) -- (q20);
\draw[black, thick] (q00) -- (q10);
\draw[black, thick] (q20) -- (q10);
\draw[black,  thick] (q00) -- (q01);
\draw[black,  thick] (q11) -- (q01);
\draw[black,  thick] (q11) -- (q21);
\draw[black,  thick] (q11) -- (q10);
\draw[black,  thick] (q20) -- (q21);
\draw[black,  thick] (q01) -- (q21);

\filldraw[black] (q00) circle (1pt) ; 
\filldraw[black] (q10) circle (1pt); 
\filldraw[black] (q20) circle (1pt) ; 
\filldraw[black] (q01) circle (1pt) ; 
\filldraw[black] (q11) circle (1pt); 
\filldraw[black] (q21) circle (1pt) ; 

  \end{tikzpicture}
      \caption{\small a prism: $\Delta^2 \times \Delta^1$}\label{fig:staircase-2}
  \end{subfigure}
  \begin{subfigure}[b]{0.3\linewidth}
  \centering
    \tdplotsetmaincoords{70}{370}
   \tikzsetnextfilename{promotion-6}
\begin{tikzpicture}[scale=1.2, tdplot_main_coords]
        \definecolor{point_color}{rgb}{ 0,0,0 }
    \tikzstyle{point_style} = [fill=point_color]

    \coordinate (q00) at (0, 0, 0);
    \coordinate (q10) at (1, 0, 0);
    \coordinate (q20) at (1, 2, 0);
    \coordinate (q21) at (1, 2, 1);

\node [left] at (q21) {\tiny{$(1,1,1)$}};
\node [right] at (q20) {\tiny{$(1,1,0)$}};
\node [right] at (q10) {\tiny{$(1,0,0)$}};
\node [below] at (q00) {\tiny{$(0,0,0)$}};

        \definecolor{edge_color}{RGB}{100,100,100}

      \definecolor{plane_color}{RGB}{249, 6, 33}
    \definecolor{shadow_color}{RGB}{100,100,100}

    \tikzstyle{plane_style} = [fill=plane_color, fill opacity=0.8, draw=plane_color, line width=1 pt, line cap=round, line join=round] 
    \tikzstyle{shadow_style} = [fill=gray, fill opacity=0.2, line width=0.05 pt, line cap=round, line join=round]
      \draw[shadow_style] (q00) -- (q10) -- (q20) -- (q00) -- cycle;
	  \draw[shadow_style] (q00) -- (q21) -- (q20) -- (q00)-- cycle;
	  \draw[shadow_style] (q00) -- (q10) --  (q21) -- (q00)-- cycle;
	  \draw[shadow_style] (q10) -- (q21) -- (q20) -- (q10)-- cycle;
	  \draw[shadow_style] (q10) -- (q21) -- (q20) -- (q10)-- cycle;

\draw[black,  thick,dashed] (q00) -- (q20);
\draw[black,  thick] (q00) -- (q10);
\draw[black,  thick] (q20) -- (q10);
\draw[black,  thick] (q10) -- (q21);
\draw[black,  thick] (q00) -- (q21);
\draw[black,  thick] (q20) -- (q21);

\filldraw[black] (q00) circle (1pt) ; 
\filldraw[black] (q10) circle (1pt); 
\filldraw[black] (q20) circle (1pt) ; 
\filldraw[black] (q21) circle (1pt) ; 
  \end{tikzpicture}
\caption{\small the lower pyramid}\label{fig:staircase-3}
\end{subfigure}
\caption{The geometry behind the  inequalities~\eqref{eq:generalziedMC1} and~\eqref{eq:generalziedMC2}.} \label{fig:staircase-ineq}
\end{figure}
Consider~\eqref{eq:generalziedMC1} and~\eqref{eq:generalziedMC2} when $(x_1,x_2) \in \{1,2,4\}\times\{1,2\}$. In this case, the cross-item effect is represented as $\psi(\z):=(1+ z_{11} + 2z_{12})(1+z_{21})$ over the vertices of the prism $\Delta:=\Delta^2 \times \Delta^1$, which is depicted in Figure~\ref{fig:staircase-2}.
To obtain the two inequalities, we will interpolate $\psi(\cdot)$ over the vertices of two pyramids, respectively. The first, depicted in Figure~\ref{fig:staircase-1}, has vertices: $(0,0,0)$, $ (0,0,1)$, $(1,0,1)$ and $(1,1,1)$, 
and relates to the inequality~\eqref{eq:generalziedMC1}. The second, depicted in Figure~\ref{fig:staircase-3}, has vertices: $(0,0,0)$, $(1,0,0)$, $(1,1,0)$ and $(1,1,1)$,
and relates to the inequality~\eqref{eq:generalziedMC2}. These special pyramids are such that the first one includes all elevated vertices of the prism (where $z_{21}$ is $1$) along with the origin and the second one includes all base vertices of the prism (where $z_{21}$ is $0$) along with the vertex $(1,1,1)$. A key property that will be exploited to show that \eqref{eq:generalziedMC} models the cross-item effect exactly is that every vertex of the prism is in at least one of the two pyramids. By expressing incremental changes of $\psi(\cdot)$ along the vertices of the first pyramid, we obtain 
\[
\psi(0,0,0) + (\psi(0,0,1)-\psi(0,0,0))z_{21}  + (\psi(1,0,1)-\psi(0,0,1))z_{11} + (\psi(1,1,1)-\psi(1,0,1))z_{12},
\]
which is $ z_{21} + 2z_{11} + 4z_{12} + 1$, exactly the right-hand-sides of \eqref{eq:generalziedMC1}. Similarly, interpolating $\psi(\cdot)$ over the vertices of the second pyramid, we obtain $z_{11}+2z_{12} + 4z_{21} + 1$, which is the right-hand-sides of~\eqref{eq:generalziedMC2}. \hfill \Halmos
\end{example}
\begin{proposition}\label{prop:simple-exact}
If $d_2 = 1$ then an MIP formulation of~\eqref{eq:example} is given by $\bigl\{(\x,\mu,\z) \in  \R^2 \times \R \times \{0,1\}^{d_1+d_2} \bigm|\eqref{eq:unary-example}\text{ and }\eqref{eq:generalziedMC} \bigr\}$.
\end{proposition}
We compare the size of our formulation and~\eqref{eq:cohen}. Our formulation introduces $d_1 + d_2$ binary variables to binarize discrete pricing variables, while~\eqref{eq:cohen} requires $(d_1+1)+(d_2 + 1)$  binary variables for binarization and $(d_1+1)\times (d_2 + 1)$ additional continuous variables to denote the marginal contribution to revenue at every promotion levels. In terms of the number of constraints, in addition to binarization, our formulation requires only four linear constraints; in contrast, formulation~\eqref{eq:cohen} uses $4 \times (d_1+1)\times (d_2 + 1)$ linear inequalities to linear  $\gamma_{jk}=\lambda_{1j}\lambda_{2k}$ where $j \in 0\cup [d_1]$ and $k \in 0 \cup [d_2]$. 

Although the formulation in Proposition~\ref{prop:simple-exact} is compact, it is not ideal either. However, the geometry of Example~\ref{ex:fourineq} allows us to derive an ideal formulation for \eqref{eq:example} as follows. 

\begin{figure}[h]
	\centering
    \tdplotsetmaincoords{70}{370}
  \tikzsetnextfilename{promotion-7}
\begin{tikzpicture}[scale=1.1, color={lightgray},tdplot_main_coords]

    \coordinate (q00) at (0, 0, 0);
    \coordinate (q10) at (1, 0, 0);
    \coordinate (q20) at (1, 2, 0);
    \coordinate (q01) at (0, 0, 1);
    \coordinate (q11) at (1, 0, 1);
    \coordinate (q21) at (1, 2, 1);

        \definecolor{edge_color}{RGB}{100,100,100}

      \definecolor{plane_color}{RGB}{249, 6, 33}
    \definecolor{shadow_color}{RGB}{100,100,100}

    \tikzstyle{plane_style} = [fill=plane_color, fill opacity=0.8, draw=plane_color, line width=1 pt, line cap=round, line join=round] 
    \tikzstyle{shadow_style} = [fill=shadow_color, fill opacity=0.2, line width=0.05 pt, line cap=round, line join=round]

      \draw[shadow_style] (q00) -- (q10) -- (q20) -- (q00) -- cycle;
	  \draw[shadow_style] (q00) -- (q01) -- (q21) -- (q20) -- (q00)-- cycle;
	  \draw[shadow_style] (q00) -- (q01) -- (q11) -- (q10) -- (q00)-- cycle;
	  \draw[shadow_style] (q10) -- (q11) -- (q21) -- (q20) -- (q10)-- cycle;

\draw[black, thick,dashed] (q00) -- (q20);
\draw[black, thick] (q00) -- (q10);
\draw[black, thick] (q20) -- (q10);
\draw[black, thick] (q00) -- (q01);
\draw[black, thick] (q11) -- (q01);
\draw[black, thick] (q11) -- (q21);
\draw[black, thick] (q11) -- (q10);
\draw[black, thick] (q20) -- (q21);
\draw[black, thick] (q01) -- (q21);

\filldraw[black] (q00) circle (1pt) ; 
\filldraw[black] (q10) circle (1pt); 
\filldraw[black] (q20) circle (1pt) ; 
\filldraw[black] (q01) circle (1pt) ; 
\filldraw[black] (q11) circle (1pt); 
\filldraw[black] (q21) circle (1pt) ; 
  \end{tikzpicture}
\qquad \quad \qquad
   \tikzsetnextfilename{promotion-8}
\begin{tikzpicture}[scale=1, color={lightgray},tdplot_main_coords]
        \definecolor{point_color}{rgb}{ 0,0,0 }
    \tikzstyle{point_style} = [fill=point_color]

    \coordinate (q00) at (0, 0, 0);
    \coordinate (q10) at (1, 0, 0);
    \coordinate (q20) at (1, 2, 0);
    \coordinate (q21) at (1, 2, 1);

        \definecolor{edge_color}{RGB}{100,100,100}

      \definecolor{plane_color}{RGB}{249, 6, 33}
    \definecolor{shadow_color}{RGB}{100,100,100}
    \tikzstyle{shadow_style} = [fill=gray, fill opacity=0.2, line width=0.05 pt, line cap=round, line join=round]

    \tikzstyle{plane_style} = [fill=plane_color, fill opacity=0.8, draw=plane_color, line width=1 pt, line cap=round, line join=round] 
\draw[black,dashed] (q00) -- (q20);
\draw[black] (q00) -- (q10);
\draw[black] (q20) -- (q10);
\draw[black] (q10) -- (q21);
\draw[black] (q00) -- (q21);
\draw[black] (q20) -- (q21);

      \draw[shadow_style] (q00) -- (q20) -- (q21) -- (q00) -- cycle;
     \draw[shadow_style] (q00) -- (q10) -- (q20) -- (q00) -- cycle;
     \draw[shadow_style] (q10) -- (q20) -- (q21) -- (q10) -- cycle;
     \draw[shadow_style] (q00) -- (q10) -- (q21) -- (q00) -- cycle;
          \draw[shadow_style] (q10) -- (q20) -- (q21) -- (q10) -- cycle;

    \tikzstyle{shadow_style1} = [fill=blue, fill opacity=0.2, line width=0.05 pt, line cap=round, line join=round]

\coordinate (q00mid) at (0, 0, 0.4);
\coordinate (q10mid) at (1, 0, 0.4);
\coordinate (q11mid) at (1, 0, 1.4);
\coordinate (q21mid) at (1, 2, 1.4);

\filldraw[black] (q00mid) circle (1pt) ; 
\filldraw[black] (q10mid) circle (1pt); 
\filldraw[black] (q11mid) circle (1pt) ; 
\filldraw[black] (q21mid) circle (1pt) ; 
\draw[shadow_style1] (q00mid) -- (q10mid) -- (q21mid) -- (q00mid) -- cycle;
\draw[shadow_style1]  (q10mid) -- (q21mid) -- (q11mid) -- (q10mid) -- cycle;
\draw[shadow_style1]  (q00mid) -- (q21mid) -- (q11mid) -- (q00mid) -- cycle;

\draw[shadow_style1] (q00mid) -- (q10mid) -- (q11mid) -- (q00mid) -- cycle;
\draw[shadow_style1] (q00mid) -- (q10mid) -- (q11mid) -- (q00mid) -- cycle;
\draw[black,ultra thick] (q00mid) -- (q11mid);
\draw[black,ultra thick] (q00mid) -- (q10mid);
\draw[black,ultra thick] (q00mid) -- (q10mid);
\draw[black,ultra thick] (q11mid) -- (q10mid);
\draw[black,ultra thick] (q11mid) -- (q21mid);
\draw[black,ultra thick] (q10mid) -- (q21mid);
\draw[black,ultra thick] (q00mid) -- (q21mid);

  \coordinate (q00up) at (0, 0, 0.8);
    \coordinate (q01up) at (0, 0, 1.8);
    \coordinate (q11up) at (1, 0, 1.8);
    \coordinate (q21up) at (1, 2, 1.8);

                    \draw[black, thick] (q00) -- (q10);
\draw[black, thick] (q10) -- (q21);
\draw[black, thick] (q00) -- (q21);
\draw[black, thick] (q20) -- (q10);
\draw[black, thick] (q20) -- (q21);
    \draw[shadow_style] (q00up) -- (q01up) -- (q21up) -- (q00up) -- cycle;
    \draw[shadow_style] (q00up) -- (q21up) -- (q11up) -- (q00up) -- cycle;
     \draw[shadow_style] (q00up) -- (q01up) -- (q11up) -- (q00up) -- cycle;
     \draw[shadow_style] (q01up) -- (q21up) -- (q11up) -- (q01up) -- cycle;
     \draw[shadow_style] (q00up) -- (q01up) -- (q11up) -- (q00up) -- cycle;

\draw[black,thick] (q00up) -- (q01up);
\draw[black,thick] (q11up) -- (q01up);
\draw[black,thick] (q11up) -- (q21up);
\draw[black,thick] (q01up) -- (q21up);
\draw[black,thick] (q00up) -- (q11up);
\draw[black,dashed] (q00up) -- (q21up);

  \end{tikzpicture}
\caption{Subdivision of prism into three pyramids, where the blue one does not appear in Figure~\ref{fig:staircase-ineq}.} \label{fig:triangulation-basic}
\end{figure}
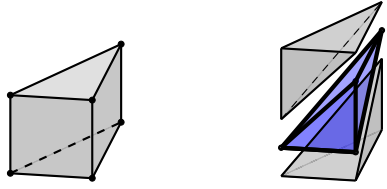
\begin{example}\label{ex:bilinearcontinues}
Consider the case $(x_1,x_2) \in \{1,2,4\}\times\{1,2\}$ again, and recall that the cross-item effect is represented as $\psi(\z):=(1+z_{11} + 2z_{12})(1+z_{21})$ over the vertices of the prism, where, in Example~\ref{ex:fourineq}, two gray  pyramids are associated with two inequalities. Clearly, these pyramids do not cover the entire prism. The region that is left is contained in the blue pyramid, whose vertices are: $(0,0,0)$, $(1,0,0)$, $(1,0,1)$ and  $(1,1,1)$.
Expressing incremental changes of $\psi(\cdot)$ along the vertices, we obtain an overestimator of $\psi(\cdot)$:
\[
\psi(0,0,0) + (\psi(1,0,0)-\psi(0,0,0) )z_{11} + (\psi(1,0,1)-\psi(1,0,0) )z_{21} + (\psi(1,1,1)-\psi(1,0,1) )z_{12},
\]
which yields  $ \mu \leq 1+z_{11}+2z_{21}+4z_{12}$, tightening the formulation given by inequalities in~\eqref{eq:generalziedMC}. \hfill  \Halmos
\end{example}
The two pyramids used to derive inequalities in  Example~\ref{ex:fourineq} left the interior of a pyramid uncovered (see Figure~\ref{fig:triangulation-basic}). When $d_1 > 2$, this region decomposes into $d_1-1$ simplices, each of which generates a valid inequality. This approach yields linear inequalities for the cross-item effect $l_1(\z_1) l_2(\z_2)$, each affinely interpolating it over a simplical subset of the simplotope $\Delta^{d_1} \times \Delta^{d_2}$:
\begin{subequations}\label{eq:hull-example}
	\begin{alignat}{3}
		\mu \leq p_{10}p_{20} &+  \sum_{j=1}^k p_{20}(p_{1j}-p_{1(j-1)})z_{1j} + p_{1k}(p_{21} - p_{20}) z_{21} \notag \\
& + \sum_{j = k+1}^{d_1}p_{21}(p_{1j}-p_{1(j-1)})z_{1j} \qquad \for k \in 0 \cup [d_1] \label{eq:hull-example-1} \\
\mu \geq p_{10}p_{21}  &+ \sum_{j=1}^kp_{21}(p_{1j}-p_{1(j-1)})z_{1j} + p_{1k}(p_{20}-p_{21})(1-z_{21})  \notag \\
& + \sum_{j = k+1}^{d_1}p_{20}(p_{1j}-p_{1(j-1)})z_{1j} \qquad \for k \in 0 \cup [d_1]. \label{eq:hull-example-2}
	\end{alignat}
\end{subequations}
\begin{proposition}\label{prop:simple-perfect}
If $d_2 = 1$ then an ideal MIP formulation of~\eqref{eq:example} is given by  $\bigl\{(\x,\mu,\z) \in  \R^2 \times \R \times \{0,1\}^{d_1+d_2} \bigm|\eqref{eq:unary-example}\text{ and }\eqref{eq:hull-example}\bigr\}$.
\end{proposition}
\section{Modeling Composite Effects}\label{section:composite-MIP}
In Section~\ref{section:example}, we motivated our study by examining a special case of the additive cross-item effect. Building on those
observations, we develop strong formulations for the price promotion problem under general composite demand functions. The key insight, formalized in Section~\ref{section:formalization}, is that the ideality in Proposition~\ref{prop:simple-perfect} arises from the
\textit{convex hull} of the bilinear term described by the linear inequalities~\eqref{eq:hull-example} over the corresponding simplotope. Subsequently, in Section~\ref{section:supermodular}, we develop explicit formulations for \textit{supermodular} compositions, capturing a wide range of practically relevant nonlinear effects, including additive and multiplicative interactions as detailed in Section~\ref{section:ppo}. Finally, in Section~\ref{section:otherbinarization}, we show that ideal formulations can be preserved under alternative binarization schemes. 
\subsection{From the unary simplotope to an ideal formulation}\label{section:formalization}

We first formalize the role of the unary simplotope in modeling
compositions of discrete univariate functions. Let $\X := \prod_{i\in[n]}
\{p_{i0},p_{i1},\ldots,p_{id_i}\}$ be a Cartesian product of finite sets of lattice points and $P$ be an $n$-dimensional polytope. Consider
the composite function
\[
(\phi\mcirc f)(\x)
=
\phi\bigl(f_1(x_1),\ldots,f_n(x_n)\bigr),
\qquad \x\in\X \cap P,
\]
where each $f_i$ is a univariate function and
$\phi:\R^n\to\R$ is an arbitrary function. We are interested in the
graph
\[
\graph (\phi \mcirc f) := \Bigl\{(\x,\mu) \in \R^n \times \R \Bigm|  \mu = \phi\bigl(f_1(x_1), \ldots, f_n(x_n) \bigr),\ \x \in P \cap \X  \Bigr\}.
\]

The key observation is that the discrete function can be represented
on the vertices of a unary simplotope. Recall that, for each $i\in[n]$,
$\Delta^{d_i}$ denotes the $d_i$-dimensional unary simplex, and
$\vertex(\Delta^{d_i})$ denotes its vertices, specified in~\eqref{eq:verticesDelta}. We lift $f_i$ to the affine function
\begin{equation}\label{eq:liftedfunction}
 \ell_i(\z_i)
 :=
 f_i(p_{i0})
 +
 \sum_{j\in[d_i]}
 \bigl(f_i(p_{ij})-f_i(p_{i,j-1})\bigr)z_{ij},
 \qquad
 \z_i\in\vertex(\Delta^{d_i}).
 \tag{\textsc{Lifting}}
\end{equation}
The lifted composite function is then defined on the vertices of the simplotope
$\Delta$ by
\[
(\phi\mcirc\ell)(\z)
:=
\phi\bigl(\ell_1(\z_1),\ldots,\ell_n(\z_n)\bigr),
\qquad
\z\in\vertex(\Delta).
\]
Consequently, every point in $\graph(\phi\mcirc f)$ has a
corresponding point in
\[
\graph(\phi\mcirc\ell)
:=
\left\{
(\z,\mu)\in\vertex(\Delta)\times\R
\;\middle|\;
\mu=(\phi\mcirc\ell)(\z)
\right\}.
\]
Our main observation is that replacing the graph of the lifted function with its convex hull yields an ideal formulation in the absence of polyhedral constraints.

\begin{lemma}\label{lemma:MIP-discrete}
Let $ \Phi$ be a polytope in $(\z,\mu)$ space such that $\Phi \cap \bigl(\{0,1\}^{\sum_{i\in [n]}d_i} \times \R \bigr)  = \graph(\phi \mcirc \ell)$.
Then, an MIP formulation of the graph of $\phi \mcirc f$ is given as follows
\begin{equation*}
\begin{aligned}
\Bigl\{(\x,\z,\mu)  \Bigm| (\z,\mu) \in  \Phi,\ \x \in P ,\ (x_i,\z_i)  \in  \eqref{eq:unary} \; \for i \in [n] \Bigr\}.
\end{aligned}
\end{equation*}
This formulation is ideal when $\Phi = \conv \bigl(\graph(\phi \mcirc \ell)\bigr)$, and the constraint $x\in P$ is relaxed. 
\end{lemma}
The key insight above is that each point in the graph of $\phi \mcirc f$ over $\X$ lifts to a corner point of the graph of $\phi \mcirc \ell$. This is accomplished using $\z$ variables ($\sum_{i\in [n]} d_i$ of them), each of which takes binary values.

\subsection{Supermodular composition}\label{section:supermodular}
{\revision
We now specialize the general framework of Lemma~\ref{lemma:MIP-discrete} to composite functions with additional structure. Our motivation comes from the multiplicative interactions that arise in the price promotion problem. In particular, under the log-log demand model given in~\eqref{eq:multdemand}, the profit contribution associated with the self-effect and historical effects of item $i$ in period $t$ is the product of several nonnegative univariate functions, as detailed in~\eqref{eq:multobj}. This structure motivates the class of \textit{supermodular compositions} studied in this section.
}

\begin{definition}[\cite{topkis2011supermodularity}]\label{defn:supermodularity}
A function $f(\x): X \subseteq \R^n \to \R$ is \emph{supermodular} if
\[
f(\x' \vee \x'') + f(\x' \wedge \x'') \geq f(\x') + f(\x'')
\]
for all $\x', \x'' \in X$, where $\x' \vee \x''$ (resp. $\x' \wedge \x''$) denotes the componentwise maximum (resp. minimum). A function is \emph{submodular} if its negation is supermodular. \hfill \Halmos
\end{definition}
{\revision
The product function $\phi(\y)=\prod_{j\in[n]} y_j$ is supermodular on the nonnegative orthant $\{\y\in\R^n\mid \y\geq 0\}$. Consequently, the profit term discussed above can be viewed as a supermodular composition in which the product function serves as the outer function and the self-effect and historical-effect functions serve as the inner univariate functions. More generally, we study compositions of a supermodular (resp. submodular) function with univariate discrete functions, which we refer to as \emph{supermodular (resp. submodular) compositions}. Importantly, the resulting composition need not be supermodular in the original variables when the inner univariate functions are non-monotone. The ordering construction developed below addresses this issue and allows us to obtain convex hull formulations for these compositions. This perspective provides a unified framework for modeling the historical and cross-item effects arising in price promotion planning.} Beyond promotion planning, supermodular compositions contain the class of $\lnatural$-convex functions, a powerful tool for analyzing economic or operations models with indivisibilities~\citep{murota1998discrete,chen2021discrete}. The connection to discrete convex analysis is discussed further in Section~\ref{section:discreteconvexity}.

Because the composite effects enter the PPO objective in a maximization problem, it suffices to treat the following hypograph 
\begin{equation*}
    \hypo(\phi \mcirc f) := \Bigl\{(\x,\mu) \in \X\times \R \Bigm| \mu \le \phi\bigl(f_1(x_1),\ldots,f_n(x_n)\bigr)\Bigr\}.
\end{equation*}
To model this set, we leverage Lemma~\ref{lemma:MIP-discrete}. More specifically, we lift the lattice points in $\X$ to  the vertices of $\Delta$ and represent each inner-function $f_i(\cdot)$ as an affine function $\ell_i(\cdot)$ in the lifted space using~\eqref{eq:liftedfunction}, obtaining the following extended formulation of the hypograph
\begin{equation*}
    \Bigl\{(\x,\z,\mu)  \Bigm| \mu \le (\phi\mcirc \ell)(\z),\ (x_i,z_i)\in \eqref{eq:unary} \, \for i \in [n] \Bigr\}.
\end{equation*}
 To invoke Lemma~\ref{lemma:MIP-discrete}, we need a convex hull description of the hypograph of the lifted function $\phi \mcirc \ell$, that is
\begin{equation*}\label{eq:liftedhypo}
	\hypo(\phi \mcirc \ell):= \bigl\{(\z, \mu) \in \vertex(\Delta) \times \R \bigm| \mu \leq (\phi \mcirc \ell)(\z) \bigr\}.
\end{equation*}
{\revision Although $\phi(\cdot)$ is supermodular, the composition $(\phi \mcirc f)(\cdot)$, and consequently its lift $(\phi \mcirc \ell)(\cdot)$, need not be supermodular when the inner functions $f_i(\cdot)$ are non-monotone, as illustrated in the following  PPO example.

\begin{example}\label{ex:reorder}
Consider a single item over three periods, with discrete price levels
$x_t \in \{1,2,3\}$ for $t=1,2,3$. Under the log-log demand model,
as described in~\eqref{eq:multobj}, the profit contribution associated
with the self-effect in period 3 and two historical effects may take the
following multiplicative form:
\[
    (x_3-1)(x_3)^{-4}x_2^3x_1^2.
\]
The self-effect is $f_3(x_3):=(x_3-1)(x_3)^{-4}$, whose values at the three price levels are  $f_3(1)=0$, $f_3(2) = \frac{1}{16}$ and $f_3(3) = \frac{2}{81}$. Thus, although the self-effect is nonnegative, it is not monotone in
the price level.  In contrast, the two historical effects,
$f_2(x_2):=x_2^3$ and $f_1(x_1):=x_1^2$, are increasing in their respective price levels. Consequently,  the lifted profit function induced by the natural ordering of the price levels 
\[
\Bigl(\frac{1}{16}z_{31} - \frac{49}{1296}z_{32} \Bigr) (1+7z_{21} + 19z_{22})(1+ 3z_{11} + 5z_{12}),
\]
which contains negative coefficients when expanded. Hence, the natural lift is not supermodular.  \hfill\Halmos
\end{example}
}

Nevertheless, we can recover supermodularity in the lifted space by reordering the inner functions according to their function values.
To do so, we construct a permutation $\sigma = (\sigma_1, \ldots, \sigma_n)$ so that inner-function values are in non-decreasing order as follows:
\begin{equation}\label{eq:IO}
	f_i(p_{i,\sigma_i(0)}) \leq f_i(p_{i,\sigma_i(1)}) \leq \cdots \leq f_i(p_{i,\sigma_i(d_i)}) \qquad \for i \in [n].\tag{\textsc{Increasing-Ordered}}
\end{equation}
For convenience of notation, we define $f^{\sigma_i}_i(\cdot)$ such that $f_i^{\sigma_i}(p_{ij}) = f_i(p_{i, \sigma_i(j)})$ for $j \in \{0\} \cup [d_i]$ and observe that $f_i^{\sigma_i}(\cdot)$ is a non-decreasing function. 
Then, let:
\[
(\phi \mcirc f^\sigma)(\x') = \phi\bigl(f^{\sigma_1}_1(x'_1), f^{\sigma_2}_2(x'_2), \ldots, f^{\sigma_n}_n(x'_n) \bigr) \quad \for \, \x' \in \X,
\]
and observe that $(\phi \mcirc f^\sigma)(\cdot)$ is supermodular by Topkis' supermodularity-preservation result~\citep[Lemma 2.6.4]{topkis2011supermodularity} as long as $\phi(\cdot)$ is supermodular.
Our interest is, however, to derive a supermodular function in the lifted space. 
For each $i \in [n]$, we define $\ell_i^{\sigma_i}:\vertex(\Delta_i) \to \R$ so that $\ell_i^{\sigma_i}(\v_{ij}) = \ell_i(\v_{i,\sigma_i(j)})$, where we recall that $\v_{i0}, \v_{i1}, \cdots \v_{id_i}$ are vertices of $\Delta^{d_i}$ described in~\eqref{eq:verticesDelta}. It follows easily from \citep[Lemma 2.6.4]{topkis2011supermodularity} that if $\phi(\cdot)$ is supermodular then $(\phi\circ \ell^\sigma)(\cdot)$ defined as follows: 
\begin{equation*}\label{eq:switched}
(\phi \mcirc \ell^\sigma)(\z') = \phi\bigl(\ell^{\sigma_1}_1(\z'_1), \ell^{\sigma_2}_2(\z'_2), \ldots, \ell^{\sigma_n}_n(\z'_n) \bigr) \qquad \for \, \z'  \in \vertex(\Delta),
\end{equation*}
is also supermodular since $\ell_i^{\sigma_i}(\cdot)$ is a non-decreasing function.
{\revision
\begin{example}\label{ex:reorder-ordering}
Consider the setup in Example~\ref{ex:reorder}. Since
$f_1(x_1)=x_1^2$ and $f_2(x_2)=x_2^3$ are already non-decreasing in the price level, we have $\sigma_1 = \sigma_2 = (0,1,2)$. For the self-effect, $f_3$ maps $\{1,2,3\}$ to $\{0,\frac{1}{16},\frac{2}{81}\}$, and hence $\sigma_3=(0,2,1)$.  Consequently, the lifted profit function under the ordering
$(\sigma_1,\sigma_2,\sigma_3)$ is
\begin{equation}\label{eq:reorder-lift}
(\phi\mcirc\ell^\sigma)(\z')
={}
(1+3z'_{11}+5z'_{12})
(1+7z'_{21}+19z'_{22})
\left(
\frac{2}{81}z'_{31}
+\frac{49}{1296}z'_{32}
\right),
\end{equation}
which is supermodular.
\hfill\Halmos
\end{example}
}

We now focus on describing the convex hull  of the hypograph of $\phi \mcirc \ell^\sigma$ since it can be shown that an affine transformation, specified later in~\eqref{eq:variable-switching-general},  of such a description yields the convex hull of $\phi\mcirc \ell$.  Our description uses the grid representation of the lattice set $\X$ and $\vertex(\Delta)$ introduced in Section~\ref{section:example}. Recall that  with each point $\j:= (j_1, \ldots,j_n)$ in the grid $\G:=\prod_{i \in [n]}\{0, \ldots, d_i\}$, we associate a lattice point $(p_{1j_1}, p_{2j_2}, \ldots,  p_{nj_n})$ and a vertex $(\v_{1j_1}, \v_{2j_2}, \ldots, \v_{nj_n})$ of the unary simplotope, where $\v_{ik}$ is the $k^{\text{th}}$ vertex of the simplex $\Delta^{d_i}$ defined as in~\eqref{eq:verticesDelta}. The grid representation is convenient since we can use $\psi^\sigma(\j)$ to capture the function values at the points of interest:
\begin{equation}\label{eq:grid-function}
\psi^\sigma(\j) : = (\phi \mcirc f^\sigma)(p_{1j_1}, p_{2j_2}, \ldots, p_{nj_n} )  = (\phi \mcirc \ell^\sigma) \bigl(\v_{1j_1}, \v_{2 j_2}, \ldots, \v_{n j_n} \bigr) \qquad \for \, \j \in \G.\tag{\textsc{GridRep}}
\end{equation}

Our convex hull result is related to certain walks over the grid $\G$, where each step along the walk does not decrease the coordinate values, and each walk yields one of the inequalities that describes the convex hull. Because of the non-decreasing property, these walks are often referred to as staircases, which we next define formally.
A \textit{staircase} in the grid $\G$ is a sequence of points (or stops) $\j^0, \ldots, \j^r$ in $\G$ satisfying
\begin{itemize}
    \item $\j^0 = (0, \ldots, 0)$ and $\j^r = (d_1, \ldots, d_n)$, and
    \item for $s \in [r] $, $\j^s - \j^{s-1} = \e_{i}$ for some $i \in [n]$, where $\e_{i}$ is the principal vector in $i^{\text{th}}$ direction. 
\end{itemize}
In words, a staircase is a walk over the grid points going from $(0, \ldots, 0)$ to $(d_1, \ldots, d_n)$ with each stop differing from the previous stop in a single coordinate direction. 

{\revision
Alternatively, we can specify a staircase by the coordinate direction
of each move. Let \(D:=\sum_{i\in[n]}d_i\) denote the total number of
moves. For each move \(s\in[D]\), let
\[
\pi_s=(i_s,j_s) 
\qquad \text{ and} \quad
j_s:=\sum_{k=1}^s\ind\{i_k=i_s\}.
\]
Thus, \(i_s\) specifies the direction of the \(s\)-th move, while
\(j_s\) specifies the number of steps taken in that direction up to
and including move \(s\). We refer to
\(\pi=(\pi_1,\ldots,\pi_D)\) as the direction representation of the
staircase, and denote the set of all such representations by \(\Pi\).
}

Now, we use these representations to construct a family of linear inequalities. With a direction representation $\pi \in \Pi$, we associate a point representation $\j_\pi$, obtained by recursively using the following relations:
\begin{equation}\label{eq:pointrep}
    \j^0_\pi  = (0, \ldots, 0 ) \quad \text{ and } \quad \j^s_\pi = \j^{s-1}_\pi + \e_{\pi_s(1)} \quad \for s \in [D]. \tag{\textsc{StairPoint}}
\end{equation}
Evaluating the grid function $\psi^\sigma(\cdot)$ at these points, we obtain consecutive differences, $\psi^\sigma(\j^1_\pi) - \psi^\sigma(\j^{0}_\pi), \ldots, \psi^\sigma(\j^D_\pi) - \psi^\sigma(\j^{D-1}_\pi)$, each of which corresponds to a move. These incremental changes can be organized to yield the following inequality:
\begin{equation}\label{eq:staircase}
    \mu \leq		\psi^\sigma(\j^0_\pi) + \sum_{s  \in [D]} \bigl( \psi^\sigma(\j^s_\pi) - \psi^\sigma(\j^{s-1}_\pi) \bigr) \cdot z'_{\pi_s} \qquad \for \pi \in \Pi. \tag{\textsc{Staircase}}
\end{equation}
In the proof of Theorem~\ref{them:supermodular-model}, we will show that this class of inequalities describes the convex hull of the hypograph of a supermodular function over the vertices of the unary simplotope.

The remaining detail is the affine mapping that relates convex hulls of $\phi \mcirc \ell^\sigma$ and $\phi \mcirc \ell$. 
Let $P^{\sigma_i}$ be a permutation matrix in $\R^{(d+1)\times (d+1)}$ such that, for all $(j,k)$, $P^{\sigma_i}_{jk}=1$ if $j=\sigma_i(k)$ and zero otherwise. The affine transformation associated with $\sigma$ is $L^\sigma(\z) := \bigl(L^{\sigma_1}(\z_1), \ldots, L^{\sigma_n}(\z_n)\bigr)$, where 
\begin{equation}\label{eq:variable-switching-general}
L^{\sigma_i}(\z_i) := \Lambda^{-1}_i P^{\sigma_i} \Lambda_i\z_i. \tag{\textsc{Switch}}
\end{equation}
Here, the mapping $\Lambda_i$,  defined as in~\eqref{eq:fromztolambda}, maps a point $\z_i$ in the unary simplex to a point $\blambda_i$ in the standard simplex; this point is then mapped back to the unary simplex via the inverse of $\Lambda_i$. In the next example, we illustrate how the mapping relates $\ell$ with $\ell^\sigma$ in the context of price promotion and how an ideal MIP formulation is derived.
{\revision
\begin{example}\label{ex:cuts}
    Consider the setup in Examples~\ref{ex:reorder} and~\ref{ex:reorder-ordering} again. Recall that, to make the self-effect function $f_3(\cdot)$ non-decreasing, we use the order $\sigma_3 = (0,2,1)$. Under this order, the transformation~\eqref{eq:variable-switching-general} maps $\z_3$ to $\z_3'$ so that $ z'_{31} = z_{31}$ and $z'_{32} = z_{31} -z_{32}$.  Indeed, substituting the switch mapping into the ordered lifted function gives
    \[
            \ell_3^{\sigma_3}(\z'_3)
        =\frac{2}{81}z_{31}
        +\frac{49}{1296}(z_{31}-z_{32})
        =\frac{1}{16}z_{31}-\frac{49}{1296} z_{32}
        =\ell_3(\z_3).
    \]

  Next, we illustrate the derivation of one of the proposed inequalities. Consider a staircase $\pi$ specified by the direction
    sequence $(1,2,3,1,2,3)$. The corresponding sequence of grid points is
    \[
    \begin{aligned}
    \j_\pi^0&=(0,0,0)&
    \j_\pi^1&=(1,0,0)&
    \j_\pi^2&=(1,1,0),\\
    \j_\pi^3&=(1,1,1)&
    \j_\pi^4&=(2,1,1)&
    \j_\pi^5&=(2,2,1)&
    \j_\pi^6&=(2,2,2).
    \end{aligned}
    \]
    We then evaluate the grid representation of the reordered profit
    function~\eqref{eq:reorder-lift} at these staircase stops:
   \[
    \begin{aligned}
        \psi^\sigma(\j_\pi^0)&=0 &
        \psi^\sigma(\j_\pi^1)&=0 &
        \psi^\sigma(\j_\pi^2)&=0,\\
        \psi^\sigma(\j_\pi^3)&=\frac{64}{81}&
        \psi^\sigma(\j_\pi^4)&=\frac{16}{9} &
        \psi^\sigma(\j_\pi^5)&=6&
        \psi^\sigma(\j_\pi^6)&=\frac{243}{16}.
    \end{aligned}
    \]
 Therefore, the consecutive differences in profit along the staircase are
    \[
        0,\quad 0,\quad
        \frac{64}{81},\quad
        \frac{80}{81},\quad
        \frac{38}{9},\quad
        \frac{147}{16},
    \]
    where the first two differences are zero because the self-effect is zero at
    the lowest price level. Following~\eqref{eq:staircase}, and noting that
    the six moves correspond respectively to $z'_{11},z'_{21},z'_{31},z'_{12},z'_{22},z'_{32}$, we obtain the following staircase inequality in the space of $\z'$ variables:
    \[
        \mu
        \leq
        \frac{64}{81}z'_{31}
        +\frac{80}{81}z'_{12}
        +\frac{38}{9}z'_{22}
        +\frac{147}{16}z'_{32}.
    \]
    Finally, substituting the switch mapping gives the corresponding
    inequality in the original lifted variables:
    \begin{equation}\label{eq:stairIneqEx}
        \mu
        \leq
        \frac{64}{81}z_{31}
        +\frac{80}{81}z_{12}
        +\frac{38}{9}z_{22}
        +\frac{147}{16}(z_{31}-z_{32}).
    \end{equation}
    Thus, the staircase inequality is first constructed in the reordered space, where the lifted function is supermodular, and then mapped back to the original lifted variables through \(L^\sigma\). \hfill \Halmos
    \end{example}
}

\begin{theorem}\label{them:supermodular-model}
Let $\sigma:=(\sigma_1, \ldots, \sigma_n )$ satisfying~\eqref{eq:IO} and $L^{\sigma}$ be the affine mapping~\eqref{eq:variable-switching-general}.  
If $\phi(\cdot)$ is supermodular then an ideal MIP formulation of $ \hypo(\phi \mcirc f) $  is 
	\[
	\Bigl\{(\x,\z,\mu) \Bigm| (x_i,z_i) \text{ satisfies } \eqref{eq:unary}\, \for i\in [n],\ \z'= L^{\sigma}(\z),\  \eqref{eq:staircase}
	 \Bigr\}.
	\]
\end{theorem}
In the proof of Theorem~\ref{them:supermodular-model}, we demonstrate that the staircase inequalities characterize the convex hull of the hypograph $\phi \mcirc \ell^\sigma$. Their validity stems from a constructive procedure that exploits supermodularity, while their tightness is proved by triangulating the simplotope using staircases \citep[Section 6.2.3]{de2010triangulations}. We note that supermodularity has been used for conic mixed-binary sets \citep{atamturk2020submodularity,kilincc2025conic} and staircase triangulations have been used for convex relaxations over continuous domains \citep{he2022tractable}. However, these works do not consider discrete lattice domains and do not allow composition of a supermodular function with non-monotone univariate functions.

Now, we specialize Theorem~\ref{them:supermodular-model} to the case when the outer-function $\phi(\cdot)$ is a bilinear term, extending the MIP formulation for~\eqref{eq:example} in Propositions~\ref{prop:simple-exact} and~\ref{prop:simple-perfect} to relax the restriction on $d_2$. More specifically, we obtain an ideal formulation for the product of two arbitrary univariate functions:
\begin{equation}\label{eq:bifunc}
    	\Bigl\{(\x, \mu) \in \R^2 \times \R \Bigm| \mu = f_1(x_1)f_2(x_2),\ x_i \in  \{p_{i0}, \ldots, p_{id_i}\}\; \for \; i = 1,2 \Bigr\}. \tag{\textsc{Bi-General}}
\end{equation}
For the promotion planning problem, this will allow us to capture cross-item effects where each item has multiple price levels. Moreover, by recursively using this technique we will be able to model higher-order effects as well, as detailed in Section~\ref{section:loglog}. 

{\revision Observe that  $f_1(x_1)f_2(x_2)$ is a supermodular composition of two univariate functions   because the bilinear term, which is supermodular, can be used as the outer function, with $f_1(x_1)$ and $f_2(x_2)$ as the inner functions. The composite function does not need to be supermodular, since $f_i(x_i)$ can be an arbitrary non-monotone function.} Thus, we  use $\sigma = (\sigma_1, \sigma_2)$ that satisfies~\eqref{eq:IO} in Theorem~\ref{them:supermodular-model} to obtain an ideal formulation for the hypograph $\mu \leq f_1(x_1)f_2(x_2)$. To obtain an ideal formulation for the epigraph $\mu \geq f_1(x_1)f_2(x_2)$, we use a different tuple of permutations $\varsigma = (\varsigma_1,\varsigma_2)$ such that 
\begin{equation}\label{eq:mix-order}
    f_1(p_{1,\varsigma_1(0)}) \leq \cdots \leq f_1(p_{1,\varsigma_1(d_1)}) \text{ and } f_2(p_{2,\varsigma_2(0)}) \geq \cdots \geq f_2(p_{2,\varsigma_2(d_2)}) \tag{\textsc{Mix-Ordered}}.
\end{equation}
This is analogous to transforming a supermodular function  $f_1f_2$ over the box $[0,1]^2$ into the submodular function $f_1(1-f_2)$ by switching  one of the arguments~\citep{tawarmalani2013explicit}.

\begin{proposition}\label{prop:bilinear-sub}
Let $\sigma = (\sigma_1, \sigma_2)$ satisfy~\eqref{eq:IO}, and let $\varsigma = (\varsigma_1, \varsigma_2)$ satisfy~\eqref{eq:mix-order}. An ideal MIP formulation of ~\eqref{eq:bifunc} is given as follows:
\[
\begin{aligned}
    \Bigl\{(\x,\z,\mu)  \Bigm|\,  &  \eqref{eq:unary},\  \z' = L^{\sigma}(\z),\ \eqref{eq:staircase},  \\
    &  \qquad \z'' = L^{\varsigma}(\z),\  \mu \geq		\psi^\varsigma(\j^0_\pi) + \sum_{s  \in [D]} \bigl( \psi^\varsigma(\j^s_\pi) - \psi^\varsigma(\j^{s-1}_\pi) \bigr) \cdot z''_{\pi_s}  \for \pi \in \Pi \Bigr\},
\end{aligned}
\]
where $\psi^\varsigma(j_1,j_2) := f_1(p_{1,\varsigma_1(j_1)})f_2(p_{2,\varsigma_2(j_2)})$ for $(j_1,j_2) \in \G$,  and  \(L^\sigma\) and \(L^\varsigma\) are defined by~\eqref{eq:variable-switching-general}.
\end{proposition}
Here, we provide a sketch of the proof. By Lemma~\ref{lemma:MIP-discrete}, it suffices to derive the convex hull of the graph of $\ell_1\ell_2$ over $\vertex(\Delta)$. The argument proceeds in two parts. For the hypograph, we apply the affine transformation associated with a permutation tuple $\sigma$ that orders the inner functions in non-decreasing order and then use staircase inequalities for the reordered function. The permutation $\sigma$ is chosen so that the reordered function is supermodular. For the epigraph, we follow an analogous procedure, but use a different permutation tuple $\varsigma$ chosen to make the reordered function submodular.

\subsection{Modeling via a logarithmic binarization}\label{section:otherbinarization}

The formulations developed in Section~\ref{section:composite-MIP} rely on the unary binarization, which represents each discrete variable $x_i$ using $d_i$ binary variables. In this section, we modify these formulations to reduce the number of binary variables to $\lceil \log_2(d_i+1) \rceil$. Before presenting the resulting formulations, we establish a general result that characterizes when ideal formulations can be preserved under alternative binarization schemes. Lemma~\ref{lemma:MIP-discrete} shows that combining the unary binarization with a convex hull description yields an ideal formulation for the graph of a discrete composite function $\phi \mcirc f$. We extend this result by showing that any ideal (resp. valid) formulation of the vertices of the unary simplotope can be combined with the same convex hull construction to produce an ideal (resp. valid) formulation for the composite function.

\begin{lemma}\label{lemma:MIP-discrete-log}
For each $i \in [n]$, let $r_i$ be a positive integer and $Q_i$ be a polytope in $(\z_i,\bdelta_i) \in \R^{d_i} \times \R^{r_i}$. Assume that $Q_i \cap \bigl(\R^{d_i} \times \{0,1\}^{r_i}\bigr)$ is an MIP formulation of $\vertex(\Delta^{d_i})$.
Then, an MIP formulation of $\graph( \phi \mcirc f )$ is given as follows 
\begin{equation*}\label{eq:MIP-discrete-log}
\begin{aligned}
\Bigl\{(\x, \z, \mu, \bdelta) \Bigm| \; (\z,\mu) &\in \conv\bigl(\graph(\phi \mcirc \ell)\bigr),\ \x \in P, \\
& x_i = \pt_{i0} + \sum_{ j \in [d_i]} (\pt_{ij} - \pt_{ij-1}) z_{ij},\ (\z_i, \bdelta_i) \in Q_i,\ \bdelta_i \in \{0,1\}^{r_i}\, \for i \in [n] \Bigr\}.
\end{aligned}
\end{equation*}
This formulation is ideal if $Q_i \cap \bigl(\R^{d_i} \times \{0,1\}^{r_i}\bigr)$ is ideal and the constraint $x\in P$ is relaxed. 
\end{lemma}

We next illustrate the use of Lemma~\ref{lemma:MIP-discrete-log} in reducing the number of binary variables in our formulations to $\sum_{i \in [n]} \lceil \log_2(d_i+1) \rceil$. To do so, we modify the formulation in Theorem~2 of~\cite{vielma2011modeling} to model $\vertex(\Delta^{d_i})$ as follows:

\begin{equation}\label{eq:log}
\begin{aligned}
\z_i \in \Delta^{d_i} \quad \gamma_{ik} \in \{0,1\} \quad  & \sum_{j \in A_{ik}}     (z_{ij} - z_{ij+1}) \leq \gamma_{ik} \quad \text{and } \\
&\sum_{j \notin A_{ik}} (z_{ij} - z_{ij+1}) \leq 1-\gamma_{ik} \quad \for  k =1, \ldots \lceil \log_2(d_i+1)\rceil,
\end{aligned}\tag{\textsc{Log}}    
\end{equation}
where $A_{ik}: = \bigl\{ j \in 0 \cup [d_i]  \bigm| \text{$k^{\text{th}}$ coordinate in the binary representation of $j$ is $1$} \bigr\}$, and for convenience, we let $z_{i0} = 1$ and $z_{i(d_i+1)} = 0$.  Replacing the binary restrictions on $\z$ in Theorem~\ref{them:supermodular-model} with the logarithmic constraints~\eqref{eq:log} yields the following logarithmic counterpart of Theorem~\ref{them:supermodular-model}.

\begin{corollary}\label{cor:supermodular-model-log}
Let $\sigma:=(\sigma_1, \ldots, \sigma_n )$ satisfying~\eqref{eq:IO}.  
If $\phi(\cdot)$ is supermodular then an ideal MIP formulation of $ \hypo(\phi \mcirc f) $  is 
	\[
	\Biggl\{(\x,\z, \boldsymbol{\gamma}, \mu) \Biggm| x_i = \pt_{i0} + \sum_{ j \in [d_i]} (\pt_{ij} - \pt_{ij-1}) z_{ij} \text{ and } ~\eqref{eq:log} \; \for i \in [n] ,\ \z'= L^{\sigma}(\z),\ \eqref{eq:staircase}
	 \Biggr\},
	\]
 where $L^{\sigma}$ is the mapping~\eqref{eq:variable-switching-general}.
\end{corollary}

\section{Strong Formulations for Price Promotion Optimization}\label{section:ppo}

We consider two demand models. The first one is a \textit{linear} model with additive linear cross-item effects:
\begin{equation}\label{eq:demand-linear}
	D^i_t(\x) = \overbrace{(a^i_t - b^i_0x^i_t)}^{\text{self effect}} + \overbrace{\sum_{m \in [M_i]} b^i_m x^i_{t-m}}^{\text{historical effects}} + \sum_{k \in [N]\setminus \{i\}} \overbrace{\delta^{ik}x^k_t}^{\text{cross-item effect}} ,
\end{equation}
and the second one is a \textit{log-log} model with additive linear cross-item effects:
\begin{equation}\label{eq:demand-loglog}
    D^i_t(\x) = \overbrace{a^i_t \cdot (x^i_t)^{-b^i_0}}^{\text{self effect}} \cdot \overbrace{\prod_{m \in [M_i]}(x^i_{t-m})^{b^i_m}}^{\text{historical effects}} + \sum_{k \in
[N]\setminus \{i\}}\overbrace{\delta^{ik}x^k_t}^{\text{cross-item effect}}.
\end{equation}
Both models are characterized by shared additively separable cross-item effects but differing in their self effects and historical effects. Specifically, the term $a^i_t$ captures seasonality effects, $b^i_0$ represents the price-sensitivity factor, $b^i_m$ quantifies the influence of historical prices on current demand, and $\delta^{ij}$ encodes cross-item effects. We assume that $b^{i}_{m}$ is non-negative for all $m \in [M_i]$, reflecting non-decreasing historical price effects. This historical effect is termed as the \textit{post-promotion dip effect}, which models the transient decline in demand following a promotion. Furthermore, for any pair of  items $i$ and $k$, the cross-item parameter $\delta^{ik}_t$ can take positive or negative values. The pair is classified as \textit{complementary} if $\delta^{ik}_t<0$ (indicating mutual demand enhancement) and \textit{substitutable} if  $\delta^{ik}_t>0$ (representing demand substitution).

Our formulation techniques, especially staircase inequalities, enable MIP formulations for PPO under the two demand models. A limitation of naively applying staircase inequalities is their exponential growth with the number of periods $M_i$ that influence demand and the number of price levels $d_i$. To mitigate this, we develop compact MIP formulations in Section~\ref{section:pop-compact}, retaining a subset of the inequalities while maintaining exactness. In Section~\ref{section:pop-cuttingplane}, we integrate the remaining inequalities using efficient separation algorithms, theoretically yielding new polynomial-time solvable cases. The computational efficacy of the proposed formulations is evaluated in Section~\ref{section:computation}.

\subsection{Compact formulations}\label{section:pop-compact}
We model the price ladder constraint $x_t^i \in \{p^i_0,p^i_1, \ldots, p^i_{d_i} \}$ using the unary binarization scheme. 
Under this representation, the feasible region of PPO can be expressed as:
\begin{subequations}\label{eq:brulesMIP}
	\begin{align}
     &  1 \geq z^i_{t1}  \geq \cdots \geq z^i_{td_i} \geq 0 \text{ and } \z^i_t \in  \{0,1\}^{d_i} && \for i \in [N] \text{ and } t \in [T] \label{eq:brulesMIP-0}  \\
&x^i_t = p^i_{0} + \sum_{j \in [d_i]} (p^i_j - p^i_{j-1})z^i_{tj}  && \for i \in [N] \text{ and } t \in [T]  \label{eq:brulesMIP-1} \\
&T- \sum_{t \in [T]}z^i_{td_i} \leq K^i &&  \for i \in [N]\label{eq:brulesMIP-2}  \\
&N - \sum_{i \in [N]}z^i_{td_i} \leq C_t &&  \for t \in [T], \label{eq:brulesMIP-3}  
	\end{align}
\end{subequations}
where constraints~\eqref{eq:brulesMIP-1}-\eqref{eq:brulesMIP-3} reformulate constraints~\eqref{eq:brules1}-\eqref{eq:brules3};  in the presence of~\eqref{eq:brulesMIP-1},  $z^i_{td_i} = 1$ if and only if $\ind \{x^i_t < p^i_{td_i}\} = 0$. Our formulations support alternative binarization schemes, such as full or logarithmic, and we evaluate their performance in Section~\ref{sec:computation-comparision} after presenting these detailed formulations.

 \subsubsection{Linear demand model}\label{sec:linearDemandModel}
 We write the objective function in PPO as:
 \[
\sum_{i \in [N]} \sum_{t \in [T]}x^i_t \cdot (a^i_t-b^i_0x^i_t) + \sum_{m \in [M_i]} b^i_m \cdot x^i_t \cdot x^i_{t-m} 
 + \sum_{k \in [N] \setminus i}  \delta^{ik} \cdot x^{i}_{t} \cdot x^k_t -  c^i_t\cdot D_t^i(\x),
 \]
{\revision where $D^i_t(\cdot)$ is the linear demand function defined as in~\eqref{eq:demand-linear}.} This, thereby, decomposes the promotion effect on revenue into three parts: self effects, historical effects, and cross-item effects. For each item $i$ and each period $t$, we introduce a variable $s^i_t$ to represent the self-effect, and further express it as the sum of incremental self-effects:
\begin{equation}\label{eq:linear-self}
    s^i_t =  w^i_t(p^i_0) + \sum_{j \in [d_i]} \Bigl(w^i_t(p^i_j) - w^i_t(p^i_{j-1})  \Bigr)  z^i_{t,j},
\end{equation}
where $w^i_t(x^i_t):=x^i_t \cdot (a^i_t-b^i_0x^i_t)$.  Unlike the self-effect, the historical effects and cross-item effects are bilinear. To construct exact MIP representations for such effects, we specialize Proposition~\ref{prop:bilinear-sub} to the bilinear terms. More specifically, due to the post-promotion dip effect, \textit{i.e.}, $b^i_k > 0$, it suffices to impose the constraint $u^i_{tm} \leq x^i_tx^{i}_{t-m}$. Instead of using all $\binom{2d_i}{d_i}$ staircase inequalities, we will show in Proposition~\ref{prop:base-linear} that the following $d_i+1$ inequalities  are sufficient to obtain an exact formulation:
\begin{equation}\label{eq:linear-model-period}
\begin{aligned}
u^i_{tm} \leq p^i_0p^i_0  + \sum_{j =1}^\tau (p^i_{j} - p^i_{j-1})p^i_{0}z^i_{t,j} & + \sum_{j =\tau + 1}^{d_i} (p^i_{j} - p^i_{j-1})p^i_{d_i}z^i_{t,j} \\
&+  \sum_{j =1 }^{d_i}p^i_{\tau}(p^i_{j} - p^i_{j-1}) z^i_{t-m,j} \quad  \for \tau \in  0 \cup [d_i] \; m \in [M_i],
\end{aligned}
\end{equation}
where the second through fourth terms model the incremental historical effects resulting from changes in the price of item $i$ in periods $t$ or $t-m$ in a sequence, determined by the corresponding staircase. 

To capture both complementary and substitutable cross-item effects, we  need to use  two separate sets of staircase inequalities. For each item $i$, we classify other items into two groups, its substitutes  $A_i:= \bigl\{ k \in [N] \setminus \{i\} \bigm| \delta^{ik} >0\bigr\}$ and complements $B_i := \bigl\{ k \in [N] \setminus \{i\} \bigm| \delta^{ik} <0 \bigr\}$. For the substitutable effects, we impose constraint $v^{ik}_t \leq x^i_tx^k_t$, and its exact formulation is given by the following staircase inequalities:
\begin{equation}\label{eq:linear-model-item-over}
\begin{aligned}
v^{ik}_{t} \leq p^i_0p^k_0  &+ \sum_{j =1}^\tau (p^i_{j} - p^i_{j-1})p^k_{0}z^i_{t,j}  + \sum_{j =\tau + 1}^{d_i} (p^i_{j} - p^i_{j-1})p^k_{d_k}z^i_{t,j} \\
&+  \sum_{j =1}^{d_k}p^i_\tau(p^k_{j} - p^k_{j-1}) z^k_{t,j} \quad  \for \tau \in  0 \cup [d_i]  \text{ and }  k \in A_i. 
\end{aligned}
\end{equation}
For the complementary effects, we impose constraint $v^{ik}_t \geq x^i_tx^k_t$, and its exact formulation is:
\begin{equation}\label{eq:linear-model-item-under}
\begin{aligned}
v^{ik}_{t} \geq p^i_0p^k_{d_k}  &+ \sum_{j =1}^\tau (p^i_{j} - p^i_{j-1})p^k_{d_k}z^i_{t,j}  + \sum_{j =\tau + 1}^{d_i} (p^i_{j} - p^i_{j-1})p^k_{0}z^i_{t,j} \\
&+  \sum_{j =1}^{d_k}p^i_\tau(p^k_{d_k-j} - p^k_{d_k-j+1}) (1- z^k_{t,d_k-j +1}) \quad  \for \tau \in  0 \cup [d_i]  \text{ and }  k \in B_i.    
\end{aligned}
\end{equation}
Incorporating these representations, we obtain a compact MIP formulation of PPO:  
 \begin{equation}\label{eq:linear-base}
	\begin{aligned}
		\max \quad & \sum_{i \in [N]} \sum_{t \in [T]}  s^i_t + \sum_{m \in [M_i]} b^i_k \cdot u^{i}_{tm} 
 + \sum_{k \in [N] \setminus \{i\}}\delta^{ik} \cdot v^{ik}_t - c^i_t \cdot D^i_t(\x)
 \\
		\text{s.t.} \quad &  \eqref{eq:brulesMIP}  \text{ and }  \eqref{eq:linear-self}\text{-}\eqref{eq:linear-model-item-under} \, \for i \in [N],\, t \in [T],
	\end{aligned}\tag{\textsc{Base-Linear}}
\end{equation}
{\revision where $D^i_t(\cdot)$ is the linear demand function defined as in~\eqref{eq:demand-linear}.} Note that this formulation requires {\revision at most $T\sum_{i \in [N]}M_i(d_i+1) + TN(N-1)\sum_{i \in [N]} d_i$ linear inequalities in addition to constraints in~\eqref{eq:brulesMIP} and \eqref{eq:linear-self}. The first term counts the inequalities associated with the $M_i$ memory length for each item $i$, while the second term provides an upper bound on the cross-item inequalities by considering the $N-1$ possible other items $k\neq i$. The actual number can be smaller since the substitute and complement sets are item-specific.}

 \begin{proposition}\label{prop:base-linear}
Assume that the demand model is given as in~\eqref{eq:demand-linear} with the presence of the post-promotion dip effect. Then, an exact MIP formulation of PPO is given by~\eqref{eq:linear-base}.
\end{proposition}
Both~$\text{App(2)}$ given in \citet{cohen2021promotion} and \eqref{eq:linear-base} are exact MIP formulations. However, $\text{App(2)}$ uses more variables than \eqref{eq:linear-base}.

\subsubsection{Log-log demand model}\label{section:loglog}
We consider the log-log demand model with additive linear cross-item effects and multiplicative historical effects. The objective function of PPO can be expressed as follows:
\begin{equation*}
\sum_{i \in [N]} \sum_{t \in [T]}  (x^i_t-c^i_t) \cdot a^i_t \cdot (x^i_t)^{-b^i_0} \cdot \prod_{m \in [M_i]}(x^i_{t-m})^{b^i_m} + \sum_{k \in [N]\setminus \{i\}} \delta^{ik}(x_t^ix^k_t - c^i_t x^k_t),
\end{equation*}
where the first term captures the multiplicative relation between the self-effect, denoted as $w^i_t(x^i_t):= (x^i_t-c^i_t) \cdot a^i_t \cdot (x^i_t)^{-b^i_0}$, and historical effects while the second term captures the cross-item effects. We will focus on treating the first term since the cross-item effects can be handled using~\eqref{eq:linear-model-item-over} and~\eqref{eq:linear-model-item-under}. 

We recursively express the multiple time horizon historical effect using consecutive period effects:
\begin{equation}\label{eq:period-recursion}
u^i_{t1} = (x^i_{t-1})^{b^i_1} \quad \text{ and } \quad u^i_{tm} \leq  (x^i_{t-m})^{b^i_m} u^i_{tm-1}  \quad \for \; m =2, \ldots, M_i.     
\end{equation}
To represent these nonlinear effects with staircase inequalities, we discretize $u^i_{tm}$ using its minimum and maximum values, denoted as $(u^i_{t m})_{\min}$ and $(u^i_{t m})_{\max}$. We then introduce recursive  grid representations for discrete values, that is,  $\zeta^i_{t1}(j) = (p^i_{j})^{b^i_1}$ for  $j \in 0 \cup [d_i]$, and for $m \geq 2$, $\zeta^i_{tm}:\{0, 1, \ldots, d_i \} \times \{0,1\}$ such that 
\[
  \zeta^i_{tm}(j,0) = (p^i_{j})^{b^i_m}(u^i_{t m-1})_{\min} \quad \text{  and } \quad \zeta^i_{tm}(j,1) = (p^i_{j})^{b^i_m}(u^i_{t m-1})_{\max} \qquad \for j \in 0 \cup [d_i].
\]
Now, specializing Proposition~\ref{prop:bilinear-sub} to the discretized nonlinear effects on two consecutive periods, we obtain the following system of staircase inequalities:
\begin{equation}\label{eq:log-crossperiod}
\begin{aligned}
& u^i_{t1} = \zeta^i_{t1}(0) + \sum_{j = 1}^{d_i}\Bigl(\zeta^i_{t1}(j) - \zeta^i_{t1}(j-1)\Bigr) \cdot z^i_{t-1, j} \\
    &u^i_{tm} \leq \zeta^i_{tm}(0,0) + \sum_{j = 1}^{\tau} \Bigl(\zeta^i_{tm}(j,0) - \zeta^i_{tm}(j-1,0)\Bigr) \cdot z^i_{t-m,j} + \sum_{j = \tau}^{d_i} \Bigl(\zeta^i_{tm}(j,1) - \zeta^i_{tm}(j-1,1)\Bigr) \cdot z^i_{t-m,j} \\ 
    & \qquad \qquad  + \Bigl(\zeta^i_{tm}(\tau,1) - \zeta^i_{tm}(\tau,0)\Bigr) \cdot \frac{u^i_{t,m-1} - (u^i_{t,m-1})_{\min}}{(u^i_{t,m-1})_{\max} - (u^i_{t,m-1})_{\min}}  \qquad \for \tau  \in 0 \cup [d_i]  \text{ and } {\revision m \geq 2}.
\end{aligned}
\end{equation}

Last, consider the multiplicative relation between the self effect and the historical effect. Note the self-effect may not be monotone. To handle this, for each period $t$ and each product $i$, we introduce a permutation $\sigma^i_t:\{0,1, \ldots, d_i \} \to \{0,1, \ldots, d_i \} $ such that the self-effect is monotone on the sorted prices, that is, 
\begin{equation}\label{eq:permute-self-effects}
    w^i_t\bigl(p^i_{\sigma^i_t(0)}\bigr) \leq w^i_t\bigl(p^i_{\sigma^i_t(1)}\bigr) \leq \cdots \leq w^i_t\bigl(p^i_{\sigma^i_t(d_i)}\bigr). 
\end{equation}
In addition,  we use this permutation to introduce a grid function  to discretize the multiplicative relation: 
\[
\eta^i_t(j,0)=w^i_t\bigl(p^i_{\sigma^i_t(j)}\bigr)  (u^i_{tM_i})_{\min} \quad \text{ and } \quad \eta^i_t(j,1)=w^i_t\bigl(p^i_{\sigma^i_t(j)}\bigr) (u^i_{tM_i})_{\max} \quad \for j \in 0 \cup [d_i].
\]
This permutation also induces a corresponding swap of the unary simplex vertices, which can be represented by the linear transformation defined in~\eqref{eq:variable-switching-general}, that is, 
\begin{equation}\label{eq:swtich-self-effect}
    L^{i}_t (\z^i_t) = \Lambda_i^{-1}P^{\sigma^i_t} \Lambda_i  \z^i_t,
\end{equation}
where $\Lambda_i$ is defined as in~\eqref{eq:fromztolambda} and $P^{\sigma^i_t}$ is the permutation matrix given by $\sigma^i_t$. Using these, we obtain
\begin{equation}\label{eq:self-crossperiod}
\begin{aligned}
& \hat{\z}^i_t = L^i_t(\z^i_t) \\ 
   & s^i_t \leq \eta^i_t(0,0) + \sum_{j = 1}^\tau\Bigl(\eta^i_t(j,0) - \eta^i_t(j-1,0) \Bigr)\cdot \hat{z}^i_{t,j} +   \sum_{j = \tau + 1}^{d_i}\Bigl(\eta^i_t(j,1) - \eta^i_t(j-1,1) \Bigr)\cdot \hat{z}^i_{t,j}  \\
     & \qquad \qquad  + \Bigl(\eta^i_{t}(\tau,1) - \eta^i_{t}(\tau,0)\Bigr) \cdot \frac{u^i_{tM_i} - (u^i_{tM_i})_{\min}}{(u^i_{tM_i})_{\max} - (u^i_{tM_i})_{\min}}  \qquad \for \tau  \in 0 \cup [d_i].
     \end{aligned}
\end{equation}
Observe that $\hat{z}^i_t$ permutes $p^i$ so that $w^i_t(\cdot)$ is non-decreasing just as, in Theorem~\ref{them:supermodular-model},  $\z'$ permutes the discrete vector $\x$ so that $f$ became non-decreasing. The constraints~\eqref{eq:log-crossperiod} and~\eqref{eq:self-crossperiod} together provide an exact representation of the multiplicative relationship between the self-effect and the history effect. In total, they comprise $T\sum_{i\in [N]}(d_i+1)(M_i+1)$ linear inequalities.   Using these expressions, we obtain:
\begin{equation}\label{eq:formulation-multi}
	\begin{aligned}
		\max \quad & \sum_{i \in [N]} \sum_{t \in [T]} s_t^i
 + \sum_{k \in [N] \setminus \{i\} }  \delta^{ik}( v^{ik}_t - c^
i_tx^i_t) \\
		\text{s.t.} \quad &  \eqref{eq:brulesMIP} \text{ and}~\eqref{eq:linear-model-item-over}, \eqref{eq:linear-model-item-under}, ~\eqref{eq:log-crossperiod},~\eqref{eq:self-crossperiod} \for i \in [N],\; t \in [T] .
	\end{aligned}
	\tag{\textsc{Base-LogLog}}
\end{equation}
{\revision Note that~\eqref{eq:linear-model-item-over} and~\eqref{eq:log-crossperiod} are hypograph inequalities. This suffices because the objective is to maximize profit and the self-effect and historical effects are nonnegative, ensuring that the auxiliary variables attain their tight values at optimality.}
\begin{proposition}\label{prop:base-loglog}
	Assume that the demand model is as given in \eqref{eq:demand-loglog}, which includes cross-item effects and  multiplicative post-promotion dip effect. Then, an MIP formulation of PPO is given by~\eqref{eq:formulation-multi}.  
\end{proposition}
For the log-log demand model, \citet{cohen2021promotion} focus on finding approximate solutions. When the demand model is as given in \eqref{eq:demand-loglog}, their integer-programming formulations, $\App(2), \ldots, \App(n)$, do not model PPO exactly while our formulation~\eqref{eq:formulation-multi} is an exact formulation. 

\subsection{Stronger formulations via fast separations and their theoretical performance}\label{section:pop-cuttingplane}
In Section~\ref{section:supermodular}, we characterized the staircase inequalities for each supermodular term, and in Section~\ref{section:pop-compact}, we used a subset of them to develop compact formulations for  PPO. This subsection details efficient separation algorithms that integrate the remaining inequalities. These algorithms,
described in Algorithms~\ref{alg:cross-item-effects}-\ref{alg:cutting-plane}, repeatedly
utilize the following operator, defined for a grid valuation $\psi: \prod_{i \in [n]}\{0,1 , \ldots, d_i\} \to \R$ and a direction representation of a staircase $\pi\in \Pi$ as
\[
\begin{aligned}
    \stair(\psi)^\pi(\z) :=  \psi(\j^0_\pi) + \sum_{s \in [D] }(\psi(\j^s_\pi)-\psi(\j^{s-1}_\pi))z_{\pi_s} ,
\end{aligned}
\]
where $D=\sum_{i \in [n]}d_i$ and $(\j^0_\pi, \j^1_\pi, \ldots, \j^D_\pi)$ represents a grid-point as in~\eqref{eq:pointrep}. 
For a given $\bar{\z}$ in the simplotope, we use $\texttt{sort-and-track}(\bar{\z})$ to generate a staircase $\pi$ where $\bar{z}_{\pi_1} \leq \bar{z}_{\pi_2} \leq \cdots \leq  \bar{z}_{\pi_N}$.

{\revision
\begin{example}\label{ex:interpret_coefficients}
We consider the setup of Examples~\ref{ex:reorder}--\ref{ex:cuts} and interpret the coefficients of the staircase inequality \eqref{eq:stairIneqEx} in the context of PPO. Specifically, the coefficient of $z_{31}$, which is $\frac{64}{81} +\frac{147}{16}$, can be understood as the marginal profit resulting from a price adjustment in period 3.

In this scenario, the item was previously offered (in periods $1$ and $2$) at a price of $2$. In period 3, the price changes from $1$ to $3$. Offering the item at a price of $1$ yields no profit in that period, as the price equals the cost. However, increasing the price to $3$ triggers a self-effect that reduces the sales volume from $2^5$ to $\frac{1}{81}2^5$. Since each unit sold at this higher price generates a profit of $2$, the total marginal profit is $2\times\frac{1}{81}2^5 = \frac{64}{81}$. This is the first term in the coefficient of $z_{31}$. However, $z_{31}$ continues to receive coefficients from price changes until the price in period $3$ reaches $2$ (second price level). Here, this occurs when the item is offered for a price of $3$ in periods $1$ and $2$ and its price is dropped in period $3$ from $3$ to $2$. As a result the marginal profit change is $\bigl(1\times \frac{1}{16}-2\times \frac{1}{81})3^5 = \frac{147}{16}$. The coefficient of $z_{31}$ aggregates these profit changes in period $3$ along the staircase. \Halmos
\end{example}

\begin{remark}\label{rmk:interpret}
We remark that the interpretation provided in Example~\ref{ex:interpret_coefficients} holds generally for the log-log demand model. First, we choose a permutation $\sigma:\{1,\ldots,d_i\}\rightarrow \{1,\ldots,d_i\}$ and explore the transformation in \eqref{eq:self-crossperiod} relating $(\hat{z}^i_{tr})_{r=0}^{d_i-1}$ to $(z^i_{tj})_{j=0}^{d_i-1}$. It follows that
\begin{equation*}
\begin{split}
\hat{z}_{tr} &= \sum_{\ell = r+1}^{d_i} \lambda^i_{t\sigma^{-1}(\ell)} = \sum_{\ell=r+1}^{d_i} \left(z^i_{t\sigma^{-1}(\ell)-1} - z^i_{t\sigma^{-1}(\ell)}\right)\\
&=\sum_{j:r+1\le \sigma(j+1)\le d_i} z^i_{tj} - \sum_{j:r+1\le \sigma(j) \le d_i} z^i_{tj}
=\sum_{j:\sigma(j)\le r < \sigma(j+1)} z^i_{tj} - \sum_{j:\sigma(j+1)\le r < \sigma(j)} z^i_{tj},
\end{split}
\end{equation*}
where $j,r\in\{0,\ldots,d_i-1\}$ and $\sigma(0)$ is assumed to be $0$. We define the index sets $J_\sigma(j) = \bigl\{r\mid r\in [\sigma(j),\sigma(j+1))\bigr\}$ and $\bar{J}_{\sigma}(j) = \bigl\{r\mid r\in [\sigma(j+1),\sigma(j))\bigr\}$. For the setting of Example~\ref{ex:interpret_coefficients} with $j=1$ and $[\sigma(i)]_{i=1}^3=[1,3,2]$, $J_\sigma(1) = \{1,2\}$ and $\bar{J}_\sigma(1)=\emptyset$. Note that for any $j$, exactly one of these sets is non-empty, depending on whether $\sigma(j) < \sigma(j+1)$ or $\sigma(j) > \sigma(j+1)$. Thus, we can write
\begin{equation}\label{eq:z-hattoz}
    \hat{z}^i_{tr} = \sum_{j:r\in J_\sigma(j)} z^i_{tj} - \sum_{j:r\in \bar{J}_\sigma(j)} z^i_{tj}.
\end{equation}

Now, we show how \eqref{eq:z-hattoz} interprets the coefficients of the staircase inequality. Let $\hat{w}_i(r)$ be the self-effect multiplied by the profit per unit sale of $i$ in period $t$ at price $p^t_{\sigma(r)}$, \textit{i.e.}, $\hat{w}_i(r) = (p^i_{\sigma(r)}-c^i_t)(p^i_{\sigma(r)})^{-b^i_0}$. We choose $\sigma$ such that $\hat{w}_i(r) \le \hat{w}_i(r')$ whenever $r < r'$. The time-period effect for product $i$ at period $t$ combines self and historical effects: $a^i_t (x^i_t-c^i_t)(x^i_t)^{-b^i_0} \prod_{m \in [M_i]}(x^i_{t-m})^{b^i_m}$. Each staircase step increases this effect by perturbing the offering price in one period $\{t,\ldots,t-M_i\}$. The coefficient $\text{coeff}(\hat{z}^i_{tr})$ represents the change in this effect when the price in period $t$ changes from $p^i_{\sigma^{-1}(r-1)}$ to $p^i_{\sigma^{-1}(r)}$, with previous prices fixed.
Using \eqref{eq:z-hattoz}, the coefficient of $z^i_{tj}$ is interpreted as follows. If $\sigma(j) < \sigma(j+1)$, then
\begin{equation*}
    \text{coeff}(z^i_{tj}) = \sum_{r\in J_\sigma(j)} \text{coeff}(\hat{z}^i_{tr}).
\end{equation*}
In Example~\ref{ex:interpret_coefficients}, this summation is taken over two profit contributions since $J_\sigma(1)=\{1,2\}$. Otherwise, if $\sigma(j) > \sigma(j+1)$,
\begin{equation*}
    \text{coeff}(z^i_{tj}) = -\sum_{r\in \bar{J}_\sigma(j)} \text{coeff}(\hat{z}^i_{tr}).
\end{equation*}
In both cases, the summation aggregates the marginal profit contributions from the time-period effect with a sequence of price changes for product $i$ in period $t$. \Halmos
\end{remark}
}

\begin{algorithm}
\fontsize{10pt}{11pt}\selectfont 
\renewcommand{\baselinestretch}{1.2}\selectfont
\caption{Separation for cross-item effects}\label{alg:cross-item-effects}
\KwData{a point $(\bar{\z}, \bar{\boldsymbol{v}})$ with $\bar{\z} \in \Delta$}
\KwResult{\texttt{cut} \tcc{a set of staircase inequalities} }
  \For{$t \in [T]$ and $i,k \in [N]$, with $i\neq k $ }{

\eIf{ the pair of items $(i,k)$ are substitutable}{
	$\pi \gets$ \texttt{sort-and-track}($\bar{\z}^i_t$, $\bar{\z}^k_t$)\;
    $\psi(j_i,j_k) = p^i_{j_i}p^k_{j_k}$ for $j_i = 0,1, \ldots, d_i$ and $j_k = 0, 1, \ldots, d_k$\;
    $v^* \gets \stair(\psi)^{\pi}(\bar{\z}^i_t,\bar{\z}^k_t)$\;
     \lIf{ $ \bar{v}^{ik}_t > v^*$ }{ push  $v^{ik}_{t} \leq \stair(\psi)^\pi(\z^i_t,\z^k_t)$ into \texttt{cut}}  }
     {
  $\bar{\boldsymbol{\zeta}}^{k}_t \gets (1-\bar{z}^k_{t,d_k}, 1- \bar{z}^k_{t,d_k-1}, \ldots,   1- \bar{z}^k_{t,1})$\;
	 $\pi \gets$ \texttt{sort-and-track}($\bar{\z}^i_t$, $\bar{\boldsymbol{\zeta}}^k_t$)\;
  $\psi(j_i,j_k) = p^i_{j_i}p^k_{d_k - j_k}$ for $j_i = 0,1, \ldots, d_i$ and $j_k = 0, 1, \ldots, d_k$\;
      $v^* \gets \stair(\psi)^{\pi}(\bar{\z}^i_t,\bar{\boldsymbol{\zeta}}^k_t)$\;
        $\boldsymbol{\zeta}^{k}_t \gets (1-z^k_{t,d_k}, 1- z^k_{t,d_k-1}, \ldots,   1- z^k_{t,1})$\;
      \lIf{ $ \bar{v}^{ik}_t < v^*$ }{ push $v^{ik}_{t} \geq \stair(\psi)^\pi(\z^i_t,\boldsymbol{\zeta}^k_t)$ into \texttt{cut}}
  }
   }
\end{algorithm}

We initially modeled cross-item effects with a limited number of staircases, as described in \eqref{eq:linear-model-item-over} and \eqref{eq:linear-model-item-under}. We then supplemented these cuts using Algorithm~\ref{alg:cross-item-effects}, which generates staircase inequalities for each item pair $(i,k)$. For substitute items (supermodular terms), we sort the pairs $(\bar{\z}^i_t$, $\bar{\z}^k_t)$. Conversely, for complementary pairs, we reorder the price ladder in decreasing order and generate staircase inequalities.

For the historical effects, we introduce Algorithms~\ref{alg:add-cross-period-effects} and~\ref{alg:mult-cross-period-effects} to handle the additive and multiplicative cases, respectively.  In the additive case, for each item $i$, period $t$, and past period $m$, we generate the staircase inequality by sorting $(\bar{\z}^i_t, \bar{\z}^i_{t-m})$. In the multiplicative case, for each item $i$ and period $t$, we use a permutation $\sigma^i_t$ to sort the self-effects alongside the corresponding $\bar{\z}^i_t$, producing the transformed vector $\bar{\boldsymbol{\zeta}}^i_t$. We then jointly sort $(\bar{\boldsymbol{\zeta}}^i_t,\bar{\z}^i_{t-1}, \ldots , \bar{\z}^i_{t-M_i})$ to construct the staircase inequality.

\begin{algorithm}
\fontsize{10pt}{11pt}\selectfont 
\renewcommand{\baselinestretch}{1.2}\selectfont
\caption{Separation for additive cross-period effects}\label{alg:add-cross-period-effects}
\KwData{a point $(\bar{\z}, \bar{\boldsymbol{u}})$ with $\bar{\z} \in \Delta$}
\KwResult{\texttt{cut} \tcc{a set of staircase inequalities} }
  \For{$t \in [T], i \in [N] \text{ and } m \in [M_i] $ }{
	$\pi \gets$ \texttt{sort-and-track}($\bar{\z}^i_t$, $\bar{\z}^i_{t-m}$)\;
    $\psi(j_i,j_k) = p^i_{j_i}p^i_{j_k}$ for $j_i = 0,1, \ldots, d_i$ and $j_k = 0, 1, \ldots, d_i$\;
    $u^* \gets \stair(\psi)^{\pi}(\bar{\z}^i_t,\bar{\z}^i_{t-m})$\;
     \lIf{ $ \bar{u}^{i}_{tm} > u^*$ }{ push  $u^{i}_{tm} \leq \stair(\psi)^\pi(\z^i_t,\z^i_{t-m})$ into \texttt{cut}}  
   }
\end{algorithm}

\begin{algorithm}
\fontsize{10pt}{11pt}\selectfont 
\renewcommand{\baselinestretch}{1.2}\selectfont
\caption{Separation for multiplicative historical effects}\label{alg:mult-cross-period-effects}alg:
\KwData{a point $(\bar{\z}, \bar{\boldsymbol{s}})$ with $\bar{\z} \in \Delta$}
\KwResult{\texttt{cut} \tcc{a set of staircase inequalities} }
  \For{$t \in [T] \text{ and } i \in [N] $ }{
     let $\sigma^i_t$  be the permutation sorting the self-effects as in~\eqref{eq:permute-self-effects}\;
     let $L^i_t$ be the linear transformation defined as in~\eqref{eq:swtich-self-effect}\;
     $\bar{\boldsymbol{\zeta}}^i_t \gets L^i_t(\bar{\z}^i_t)$ \text{ and } $\boldsymbol{\zeta}^i_t \gets   L^i_t(\z^i_t)$\;
	$\pi \gets$ \texttt{sort-and-track}$(\bar{\boldsymbol{\zeta}}^i_t, \bar{\z}^i_{t-1}, \ldots,  \bar{\z}^i_{t-M_i})$\;
    $\psi(\boldsymbol{j}) = w^i_t(p^i_{\sigma(j_0)}) \prod_{m \in [M_i]}(p^i_{j_m})^{b^i_m} $ for $\boldsymbol{j} = (j_0,j_1, \ldots, j_{M_i}) \in \{0,1,\ldots, d_i\}^{M_i+1}$ \;
    $s^* \gets \stair(\psi)^{\pi}(\bar{\boldsymbol{\zeta}}^i_t, \bar{\z}^i_{t-1}, \ldots,  \bar{\z}^i_{t-M_i})$\;
     \lIf{ $ \bar{s}^{i}_{t} > s^*$ }{ push  $s^{i}_{t} \leq \stair(\psi)^\pi(\boldsymbol{\zeta}^i_t ,\z^i_{t-1}, \ldots,  \z^i_{t-m})$ into \texttt{cut}}  
   }
\end{algorithm}

Algorithm~\ref{alg:cutting} starts with the base formulation \textsf{Stair-0} (detailed in \eqref{eq:linear-base} for the additive case and \eqref{eq:formulation-multi} for the multiplicative case). In each iteration, it generates a staircase inequality for each nonlinear term, producing \textsf{Stair-$k$} after $k$ iterations, and finally terminating with \textsf{Stair-K}. Because the base formulation is exact and all cuts are valid, the models \textsf{Stair-0} through \textsf{Stair-K} are increasingly tighter reformulations of PPO.

\begin{algorithm}
\caption{A cutting-plane implementation}\label{alg:cutting}\label{alg:cutting-plane}
\fontsize{10pt}{11pt}\selectfont 
\renewcommand{\baselinestretch}{1.1}\selectfont
\KwData{a base formulation \textsf{Stair-0} given by~\eqref{eq:linear-base} or~\eqref{eq:formulation-multi},  and a positive integer $K$}
\KwResult{formulation \textsf{Stair-K}}
$k \gets 1$\;
\While{$k \le K$}
{
    $(\bar{\x}, \bar{\z}, \bar{\s}, \bar{\boldsymbol{u}}, \bar{\boldsymbol{v}}) \gets $  an optimal solution to the continuous relaxation of $\textsf{Stair-(k-1)}$\;
    
    {\revision $\texttt{cut} \gets$ call Algorithm~\ref{alg:cross-item-effects} and Algorithm~\ref{alg:add-cross-period-effects} for the additive model or Algorithm~\ref{alg:mult-cross-period-effects} for the multiplicative model to separate $(\bar{\z}, \bar{\s}, \bar{\boldsymbol{u}}, \bar{\boldsymbol{v}})$}\;

    $\textsf{Stair-k} \gets $ add $\texttt{cut}$ into formulation $\textsf{Stair-(k-1)}$\;
    $k = k+1$\;
}
\end{algorithm}

Last, we present theoretical guarantees of proposed algorithms.  Let $\textsf{Stair}$ be the final formulation incorporating all staircase inequalities in Algorithms~\ref{alg:cross-item-effects}--\ref{alg:mult-cross-period-effects}. We characterize conditions when the LP relaxation of $\textsf{Stair}$ exactly solves PPO. 
{\revision The key idea of the proof is to invoke a decomposition result, Corollary~\ref{claim:simutanuoushull}, which shows that individual ideal formulations for supermodular compositions yield an ideal formulation for a vector of supermodular compositions.}
\begin{theorem}\label{thm:ppo-lp}
Assume that business rules in PPO are absent and the post-promotion dip effect is present. The LP relaxation of $\textsf{Stair}$ solves PPO when any of the following hold 
\begin{enumerate}
	\item the demand model is as in~\eqref{eq:demand-linear} and all items are substitutable;
	\item the demand model is as in~\eqref{eq:demand-loglog}, all items are substitutable, and the self-effect is non-decreasing. 
\end{enumerate}
\end{theorem}
{ \revision When these business constraints are included and promotion effects are general, the problem may no longer be polynomially solvable. Nevertheless, the proposed formulations remain useful because they provide strong dual bounds that can accelerate the convergence of the branch-and-bound procedure, as illustrated in Section~\ref{sec:formulation-efficiency}.}

\cite{cohen2021promotion} show that the LP relaxation of the $\App(2)$ formulation solves PPO under conditions that are more restrictive than those in Case 1 of Theorem~\ref{thm:ppo-lp}. Specifically, their analysis requires that the number of promotion price levels, $d_i$, is equal to $1$. Example~\ref{ex:mcweak} illustrates why the $\App(2)$ formulation fails to provide an LP reformulation when $d_i > 1$, highlighting the strength of our reformulation in this setting.

\section{Computational Results}\label{section:computation}

This section presents the results of our numerical experiments. In Section~\ref{sec:formulation-efficiency}, we evaluate the performance of our formulations for solving PPO under the log-log  model, with multiplicative historical effects. {\revision In Section~\ref{sec:profit-performance}, we demonstrate that stronger formulations also translate into superior profit performance.} All experiments were conducted on a Mac Studio (16-core M4 Max CPU, 64 GB of memory) with \texttt{Julia 1.11}~\citep{bezanson2017julia}, \texttt{JuMP 1.25}~\citep{lubin2023jump} and \texttt{Gurobi 12}~\citep{gurobi}. We imposed a 3600s time limit and {\revision used \texttt{Gurobi}'s default optimality tolerance of 0.01\%.}

\subsection{Alternative formulations}\label{sec:computation-comparision}
We solve PPO under business rules \eqref{eq:brulesMIP} using four alternative formulations described below: 
\begin{itemize}
    \item $\mathsf{NL}$:  Solves PPO directly using the global optimization solver \texttt{Gurobi 12}. The objective function is expressed as a nonlinear expression {\revision in the space of $\z$ variables} using the demand model from \eqref{eq:demand-loglog}.
    \item $\mathsf{Stair}$-K: For $K = 0$, this corresponds to the base-formulation given by~\eqref{eq:formulation-multi}. For $K \geq 1$, it is generated by Algorithm~\ref{alg:cutting} with  $\mathsf{Stair}$-0 as the initial formulation and $K$ iterations. Each iteration introduces additional constraints to approximate the nonlinear objective.
    \item $\mathsf{Stair}\text{-}\mathsf{Full}$-K: For $K = 0$, this formulation replaces the unary binarization~\eqref{eq:brulesMIP-0} in~\eqref{eq:formulation-multi} with the full binarization scheme. The following constraints are added:
    \[
z^i_{tj} = \sum_{k = j}^{d_i} \lambda^i_{tk} \, \for j \in [d_i],\qquad \boldsymbol{\lambda}^i_t \geq 0, \quad \text{ and }   \sum_{j =0}^{d_i}    \lambda^i_{tj} = 1 . 
    \]
    For $K\geq 1$, it is generated by Algorithm~\ref{alg:cutting-plane} using $\mathsf{Stair}\text{-}\mathsf{Full}$-0 as the base formulation. 
    \item $\mathsf{Stair}\text{-}\mathsf{Log}$-K: For $K = 0$, this formulation replaces the unary binarization~\eqref{eq:brulesMIP-0} in~\eqref{eq:formulation-multi} with the logarithmic formulation, \eqref{eq:log}, of the unary binarization adapted as follows:
    \[
 \gamma^i_{tk} \in \{0,1\}, \quad   \sum_{j \in A^i_{k}}     (z^i_{tj} - z^i_{tj+1}) \leq \gamma^i_{tk},  \text{ and } \sum_{j \notin A^i_{k}} (z^i_{tj} - z^i_{tj+1}) \leq 1-\gamma^i_{tk} \, \for  k  \in \bigl[ \lceil \log_2(d_i+1)\rceil\bigr],
    \]
    where $A^i_{k}=\bigl\{j \in 0 \cup [d_i]  \bigm| \text{$k^{\text{th}}$ coordinate in the binary representation of $j$ is $1$} \bigr\}$, $z^i_{t0} = 1$ and $z^i_{td_i} = 0$. 
    For $K\geq 1$, it is generated by Algorithm~\ref{alg:cutting-plane} using $\mathsf{Stair}\text{-}\mathsf{Log}$-0 as the base formulation. 
\end{itemize}
\subsection{Instance generation}
We follow the instance-generation procedure proposed in~\cite{cohen2021promotion}, which perturbs the parameters of an estimated demand model. Specifically, for each $N \in \{10, 50,75,100,125\}$ and {\revision $T \in \{10,20\}$}, we generate $10$ instances using the following steps.
\begin{enumerate}
    \item Seasonality Effects: The seasonality effects $a_t^i$ are sampled uniformly from the interval $[500, 1000]$.
    \item Price Factors: The current price factors $b_0^i$ are drawn from $[2, 7]$. Each past price factor $b_m^i$ is drawn from $[0, b_0^i]$ with $m\in \{1,2,3\}$ (history window of $3$ periods).
    \item Cross-item Effects: In a typical product category, not all items exhibit cross-item effects. Each pair of products $i$ and $j$ is classified as substitutes with probability $0.1$, complementary with probability $0.1$, and unrelated with probability $0.8$. Conditional on this assignment, the cross-item parameter $\gamma^{ij}$  is drawn from $\frac{1}{5}[0,a^i_t]$ if the items are substitutes, or $\frac{1}{5}[-a_t^i,0]$ if the items are complements. 
    \item Unit Costs and Price Ladders: For each item $i$, the unit cost $c_t^i$ is sampled uniformly from $[0.5, 0.75]$. The price ladder $\{p_0^i, \ldots, p_{4}^i\}$ (with $d_i = 4$) is generated by sampling each price $p_j^i$ uniformly from $[c_t^i, 1]$.
\end{enumerate} 
{\revision To ensure nonnegative demand values, we restrict attention to parameter combinations that yield nonnegative demand in the absence of promotions. Specifically, we randomly generate instances using the parameter ranges described above and discard any instance that violates the nonnegativity condition, replacing it with a newly generated instance.}
\subsection{Formulation efficiency}\label{sec:formulation-efficiency}
In this subsection, we report the computational performance of our formulation $\mathsf{Stair}$-$K$. In particular, we consider $K=5$, since $\mathsf{Stair}$-$5$ achieves the best trade-off between formulation strength and size, as discussed in Appendix~\ref{app:bounds}. We therefore use $\mathsf{Stair}$-$5$ for the subsequent computational experiments.  We evaluate the performance of the formulations using the following metrics: 
\begin{itemize}
    \item \textit{Solved}: the number of solved instances within the time limit
    \item \textit{Time (s)}: Average computation time in seconds
    \item \textit{End Gap (\%)}: The final gap for unsolved instances reported by the solver.
    \item \textit{Nodes}: The number of nodes explored by the solver.
\end{itemize}
As shown in Table~\ref{tab:computation}, $\mathsf{Stair}$-5 consistently outperforms other formulations across all performance metrics. The nonlinear formulation $\mathsf{NL}$ performs the worst and is unable to solve even the smallest instances with 10 products and 10 periods. While the base $\mathsf{Stair}$-0 formulation improves upon $\mathsf{NL}$ and scales better, $\mathsf{Stair}$-5 further enhances its performance by integrating the cutting-plane procedure from Algorithm~\ref{alg:cutting-plane}. This results in solving more instances to optimality, reduced computation time, and smaller end gaps on unsolved problems, while also exploring fewer nodes, demonstrating the cutting-plane approach’s efficiency.

\renewcommand{\arraystretch}{1}
\begin{table}[htbp]
\TABLE 
{{Computational Comparison of Three Alternative Formulations}\label{tab:computation}}
{\begin{tabular}{cc l rrrrrrr}
\toprule
\multicolumn{2}{c}{Instance} & \multirow{2}{*}{Model} & \multirow{2}{*}{Solved} 
& \multicolumn{2}{c}{Time (s)} & \multicolumn{2}{c}{End Gap (\%) } & \multicolumn{2}{c}{Nodes}   \\
\cmidrule(r){1-2} \cmidrule(lr){5-6} \cmidrule(lr){7-8} \cmidrule(lr){9-10}
$N$ & $T$ & & & Mean &Std &  Mean &Std &  Mean & Std  \\
\midrule
\multirow{4}{*}{10} & \multirow{4}{*}{10} 
& $\mathsf{NL}$ & 10 & 3600 & 0 & 15.20 & 3.32 & 10406744 & 4778131   \\
& & $\mathsf{Stair}$-0 & 10 & 0.21 & 0.27 & $\dagger$ & $\dagger$ &17  & 50 \\
& & $\mathsf{Stair}$-5 & 10 & 0.18 & 0.07 & $\dagger$ & $\dagger$ & 1 & 0  \\
\midrule

\multirow{4}{*}{10} & \multirow{4}{*}{20} 
& $\mathsf{NL}$ & 0 & 3600 & 0 & 47.9 & 4.1 & 600591 & 103873   \\
& & $\mathsf{Stair}$-0 & 10 & 0.7 & 0.3 & $\dagger$& $\dagger$ & 512 &  522  \\
& & $\mathsf{Stair}$-5 & 10 & 0.5 &0.1 &  $\dagger$ & $\dagger$ & 69 & 149   \\
\midrule

\multirow{4}{*}{50} & \multirow{4}{*}{10} 
& $\mathsf{NL}$ & 0 & 3600 & 0 &36.7 &6.9 & 577950 & 55821   \\
& & $\mathsf{Stair}$-0 & 10 & 3.1 & 3.4 & $\dagger$ & $\dagger$ & 1824 & 1696  \\
& & $\mathsf{Stair}$-5 & 10 & 2.3 & 1.5 & $\dagger$ & $\dagger$ & 1324 & 3316  \\
\midrule

\multirow{4}{*}{50} & \multirow{4}{*}{20} 
& $\mathsf{NL}$ & 0 & 3600 & 0 & 51.8 & 5.1& 139036 & 12926  \\
& & $\mathsf{Stair}$-0 & 10 & 45 & 40 & $\dagger$& $\dagger$ & 20750 &  29530  \\
& & $\mathsf{Stair}$-5 & 10 & 17 &10 & $\dagger$ & $\dagger$ & 4631 & 2285   \\
\midrule
\multirow{4}{*}{75} & \multirow{4}{*}{10} 
& $\mathsf{NL}$ &0 & 3600 & 0 & 46.9 & 3.6 & 230546 & 29900  \\
& & $\mathsf{Stair}$-0 & 10 & 11.2 & 6.6 & $\dagger$ & $\dagger$ & 3443 & 1777   \\
& & $\mathsf{Stair}$-5 & 10 & 4.5 & 1.4 & $\dagger$ & $\dagger$ & 1197 & 1566 \\
\midrule
\multirow{4}{*}{75} & \multirow{4}{*}{20} 
& $\mathsf{NL}$ &0 & 3600 & 0 & 58.5 &2.6 & 28760 & 9263  \\
& & $\mathsf{Stair}$-0 & 10 & 230 & 207 & $\dagger$ & $\dagger$  & 100918 & 176270 \\
& & $\mathsf{Stair}$-5 & 10 & 117 & 93& $\dagger$ & $\dagger$ & 23118 & 19706 \\
\midrule
\multirow{4}{*}{100} & \multirow{4}{*}{10} 
& $\mathsf{NL}$ & 0 & 3600 & 0 &51.8 & 4.4 &  114974 & 19010  \\
& & $\mathsf{Stair}$-0 & 10 & 36.6 & 42.4 & $\dagger$ & $\dagger$ & 10132 & 14362\\
& & $\mathsf{Stair}$-5 & 10 & 14.1 & 11.8& $\dagger$ & $\dagger$ & 2439 & 2593\\
\midrule
\multirow{4}{*}{100} & \multirow{4}{*}{20} 
& $\mathsf{NL}$ & 0 &3600& 0 & 57.3 & 3.2 & 12416 & 1986  \\
& & $\mathsf{Stair}$-0 & 10 & 288&102 & $\dagger$  & $\dagger$  & 66231 & 42432 \\
& & $\mathsf{Stair}$-5 &10& 135 &86& $\dagger$ & $\dagger$ &  18455 & 17427\\
\midrule
\multirow{4}{*}{125} & \multirow{4}{*}{10} 
& $\mathsf{NL}$ & 0 & 3600& 0 & 56.3 & 2.8 & 61025 & 5361  \\
& & $\mathsf{Stair}$-0 & 10 & 39.3 & 9.2& $\dagger$ & $\dagger$ & 5544 & 3343 \\
& & $\mathsf{Stair}$-5 & 10 & 21.6 & 12.9 & $\dagger$ & $\dagger$ & 2502 & 1609 \\
\midrule
\multirow{4}{*}{125} & \multirow{4}{*}{20} 
& $\mathsf{NL}$ & 0 & 3600 & 0& 58.4 & 3.9 & 10940 & 114 \\
& & $\mathsf{Stair}$-0 & 7 & 1463 & 1500 & 0.13  & 0.02  & 411086 & 526822\\
& & $\mathsf{Stair}$-5 & 10 & 317 & 216 & $\dagger$  & $\dagger$  & 32502 & 22466\\
\bottomrule
\end{tabular}}\\
{\textit{Note}. ``$\dagger$'' indicates that ten instances were solved to optimality. }
\end{table}

We investigate the computational benefits of two formulations: the logarithmic formulation $\mathsf{Stair\text{-}Log}$-K and the full binarization formulation $\mathsf{Stair\text{-}Full}$-K. As detailed in Table~\ref{tab:logformulation}, for instances with five price levels, $\mathsf{Stair\text{-}Log}$-5 shows no computational advantage, likely due to the additional inequalities; however, the logarithmic formulation slightly outperforms the unary formulation, $\mathsf{Stair}$-5 at 30 price levels. Notably, the logarithmic formulation of unary binarization, $\mathsf{Stair\text{-}Log}$-5, consistently outperforms both other formulations across all instances. While root relaxation gaps are similar, reflecting equivalence up to a linear transformation, their formulation sizes and branching directions likely account for the varying computational performance.
\renewcommand{\arraystretch}{0.9}
\begin{sidewaystable}[htbp]
\TABLE 
{{Computational Benefits of Full Binarization and Logarithmic Formulation of Unary Binarization}\label{tab:logformulation}}
{\begin{tabular}{cc l rrrr rrrr rrrr}
\toprule
\multicolumn{2}{c}{Instance} & \multirow{2}{*}{Level} 
& \multicolumn{4}{c}{$\mathsf{Stair}$-5} & \multicolumn{4}{c}{$\mathsf{Stair}$-$\mathsf{Log}$-5} & \multicolumn{4}{c}{$\mathsf{Stair}$-$\mathsf{Full}$-5}   \\
\cmidrule(r){1-2} \cmidrule(lr){4-7}  \cmidrule(lr){8-11}  \cmidrule(lr){12-15} 
$N$ & $T$ & & Sol & Time(s) & Gap(\%) & Nodes & Sol & Time(s) & Gap(\%) & Nodes & Sol & Time(s) & Gap(\%) & Nodes \\
\midrule
\multirow{2}{*}{10} & \multirow{2}{*}{20} 
& 5 & 10 & 0.5 & $\dagger$ & 69 & 10& 0.46 & $\dagger$ & 73 & 10 & 0.60 & $\dagger$& 61 \\
& & 30 & 10 & 16.8 & $\dagger$ & 207 & 10 & 10.0 & $\dagger$ & 299 & 10 & 12.2 & $\dagger$& 492 \\
\midrule
\multirow{2}{*}{50} & \multirow{2}{*}{20} 
& 5 & 10 & 17 & $\dagger$ & 4631 & 10& 14.0 & $\dagger$ & 3597 & 10 & 13.6 & $\dagger$& 4339\\
& & 30 & 7 & 2245.5& 0.35 & 30298 & 10 & 1022.8 & $\dagger$ & 16330 & 9 & 2006.8 & 0.04& 43742 \\
\midrule
\multirow{2}{*}{75} & \multirow{2}{*}{20} 
& 5 & 10 & 117 & $\dagger$  & 23118  & 10 & 100 & $\dagger$ &19264 & 10 & 58 & $\dagger$& 11565  \\
& & 30 & 0 & 3600 & 0.42  & 41712  &4   & 2716.2 & 0.32 & 33367 & 0 & 3600 & 0.16 & 43473  \\
\midrule
\multirow{2}{*}{100} & \multirow{2}{*}{20} 
& 5 & 10 &135 & $\dagger$  & 18455 & 10 & 139 & $\dagger$ &16102 & 10 & 127 & $\dagger$ &13277\\
& & 30 & 0 & 3600 & 0.42 & 26698  & 1 &3500 & 0.24 & 40212 & 0 & 3600 & 0.22 &29429 \\
\midrule
\multirow{2}{*}{125} & \multirow{2}{*}{20} 
& 5 & 10 & 317 & $\dagger$ & 32502 & 10 & 351 & $\dagger$ & 30247 & 10 & 258 & $\dagger$ & 23272\\
& & 30 & 0 & 3600 & 0.57  & 9031  & 1 & 3451 & 0.35 & 21845 & 0 & 3600 & 0.33 & 9900 \\
\bottomrule
\end{tabular}}\\
{\textit{Note}. ``$\dagger$'' indicates that ten instances were solved to optimality. }
\end{sidewaystable}

\subsection{Profit Performance of Promotion Plans}~\label{sec:profit-performance}
{\revision
We evaluate the practical economic value of the promotion plans generated by our $\mathsf{Stair}$ formulation relative to the benchmark approximation $\App$(2) from Cohen et al. (2021). While previous experiments focused on computational efficiency, this subsection quantifies whether the stronger formulations translate into superior profit performance. To measure the baseline value of optimization, we first compute the percentage profit improvement of the $\App$(2) plan over a no-promotion baseline:
\[
\text{Gain}_{\App(2)} = \frac{R_{\textsc{App2}}-R_{\textsc{NoPromo}}}{R_{\textsc{NoPromo}}}\times 100\%,
\]
where $R_{\textsc{App2}}$ and $R_{\textsc{NoPromo}}$ denote the profits from $\App$(2) and the no-promotion baseline, respectively. We then quantify the additional value recovered by our method as the percentage increase over the base gain:
\[
\text{Addl. Gain}_{\mathsf{Stair}} = \frac{R_{\textsc{Stair}}-R_{\textsc{App2}}}
{R_{\textsc{App2}}-R_{\textsc{NoPromo}}}
\times 100\%,
\]
where $R_{\textsc{Stair}}$ is the profit from the $\mathsf{Stair}$ plan. This metric represents the fraction of the potential promotion value (captured by $\App$(2)) that is additionally recovered by solving the exact formulation.

Table~\ref{tab:revgain} reports these metrics across varying levels of cross-item effects and business constraints. Each configuration averages results over 20 instances with 75 items, 7 periods, and 5 promotion levels. We vary the magnitude of cross-item effects by scaling the vector $\gamma$ (where each component corresponds to an item-period pair) by factors $0.8$, $1.0$, and $1.2$. We also vary the business constraints $(K^i, C_t)$, which bound the number of periods an item $i$ can be promoted and the number of items promoted in period $t$, respectively, to assess performance under tight versus loose restrictions.

\renewcommand{\arraystretch}{0.9}
\begin{table}[h]
\centering
\caption{Profit Performance Comparison: $\App$(2) vs. $\mathsf{Stair}$. The $\text{Gain}_{\App(2)}$ column shows the improvement of $\App$(2) over no promotion. The $\text{Addl. Gain}_{\mathsf{Stair}}$ column shows the additional percentage gain achieved by $\mathsf{Stair}$ relative to the $\App$(2) gain. For example, a 11.75\% gain for $\App$(2) and a 29.90\% additional gain for $\mathsf{Stair}$ implies a total profit increase of $11.75 \times (1 + 0.299) \approx 15.26\%$ over no promotion.}
\label{tab:revgain}
\begin{tabular}{c c cc cc}
\toprule
\multirow{2}{*}{Cross-Item Effect } 
  & \multirow{2}{*}{$(K^i,C_t)$}  &  \multicolumn{2}{c}{$\mathsf{App2}$ } & \multicolumn{2}{c}{$\mathsf{Stair}$ } \\
 \cmidrule(lr){3-4}  \cmidrule(lr){5-6}  
& &  Time(s) & $\text{Gain}_{\App(2)}$ &  Time(s) & $\text{Addl. Gain}_{\mathsf{Stair}}$ (\%) \\
\midrule
\multirow{3}{*}{$0.8\gamma$} 
& $(3,15)$         &56.6 & 11.75 & 2.6 & 29.90 \\
& $(5,45)$        &56.7 & 11.72 &2.5 & 46.86 \\
& $(7,75)$ & 56.6 & 11.69 & 2.3 & 48.17 \\
\midrule
\multirow{3}{*}{$\gamma$} 
& $(3,15)$          & 57.1 & 12.05 & 2.6 & 27.84 \\
& $(5,45)$         & 57.0 & 12.25 & 2.5 & 42.31 \\
& $(7,75)$ & 56.9 & 12.06 & 2.6 & 45.18 \\
\midrule
\multirow{3}{*}{$1.2\gamma$} 
& $(3,15)$          & 56.4  &  11.82      &  2.5    &  25.42  \\
& $(5,45)$        & 56.6   &  12.03    &    2.5  &  39.76  \\
& $(7,75)$ & 56.5  &   11.85    &   2.5 &    42.24  \\
\bottomrule
\end{tabular}
\end{table}

Table~\ref{tab:revgain} yields three key observations. First, $\App$(2) consistently generates substantial profit improvements over the no-promotion baseline, with gains ranging from 11.69\% to 12.25\%, confirming the economic value of optimized promotion decisions even under approximation. Second, the $\mathsf{Stair}$ formulation significantly outperforms $\App$(2) in terms of profit recovery. The additional profit value recovered by $\mathsf{Stair}$ ranges from 25.42\% to 48.17\% of the baseline gain achieved by $\App$(2). These improvements are most pronounced under medium and unconstrained business rules (e.g., $(K^i,C_t) = (5,45)$ or $(7,75)$), where the formulation has greater flexibility to exploit cross-item effects. For instance, under low cross-item effects ($0.8\gamma$), $\mathsf{Stair}$ recovers an additional 46.86\% and 48.17\% of the promotion value compared to $\App$(2) in the less constrained settings. Last,  our $\mathsf{Stair}$ formulation solves faster than $\App$(2). This occurs for two reasons. First, our formulation size grows much slower than $\App$(2) with the number of discretization levels. Specifically, the number of variables in our formulation is linear in the number of discretization variables while that in $\App$(2) is quadratic. Second, our formulation captures convex hull structures for supermodular compositions so that the linear relaxations are much tighter.
}

\section{Connections with \texorpdfstring{$\lnatural$}{L-natural}-convex functions}\label{section:discreteconvexity}
The concept of $\lnatural$-convex functions provides a powerful framework for analyzing operations models with indivisibilities \citep{murota1998discrete}. \cite{chen2021discrete} survey a wide range of applications of this concept in the design of operations systems. In particular, \cite{chen2018preservation} leverage $\lnatural$-convexity to address capacity allocation problems in network revenue management, where fixed resource capacities are dynamically allocated across multiple products facing stochastic demand. In this section, we discuss its connection to our formulation of submodular compositions, studied in Section~\ref{section:supermodular}. 

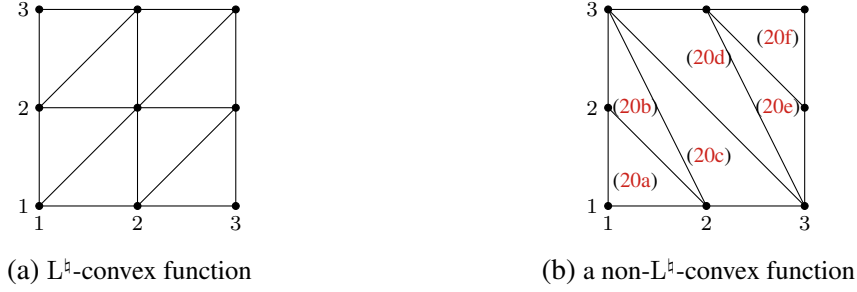
\begin{figure}
    \centering
    \begin{subfigure}[b]{0.45\linewidth}
    \centering
    \tikzsetnextfilename{promotion-9}
\begin{tikzpicture}[scale=1.3]
\begin{scope}
    \node[below] at (0 ,0) {\scriptsize $1$};
    \node[left] at (0, 0) {\scriptsize $1$};
    \node[left] at (0, 1) {\scriptsize $2$};
    \node[below] at (1 ,0) {\scriptsize $2$};
    \node[below] at (2 ,0) {\scriptsize $3$};
    \node[left] at (0 ,2) {\scriptsize $3$};
\filldraw[black] (0,0) circle (1pt) ; 
\filldraw[black] (0,1) circle (1pt) ; 
\filldraw[black] (0,2) circle (1pt) ; 
\filldraw[black] (1,0) circle (1pt) ; 
 \filldraw[black] (1,1) circle (1pt) ; 
\filldraw[black] (1,2) circle (1pt) ; 
\filldraw[black] (2,2) circle (1pt) ; 
\filldraw[black] (2,1) circle (1pt) ; 
\filldraw[black] (2,0) circle (1pt) ; 
\draw[black] (0,0) -- (0,1);
\draw[black] (0,0) -- (1,0);
\draw[black] (0,0) -- (1,1);

\draw[black] (1,0) -- (1,1);
\draw[black] (1,1) -- (0,1);

\draw[black] (0,1) -- (1,2);

\draw[black] (1,2) -- (1,1);

\draw[black] (0,2) -- (0,1);

\draw[black] (2,0) -- (2,1);
\draw[black] (1,1) -- (2,2);

\draw[black] (1,0) -- (2,1);
\draw[black] (1,0) -- (2,0);

\draw[black] (0,2) -- (1,2);
\draw[black] (1,1) -- (2,1);

\draw[black] (2,2) -- (2,1);
\draw[black] (2,2) -- (1,2);

    \end{scope}
\end{tikzpicture} 
    \caption{ \small $\lnatural$-convex function}
    \label{fig:triangulation-Lnatural}
\end{subfigure}
    \begin{subfigure}[b]{0.45\linewidth}
    \centering
    \tikzsetnextfilename{promotion-10}
\begin{tikzpicture}[scale=1.3]
\begin{scope}
    \node[below] at (0 ,0) {\scriptsize $1$};
    \node[left] at (0, 0) {\scriptsize $1$};
    \node[left] at (0, 1) {\scriptsize $2$};
    \node[below] at (1 ,0) {\scriptsize $2$};
    \node[below] at (2 ,0) {\scriptsize $3$};
    \node[left] at (0 ,2) {\scriptsize $3$};

        \node at (0.25 ,0.25) {\scriptsize~\eqref{eq:exlnatural-1}};
        \node at (0.25 ,1) {\scriptsize~\eqref{eq:exlnatural-2}};
        \node at (1 ,0.5) {\scriptsize~\eqref{eq:exlnatural-3}};
        \node at (1 ,1.5) {\scriptsize~\eqref{eq:exlnatural-4}};
        \node at (1.7 ,1) {\scriptsize~\eqref{eq:exlnatural-5}};
        \node at (1.7 ,1.7) {\scriptsize~\eqref{eq:exlnatural-6}};

\filldraw[black] (0,0) circle (1pt) ; 
\filldraw[black] (0,1) circle (1pt) ; 
\filldraw[black] (0,2) circle (1pt) ; 
\filldraw[black] (1,0) circle (1pt) ; 
\filldraw[black] (1,2) circle (1pt) ; 
\filldraw[black] (2,2) circle (1pt) ; 
\filldraw[black] (2,1) circle (1pt) ; 
\filldraw[black] (2,0) circle (1pt) ; 
\draw[black] (0,0) -- (0,1);
\draw[black] (0,0) -- (1,0);
\draw[black] (0,1) -- (1,0);

\draw[black] (1,0) -- (2,0);
\draw[black] (0,2) -- (2,0);
\draw[black] (0,2) -- (1,0);

\draw[black] (0,2) -- (0,1);

\draw[black] (2,0) -- (2,1);
\draw[black] (1,2) -- (2,1);
\draw[black] (1,2) -- (2,0);

\draw[black] (0,2) -- (1,2);

\draw[black] (2,2) -- (2,1);
\draw[black] (2,2) -- (1,2);

    \end{scope}
\end{tikzpicture} 
    \caption{ \small a non-$\lnatural$-convex function}
    \label{fig:triangulation-nonLnatural}
\end{subfigure}
    \caption{Simplicial subdivisions of a box domain associated with the locally Lov\'asz extension and staircase extension}
    \label{fig:triangulation}
\end{figure}

We begin by reviewing definitions of  $\lnatural$ functions. A function $\phi: \Z^n \to \R$ is  \textit{$\lnatural$-convex} if $\phi(x_1 -x_0, x_2 -x_0, \ldots, x_n -x_0 )$ is submodular on $(x_0,\x) \in \Z \times  \Z^n$. Clearly, an $\lnatural$-convex function is a submodular function,  as is seen by setting $x_0=0$, but the converse does not hold.  
One of the key properties of $\lnatural$-convex functions is that they can be extended to convex functions in real variables. Such an extension can be obtained by gluing together Lov\'asz extensions of $\phi(\cdot)$ over sub-boxes as follows. 
For a point $\x \in \R^n$, let $\sigma$ be a permutation of $[n]$ such that $x_{\sigma(1)} - \lfloor x_{\sigma(1)} \rfloor \geq \cdots \geq  x_{\sigma(n)} - \lfloor x_{\sigma(n)} \rfloor $, and consider a simplical neighborhood of $\x$ defined as $\Delta(\x) := \bigl\{ \lfloor \x \rfloor + \sum_{i = 1}^k\e_{\sigma(i)} \bigm| k \in \{0\} \cup [n] \bigr\}$, where $\e_j$ is the $j^\text{th}$ principal vector in $\R^n$. For a function $\phi:\Z^n \to \R$, its \textit{locally Lov\'asz extension} is given as the interpolation of $\phi(\cdot)$ over $\Delta(\x)$: 
\[
\hat{\phi}(\x) :=  \min \Biggl\{ \sum_{\v \in \Delta(\x)} \lambda_{\v} \phi(\v) \Biggm| \sum_{\v \in \Delta(\x)} \lambda_{\v} \v =\x,\ \sum_{\v \in \Delta(\x)} \lambda_{\v} = 1,\ \lambda_{\v} \geq 0 \, \for \v \in \Delta(\x) \Biggr\}.
\]
Figure~\ref{fig:triangulation-Lnatural} provides a geometric illustration of the regular simplical subdivision used to construct the locally Lov\'asz extension of a two-dimensional function. Although the extension function is not convex for all discrete functions, $\lnatural$-convex functions are those functions for which the locally Lov\'asz extension is convex. 
\begin{lemma}[\cite{murota1998discrete}]\label{lemma:murota}
A function $\phi:\Z^n \to \R$ is $\lnatural$-convex if and only if  $\hat{\phi}(\cdot)$ is convex. 

\end{lemma}

Next, we discuss how our formulation for submodular compositions, derived in Theorem~\ref{them:supermodular-model}, is related to $\lnatural$-convexity. For a function $\phi: \prod_{i \in [n]} \{0,1, \ldots d_i\} \to \R^n$ with $d_i \in \Z_+$, we introduce a continuous function $\check{\phi}: \prod_{i \in [n]}[0,d_i] \to \R$ so that, for each $\x$, there exists a $\z$ obtained by lifting $\x$ to $\Delta$ such that $\check{\phi}(\x)$ dominates all staircase inequalities evaluated at $\z$. Formally, we let $\check{\phi}(\x):= \min\bigl\{\mu \bigm| (\mu,\z, \x) \in \Phi \bigr\}$, where $\Phi$ is the LP relaxation of our MIP formulation of the epigraph of $\phi$ given as follows: 
\[
\begin{aligned}
\Phi:=\biggl\{ (\mu,\x,\z) &\biggm| \mu \geq \phi(\j^0_\pi) + \sum_{s  \in [D]} \bigl( \phi(\j^s_\pi) - \phi(\j^{s-1}_\pi) \bigr) \cdot z_{\pi_s}  \for \pi \in \Pi, \\
    & \qquad \qquad \qquad \z \in \Delta,\ x_i = \sum_{j \in [d_i]} z_{ij} \for i \in [n] \biggr\}.
\end{aligned}
\]
This continuous function will be referred to as the \textit{staircase extension} of $\phi(\cdot)$. We note  that the convex function $\check{\phi}(\cdot)$ is always at least as large as the \textit{convex envelope}, since at least one staircase inequality at each $\z$ interpolates the grid points.
The staircase extension differs from the locally Lov\'asz extension in two key ways. First, the staircase extension is convex and coincides with the convex envelope of $\phi(\cdot)$ if and only if $\phi(\cdot)$ is submodular. Second, unlike the locally Lov\'asz extension, the staircase extension $\check{\phi}(\cdot)$ may not always extend $\phi(\cdot)$ from grid points, \textit{i.e.}, it can be above or below $\phi(\cdot)$ at a grid point. For example, for $\phi(\x) = -x_1x_2$ over $\{0,1,2\}^2$, $\check{\phi}(1,1) < \phi(1,1)$. 
Notably, $\check{\phi}(\cdot)$ extends $\phi(\cdot)$ from grid points if and only if $\phi(\cdot)$ is $\lnatural$-convex, in which case it matches the locally Lov\'{a}sz extension.

\begin{theorem}\label{them:Lnatural}
For a function $\phi: \prod_{i \in [n]} \{0,1, \ldots d_i\} \to \R^n$ with $d_i \in \Z_+$, the staircase extension  $\check{\phi}(\cdot)$ is the convex envelope of $\phi(\cdot)$ if and only if $\phi(\cdot)$ is submodular. Moreover, the locally Lov\'asz extension $\hat{\phi}(\cdot)$ coincides with the staircase extension $\check{\phi}(\cdot)$  if and only if $\phi(\cdot)$  is $\lnatural$-convex.
\end{theorem}

When $\phi(\cdot)$ is a submodular composition, but not $\lnatural$-convex, its locally Lov\'asz extension $\hat{\phi}(\cdot)$ is not convex (see Lemma~\ref{lemma:murota}). However, the staircase extension, by Theorem~\ref{them:supermodular-model}, remains the convex envelope of $\phi(\cdot)$. Thus, minimizing $\phi(\cdot)$ can be achieved by minimizing its staircase extension $\check{\phi}(\cdot)$. Minimizing $\hat{\phi}(\cdot)$ would also suffice but does not provide a tractable alternative. Theorem~\ref{thm:ppo-lp} shows how this allows us to handle pricing decisions for substitutable products with different demand models. We now illustrate this difference using a numerical example, demonstrating that distinct simplicial subdivisions arise from the staircase extension and the locally Lov\'asz extension.

\begin{example}\label{ex:MonomialNotLnatural}
Consider $x_1^3x_2^2$ with $x_i \in \{1,2,3\}$ for $i = 1,2$. Using Proposition~\ref{prop:bilinear-sub} and projecting out the $\z$ variables, the convex envelope of $x_1^3x_2^2$ over $[1,3]^2$ is the pointwise maximum of the following affine functions:
\begin{subequations}
    \begin{align}
        &7x_1+3x_2 - 9   \label{eq:exlnatural-1}    \\
         &9x_1 + 5x_2 - 15   \label{eq:exlnatural-2}   \\
        &19x_1 + 10x_2- 40  \label{eq:exlnatural-3}   \\
         &171x_1 + 135x_2 - 675  \label{eq:exlnatural-4}   \\
         &63x_1 +  54 x_2 - 216 \label{eq:exlnatural-5}    \\
         &117x_1 + 81x_2 - 405  \label{eq:exlnatural-6}.   
    \end{align}
\end{subequations}
Each affine function interpolates $x_1^3x_2^2$ over three affinely independent points:  corner points of the labeled region in Figure~\ref{fig:triangulation-nonLnatural}. As is clear, this subdivision is quite different from that in  Figure~\ref{fig:triangulation-Lnatural}. The former is associated with the convex envelope of $x_1^3x_2^2$ and the latter is the simplicial subdivision for the convex envelope of all $\lnatural$ functions in two-dimensions. The function $x_1^3x_2^2$ is not $\lnatural$ and it remains unclear how its convex envelope can be derived without the lifting and projection techniques proposed here.
\hfill \Halmos
\end{example}

\section{Conclusion}
Our findings highlight the potential of combining supermodular function theory with geometric techniques, specifically staircase triangulations, for optimization in discrete domains. The characterization of the convex hull of supermodular compositions significantly advances the treatment of realistic and effective pricing models and the solution of practically-sized pricing models important for supermarket operations. Unlike previous work that relies on simplifying assumptions such as two price levels and additively-separable historical effects on demand for polynomial algorithms and exact formulations respectively, our approach incorporates the practical constraint of discrete price ladders and allows for the consideration of consumer behavior reflected in more sophisticated consumer demand models. 

Our numerical results demonstrate a substantial computational advantage, enabling scalability to 100 products and 20 periods and handling up to 30 price levels. {\revision For the classes of promotion planning problems covered by our polynomial-time results, the proposed approach yields exact solutions via linear programming. More generally, the resulting formulations provide strong relaxations that can improve the computational performance of mixed-integer optimization methods. These results demonstrate the potential of the approach for practically sized promotion planning problems.}

Future work will focus on integrating business rule sets within the convexification scheme, leading to stronger relaxations and improved solution bounds. More broadly, this work provides a foundation for applying the developed techniques to a wider range of problems, including polynomial discrete optimization problems, and pricing optimization under consumer choice models.

\begingroup \parindent 0pt \parskip 0.0ex \def\enotesize{\normalsize} \theendnotes \endgroup

\bibliographystyle{plainnat} 
 
\ECSwitch

\ECDisclaimer

\ECHead{Additional discussions, missing proofs, and formulations}

\section{Proof of Proposition~\ref{prop:simple-exact}}
Using the unary binarization scheme, a lifted bilinear representation of~\eqref{eq:example} is given as follows:
\[
\Bigl\{(\x,\z,\mu) \Bigm| \mu = l_1(\z_1)l_2(\z_2),\ x_i = l_i(\z_i),\ \z_i \in \Delta^{d_i}\cap\{0,1\}^{d_i} \for i = 1,2  \Bigr\},
\]
where $l_i(\z_i) := p_{i0} + \sum_{j \in [d_i]}(p_{ij}-p_{ij-1})z_{ij}$. To complete the proof, we need to show that replacing $\mu = l_1(\z_1)l_2(\z_2)$ with inequalities~\eqref{eq:generalziedMC} yields an exact formulation. Consider a partition $(A,B)$ of $\vertex(\Delta)$ defined as follows:
\[
A =\bigl\{(\v_{10},0), (\v_{10},1), (\v_{11},1) \ldots,  (\v_{1d_1},1) \bigr\} \text{ and } B = \bigl\{(\v_{10},0), (\v_{11},0), \ldots, (\v_{1d_1},0), (\v_{1d_1},1) \bigr\}. 
\]
It can be verified that the inequality~\eqref{eq:generalziedMC1} (resp.~\eqref{eq:generalziedMC2}) is tight over $\bigl(\z_1,\z_2,l_1(\z_1)l_2(\z_2)\bigr)$ for $(\z_1,\z_2) \in  A$ (resp. $(\z_1,\z_2) \in  B$). Moreover, they are valid for $\{(\z,\mu) \mid \mu =l_1(\z_1)l_2(\z_2),\ \z \in  \vertex(\Delta) \}$ since evaluating the right-hand-side (rhs) of~\eqref{eq:generalziedMC1} at the vertex $(\v_{1t},0)$ yields $p_{10}p_{20} + p_{21}(p_{1t} - p_{10})$, which is greater than the function value $l_1(\v_{1t})l_2(0) = p_{1t}p_{20}$, and evaluating the rhs of~\eqref{eq:generalziedMC2} at the vertex $(\v_{1t},1)$  yields $p_{1t}p_{20} + p_{1d_1}(p_{21} - p_{20})$ which is greater than the function value $l_1(\v_{1t})l_2(1) = p_{1t}p_{21}$. Thus,~\eqref{eq:generalziedMC1} and~\eqref{eq:generalziedMC2} model $\mu \leq l_1(\z_1)l_2(\z_2)$ over $\z \in \vertex(\Delta)$. 

To see that inequalities~\eqref{eq:generalziedMC3} and~\eqref{eq:generalziedMC4} model $\mu \geq l_1(\z_1)l_2(\z_2)$, we consider another partition $(C,D)$ of $\vertex(\Delta)$ given as follows:
\[
C = \bigl\{(\v_{10},1), (\v_{10},0), (\v_{11},0) \ldots, (\v_{1d_1},0) \bigr\} \text{ and } D = \bigl\{(\v_{10},1), (\v_{11},1), \ldots, (\v_{1d_1},1) , (\v_{1d_1},0) \bigr\}.
\]
Observe that the inequality~\eqref{eq:generalziedMC3} (resp.~\eqref{eq:generalziedMC4}) is tight over $\bigl\{\bigl(\z_1,\z_2,l_1(\z_1)l_2(\z_2)\bigr)\bigr\}$ for $(\z_1,\z_2) \in  C$ (resp. $(\z_1,\z_2) \in  D$). Moreover, they are valid over $\vertex(\Delta)$ as evaluating the rhs of~\eqref{eq:generalziedMC3} at the vertex $(\v_{1t},1)$ yields $p_{10}p_{21} + p_{20}(p_{1t}-p_{10})$, which is less than the function value $p_{1t}p_{21}$, and evaluating the rhs of~\eqref{eq:generalziedMC4} at the vertex $(\v_{1t},0)$ yields $p_{1t}p_{21} + p_{1d_1}(p_{21}-p_{20})$, which is less than the function value $p_{1t}p_{20}$.      \hfill \Halmos

\section{Proof of Proposition~\ref{prop:simple-perfect}}
Let $E:=\bigl\{(\x,\mu,\z) \in  \R^2 \times \R \times \{0,1\}^{d_1+d_2} \bigm|\eqref{eq:unary-example}\text{ and }\eqref{eq:hull-example}\bigr\}$. Due to Proposition~\ref{prop:simple-exact}, to prove  the validity of $E$, it suffices to show that inequalities~\eqref{eq:hull-example} are valid for $\mu = l_1(\z_1)l_2(\z_2)$ over the vertices of $\Delta$. For $k \in 0 \cup [d_1]$, let $\psi_k(\z)$ denote the rhs of~\eqref{eq:hull-example-1} and observe that 
\[
\begin{aligned}
	\psi_k(\v_{1t},0) &= p_{10}p_{20} + p_{20}(p_{1 \min\{k,t\}} - p_{10}) +p_{21}( p_{1\max \{k,t\}} -p_{1k})  \\
	& = p_{20}p_{1 \min\{k,t\}} + p_{21}(p_{1\max \{k,t\}}- p_{1k}) \\
	& \geq p_{20}(p_{1 \min\{k,t\}} + p_{1\max \{k,t\}}- p_{1k}) = p_{1t}p_{20} = l_1(\v_{1t})l_2(0) \qquad \for t \in 0 \cup [d_i],\\
	\psi_k(\v_{1t},1) &= p_{10}p_{20} + p_{20}(p_{1 \min\{k,t\}} - p_{10})+ p_{1k}(p_{21}-p_{20}) +p_{21}( p_{1\max \{k,t\}} -p_{1k}) \\
	& = p_{1 \min\{k,t\}}p_{20} - p_{1k}p_{20} + p_{1\max\{k,t\}}p_{21} \\
	& \geq p_{1 \min\{k,t\}}p_{21} - p_{1k}p_{21} + p_{1\max\{k,t\}}p_{21} = p_{1t}p_{21} = l_1(\v_{1t})l_2(1) \qquad \for t \in 0 \cup [d_i].
\end{aligned}
\]
For $k \in 0 \cup [d_1]$, let $\varphi_k(\z)$ denote the rhs of~\eqref{eq:hull-example-2} and observe that 
\[
\begin{aligned}
\varphi_k(\v_{1t},1) &= p_{10}p_{21} + p_{21}(p_{1 \min\{k,t\}} -p_{10})+p_{20}( p_{1\max \{k,t\}} -p_{1k}) \\
& = p_{21}p_{1 \min\{k,t\}} +p_{20}( p_{1\max \{k,t\}} -p_{1k}) \\
& \leq p_{21}(p_{1 \min\{k,t\}} + p_{1\max \{k,t\}}- p_{1k}) = p_{1t}p_{21} = l_1(\v_{1t})l_2(1) \qquad \for t \in 0 \cup [d_i],\\
\varphi_k(\v_{1t},0) &= p_{10}p_{21} + p_{21}(p_{1 \min\{k,t\}} -p_{10})+ p_{1k}(p_{20}-p_{21})+p_{20}( p_{1\max \{k,t\}} -p_{1k}) \\
& = p_{1\min\{k,t\}} p_{21} - p_{1k}p_{21} + p_{1\max\{k,t\}}p_{20} \\
& \leq p_{1\min\{k,t\}} p_{20} - p_{1k}p_{20} + p_{1\max\{k,t\}}p_{20} =p_{1t}p_{20} = l_1(\v_{1t})l_2(0) \qquad \for t \in 0 \cup [d_i].\end{aligned}
\]

Next, we show that $E$ is ideal. To do so, we consider the LP relaxation $R$ of $E$, and consider a vertex $(\x,\mu,\z)$ of $R$. Suppose that $\z$ is not binary, that is, there exists a pair $(i,j)$ such that $z_{ij}$ is not binary. We will argue that this implies that $(\x,\mu,\z)$ is not an extreme point, a contradiction. Since $(\x,\mu,\z)$ is an extreme point of $R$, either $\mu = \min\bigl\{\psi_k(\z)\bigm| k \in 0 \cup [d_1]\bigr\}$ or  $\mu = \max\bigl\{\varphi_k(\z)\bigm| k \in 0 \cup [d_1]\bigr\}$. Let us consider the first case. Let $k^* \in 0\cup [d_i]$ such that $z_{1k^*} \geq z_{21} \geq z_{1k^*+1}$, where $z_{10}=1$ and $z_{1d_1+1}=0$. Consider a subset of $d_1+2$ vertices of $\Delta$ given as 
\[
\w_j:=(\v_{1j},0)\, \for j \in \{0,1,\ldots, k^*\} \quad \text{ and } \quad  \boldsymbol{m}_j:=(\v_{1j},1) \,  \for j \in \{k^*, \ldots, d_1 \}.
\] 
and observe that $\z$ can be expressed as a convex combination of these points, that is, $\z = \sum_{j =0}^{k^*}\alpha_j\w_j + \sum_{j =k^*}^{d_1}\beta_j\boldsymbol{m}_j$, where  $\alpha_0 = 1-z_{11}$,  $\alpha_j=z_{1j}-z_{1(j+1)}$ for $j \in \{1, \ldots, k^*-1\}$ and $\alpha_{k^*}=z_{1k^*} - z_{21}$, and $\beta_{k^*} = z_{21} - z_{1(k^*+1)}$, $\beta_{j} = z_{1j} - z_{1(j+1)}$ for $j \in \{k^*+1, \ldots, d_1-1\}$, and $\beta_{d_1} =z_{1d_1}$. 
Moreover, $\psi_{k^*}(\cdot)$ achieves function values over these points, that is, 
\[
\begin{aligned}
\psi_{k^*} (\v_{1j},0) &= p_{10}p_{20} + \sum_{t = 1}^jp_{20}(p_{1t} - p_{1(t-1)}) = l_1(\v_{1j})l_2(0) \quad \for j \leq k^* \\
\psi_{k^*} (\v_{1j},1) &= p_{1k^*}p_{20} +p_{1k^*}(p_{21} - p_{20}) + \sum_{t = k^*+1}^jp_{21}(p_{1t} - p_{1(t-1)}) = l_1(\v_{1j})l_2(1) \quad \for j \geq k^*.
\end{aligned}
\]
Then, we obtain 
\[
\psi_{k^*}(\z) =  \sum_{j =0}^{k^*}\alpha_jl_1(\v_{1j})l_2(0)  + \sum_{j =k^*}^{d_1}\beta_jl_1(\v_{1j})l_2(1) \leq \min\bigl\{\psi_k(\z)\bigm| k \in 0 \cup [d_1]\bigr\} = \mu \leq  \psi_{k^*}(\z),
\]
where the first equality follows from the linearity of $\psi_{k^*}(\cdot)$ and the convex decomposition of $\z$, and the first inequality holds since the point-wise minimum function is a concave overestimator of $l_1(\z_1)l_2(\z_2)$. Therefore, we can conclude that $(\mu,\z)$ can be expressed as a convex combination of points in the graph of $l_1(\cdot)l_2(\cdot)$. This, together with $\x_i = l_i(\z_i)$ for $i =1,2$, shows that $(\x,\mu,\z)$ can be expressed as a convex combination of points in $R$. In other words, it is not a vertex of $R$. 

To treat the case $\mu = \max\{\varphi_k(\z) \mid k \in 0 \cup [d_1]\}$, we find $k_* \in 0 \cup [d_i]$ such that $z_{1k_*} \geq (1-z_{21}) \geq z_{1k_*+1}$. Consider a subset of $d_1+2$ vertices of $\Delta$ given as 
\[
(\v_{1j},1)\, \for j \in \{0,1,\ldots, k^*\} \quad \text{ and } \quad  (\v_{1j},0) \,  \for j \in \{k^*, \ldots, d_1 \}.
\] 
It can be verified that $\z$ can be expressed as a convex combination of these points, that is, $\z = \sum_{j=0}^{k_*} \alpha_j(\v_{1j},1) + \sum_{j = k_*}^{d_1}\beta_j(\v_{1j},0)$, where $\alpha_0 = 1-z_{11}$,  $\alpha_j=z_{1j}-z_{1(j+1)}$ for $j \in \{1, \ldots, k_*-1\}$ and $\alpha_{k_*}=z_{1k_*} -(1- z_{21})$, and $\beta_{k_*} = (1- z_{21}) - z_{1(k_*+1)}$, $\beta_{j} = z_{1j} - z_{1(j+1)}$ for $j \in \{k_*+1, \ldots, d_1-1\}$, and $\beta_{d_1} =z_{1d_1}$. Moreover, $\psi_{k_*}(\cdot)$ achieves function values over these points, that is, 
\[
\begin{aligned}
\varphi_{k_*} (\v_{1j},1) &= p_{1j}p_{21} = l_1(\v_{1j})l_2(1) \quad \for j \leq k_* \\
\varphi_{k^*} (\v_{1j},0) &= p_{1k_*}p_{21} +p_{1k_*}(p_{20} - p_{21}) + \sum_{t = k_*+1 }^{j}p_{20}(p_{1j}-p_{1(j-1)}) = l_1(\v_{1j})l_2(0) \quad \for j \geq k_*.
\end{aligned}
\]
Then, we obtain 
\[
\varphi_{k_*}(\z) = \sum_{j=0}^{k_*} \alpha_jl_1(\v_{1j})l_2(1) + \sum_{j = k_*+1}^{d_1}\beta_jl_1(\v_{1j})l_2(1) \geq \max \bigl\{\varphi_k(\z) \bigm| k \in 0 \cup[d_1]\bigr\} =\mu \geq \varphi_{k_*}(\z),
\]
where the first equality holds by the linearity of $\varphi_{k_*}(\cdot)$, and the first inequality holds because the point-wise maximum function is a convex underestimator of $l_1(\cdot)l_2(\cdot)$. Therefore, we can conclude that $(\x, \mu,\z)$ can be expressed as a convex combination of points in $\R$.  This shows that $(\x,\mu,\z)$ is not a vertex of $R$. \hfill \Halmos

\section{Proof of Lemma~\ref{lemma:MIP-discrete}}
Here, we prove a generalization of Lemma~\ref{lemma:MIP-discrete}, which considers  a vector of composite functions $\phi \mcirc f : \X \to \R^m$ defined as 
            \[
    (\phi \mcirc f)(\x) = \Bigl(\phi_1\bigl(f_{11}(x_1), 
    \ldots, f_{1n}(x_n) \bigr), \ldots, \phi_m\bigl(f_{m1}(x_1), 
    \ldots, f_{mn}(x_n) \bigr)\Bigr) \qquad \for \x \in \X \cap P.
    \]
Let $\ell(\cdot)$ be the lift of $f(\cdot)$, that is, $\ell_{qi}(\z_i):= f_{qi}(p_{i0}) + \sum_{j \in [d_i]}\bigl(f_{qi}(p_{ij}) - f_{qi}(p_{ij-1})\bigr)z_{ij}$. 
\begin{lemma}\label{eclemma:MIP-discrete}
Let $ \Phi$ be a polytope such that $\Phi \cap \bigl(\{0,1\}^{\sum_{i\in [n]}d_i} \times \R^m \bigr)  = \graph(\phi \mcirc \ell)$.
Then, an MIP formulation of the graph of $\phi \mcirc f$ is given as follows
\begin{equation*}\label{eq:MIP-discrete}
\begin{aligned}
E:=\Bigl\{(\x,\z,\bmu)  \Bigm| (\z,\bmu) \in  \Phi,\ \x \in P ,\ (x_i,\z_i)  \in  \eqref{eq:unary} \; \for i \in [n] \Bigr\}, 
\end{aligned}
\end{equation*}
which is ideal when $\Phi = \conv \bigl(\graph(\phi \mcirc \ell)\bigr)$ and the constraint $\x\in P$ is relaxed.
\end{lemma}
The validity of $E$ follows by observing 
\[
\begin{aligned}
\proj_{(\x, \bmu)}(E) &=  \bigl\{(\x,  \bmu) \bigm| \exists \z  \text{ s.t. } (\z,\bmu) \in \Phi ,\ \x \in P,\  (x_i,\z_i) \in~\eqref{eq:unary} \; \for i \in [n] \bigr\} \\
 &= \bigl\{(\x, \bmu) \bigm| \exists \z \text{ s.t. }  \bmu = (\phi \mcirc \ell)(\z),\  \x \in P,\ (x_i,\z_i) \in~\eqref{eq:unary}\; \for i \in [n] \bigr\} \\
 & =  \bigl\{(\x,\bmu) \bigm| \bmu = (\phi \mcirc f)(\x),\ \x \in P \cap \X\bigr\},
 \end{aligned}
\]
where the first equality holds by definition, the second equality holds because if $\z$ is binary then $(\z, \bmu) \in \Phi$ if and only if  $(\z,\bmu ) \in \graph(\phi \mcirc \ell )$, and the last equality follows because the unary binarization scheme is an exact formulation of the discrete set $P \cap \X$ and $\ell_i(\cdot)$ is a lift of $f_i(\cdot)$. 

To show that $E$ is ideal when $\Phi = \conv(\phi \mcirc \ell )$  and $P = \R^n$, we consider its continuous  relaxation $R = \bigl\{(\x,\z,\bmu)\bigm| (\z, \bmu) \in  \conv\bigl(\graph(\phi \mcirc \ell)\bigr),\ x_i = p_{i0} + \sum_{j \in [d_i]}(p_{ij}-p_{ij-1}) z_{ij}  \for i \in [n] \bigr\}$ and observe that $R$ is a linear transformation of $\conv\bigl(\graph(\phi \mcirc \ell)\bigr)$. Since convexification commutes with affine transformation, it follows that 
\[
R = \conv\biggl( \Bigl\{(\x,\z,\bmu)\biggm| (\z, \bmu) \in  \graph(\phi \mcirc \ell),x_i = p_{i0} + \sum_{j \in [d_i]}(p_{ij}-p_{ij-1}) z_{ij}  \for i \in [n] \Bigr\} \biggr).
\]
Therefore, $\vertex(R) \subseteq \bigl\{(\x,\z,\mu)\bigm| (\z, \bmu) \in  \graph(\phi \mcirc \ell),\ x_i = p_{i0} + \sum_{j \in [d_i]}(p_{ij}-p_{ij-1}) z_{ij}  \for i \in [n] \bigr\}$. In other words, for every $(\x,\z,\bmu) \in \vertex(R)$, $\z \in \vertex(\Delta) \subseteq \{0,1\}^{\sum_{i\in [n]}d_i+1}$. Hence, $E$ is ideal. \hfill \Halmos

\section{Proof of Theorem~\ref{them:supermodular-model}}
\begin{lemma}\label{lemma:switchtosupermodular}
Consider a supermodular composition $\phi \mcirc f: \X \to \R$. Let $\sigma = (\sigma_1, \ldots, \sigma_n)$ satisfy~\eqref{eq:IO}. Then, the function $\phi \mcirc \ell^\sigma: \vertex(\Delta ) \to \R$ is supermodular. 
\end{lemma}

\begin{lemma}\label{lemma:staircase}
    Consider a function $f: \vertex(\Delta) \to \R$, where $\Delta = \Delta^{d_1} \times \cdots \times \Delta^{d_n}$, and for $\j \in \G:=\prod_{i \in [n] }\{0,1, \ldots,d_i\}$, let $\alpha(\j) = f(\v_{1j_1}, \ldots, \v_{nj_n})$,  where $\v_{ij}$ is the $j^{\text{th}}$ vertex of $\Delta^{d_i}$ defined as in~\eqref{eq:verticesDelta}.   Then, if $f(\cdot)$ is supermodular,  the convex hull of the hypograph of $f$ can be described by the following system of valid inequalities:
    \begin{equation}\label{ec:staircaseineq}
        \mu \leq f^\pi(\z):= \alpha(\j^0) + \sum_{s  \in [D]} \bigl( \alpha(\j^s) - \alpha(\j^{s-1}) \bigr) \cdot z_{\pi_s} \qquad \for \pi \in \Pi,
    \end{equation}
    where $(\j^0, \ldots, \j^D)$ be the point representation of a staircase $\pi$ given in~\eqref{eq:pointrep}. Moreover, if $f(\cdot)$ is submodular then the convex hull of epigraph of $f$ is $\bigl\{(\z,\mu) \in \Delta \times \R \bigm| \mu \geq f^\pi(\z) \for \pi \in \Pi \bigr\}$.
\end{lemma}
\noindent{\bf Proof of~Lemma \ref{lemma:staircase}}.
Note that it suffices to treat the supermodular case since $(\z,\mu)$ belongs to the convex hull of the epigraph of $f$ if and only if $(\z, -\mu)$ belongs to the convex hull of the hypograph of $-f$, and $f$ is submodular if  $-f$ is supermodular.

 We start with showing that $\conv(\hypo(f)) \subseteq R:=\{(\z,\mu) \in \Delta \times \R \mid \mu \leq f^{\pi}(\z) \for \pi \in \Pi\}$. To prove this, it suffices to show that $f(\z) \leq f^\pi(\z)$ for every $\z \in \vertex(\Delta)$, showing $R$ is a convex relaxation of  $\hypo(f)$, and thus $\conv(\hypo(f)) \subseteq R$.   For a given staircase $\pi \in \Pi$, let $(\j^0, \ldots, \j^D)$ be its point representation given in~\eqref{eq:pointrep} and let $\v^t$ denote the vertex associated with $\j^t$, that is $\v^t = (\v_{1j^t_1}, \ldots \v_{nj^t_n} )$.  First, we telescope $f$ using the staircase, that is, 
 \[
 f (\z) = f(\z \wedge \v^0) + \sum_{s \in [D]} f(\z \wedge\v^{s}) - f(\z \wedge \v^{s-1}) \quad \for \z \in \vertex(\Delta).
 \]
Then, we relax the difference terms associated with $s^{\text{th}}$ step as follows. Recall that for each step $s \in [D]$, $\pi_s(1)$ denote the moving direction of the $s^{\text{th}}$ step and $\pi_s(2)$ denote the number of steps already taken along that direction (including the current step). Thus $\v^s_i = \v^{s-1}_i$ if  $i \neq \pi_s(1)$ and $\v^s_i = \v^{s-1}_i + \e_{i \pi_s(2)}$ if $i = \pi_s(2)$, where $\e_{ij}$ is the $j^{\text{th}}$ principal vector in the space of $\z_i$ variables. This, together with the supermodularity of $f$, implies that for every $\z \in \vertex(\Delta)$:
\[
\begin{aligned}
f(\z \wedge\v^{s})  - f(\z \wedge \v^{s-1}) & \leq f\bigl(\v^s_1, \ldots, \v^s_{\pi_s(1)-1},  \z_{\pi_s(1)}\wedge \v^s_{\pi_s(1)} , \v^s_{\pi_s(1)+1}, \ldots, \v^s_n \bigr) \\
&  \qquad - f\bigl(\v^{s-1}_1, \ldots, \v^{s-1}_{\pi_s(1)-1},  \z_{\pi_s(1)} \wedge \v^{s-1}_{\pi_s(1)}, \v^{s-1}_{\pi_s(1)+1}, \ldots, \v^{s-1}_n\bigr). 
\end{aligned}
\]
Let $h^\pi_s(\z_{\pi_s(1)})$ denote the difference function on the right-hand-side. Observe that this difference function takes two values, zero and $f(\v^s) - f(\v^{s-1})$, and this incremental change can be activated using the variable  $z_{\pi_s(1), \pi_s(2)}$, that is, $h^\pi_s(\z_{\pi_s(1)}) = (f(\v^s) - f(\v^{s-1})) z_{\pi_s(1), \pi_s(2)}$. Therefore, we conclude that 
\[
f(\z) \leq f(\v^0) + \sum_{s\in [D]}\bigl(f(\v^s) - f(\v^{s-1})\bigr)z_{\pi_s(1),\pi_s(2)} \quad \for \z \in \vertex(\Delta),
\]
where the function on the right-hand-side is $f^\pi(\z)$. 

Now, we prove that $R \subseteq \conv(\hypo(f))$.  Given a staircase $\pi \in \Pi$, we associate it with a subset  of vertices $\{ \v^0,\v^1, \ldots, \v^D \}$ of $\Delta$ as above, and let $\Upsilon_\pi$ be a simplex defined as follows:
\begin{equation}\label{staircase-triangulation}
    \Upsilon_\pi:= \conv( \v^0,\v^1, \ldots, \v^D)
\end{equation}
 It follows from Section 2.2 in~\citet{he2022tractable} or Theorem 6.2.13 in~\citet{de2010triangulations} that $\{\Upsilon_\pi\}_{\pi \in \Pi}$ is a triangulation of $\Delta$, \textit{e.g.}, $\Delta = \cup_{\pi\in \Pi} \Upsilon_{\pi}$ and for $\pi', \pi'' \in \Pi$ $\Upsilon_{\pi'} \cap \Upsilon_{\pi''}$  is a face of both $\Upsilon_{\pi'}$ and $\Upsilon_{\pi''}$. For a given $(\bar{\z}, \bar{\mu}) \in R$, we will express it as a convex combination of points in $\hypo(f)$. Since $\{\Upsilon_\pi\}_{\pi \in \Pi}$ is  a triangulation of $\Delta$, there exists $\omega \in \Pi$ such that $\bar{\z} \in \Upsilon_\omega$. Observe that 
\[
f^\omega(\v) = f(\v) \qquad \for \v \in \vertex(\Upsilon_\omega).
\]
It follows readily that $\bar{\z}$ can be expressed as a convex combination of vertices of $\Upsilon_\omega$, that is, there exist a convex multiplier $\{\lambda_{\v}\}$ such that  $\bar{\z} = \sum_{\v \in \vertex(\Upsilon_\omega)} \lambda_{\v} \v$, and in addition,
\[
 f^\omega(\bar{\z})  = \sum_{\v \in \vertex(\Upsilon_\omega)} \lambda_{\v}f^\omega(\v) = \sum_{\v \in \vertex(\Upsilon_\omega)} \lambda_{\v}f(\v), 
\]
where the first equality holds since $f^\omega$ is an affine function, and the second equality holds since $f^\omega(\v)=f(\v)$ for $\v \in \vertex(\Upsilon_\omega)$. In other words, $(\bar{\z}, f^\omega(\bar{\z}))$ can be expressed as a convex combination of points in $\hypo(f)$. Now, let $\epsilon =  f^\omega(\bar{\z}) -\bar{\mu}$ and observe that $\epsilon \geq 0$ as $\bar{\mu} \leq \min_{\pi \in \Pi} f^\pi(\bar{\z}) \leq  f^\omega(\bar{\z})$. Hence, we obtain
\[
(\bar{\z},\bar{\mu}) = \bigl(\bar{\z}, f^\omega(\bar{\z})-\epsilon\bigr) = \sum_{\v \in \vertex(\Upsilon_\omega)}\lambda_{\v}\bigl(\v, f(\v)-\epsilon\bigr) \in \conv(\hypo(f)),
\]
where the second equality holds since $\{\lambda_{\v}\}_{\v}$ is a convex multiplier to express $\bigl(\bar{\z},f^\omega(\bar{\z}) \bigr)$, and the last relation holds since $(\v, f(\v)-\epsilon)$ belongs to $\hypo(f)$.  \hfill \Halmos

\noindent{\bf Proof of Theorem~\ref{them:supermodular-model}:} Due to Lemma~\ref{lemma:MIP-discrete}, it suffices to show that 
\[
\conv\bigl(\hypo(\phi \mcirc \ell)\bigr) = \bigl\{(\z,\mu) \in \Delta \times \R \bigm| \z'= (L^\sigma)^{-1}(\z),\ \eqref{eq:staircase} \bigr\}.
\]
Given a permutation $\sigma_i$ of $\{0,1, \ldots, d_i\}$, we have  $\ell_i(\z_i) = \ell^{\sigma_{i}}_i\bigl(L^{\sigma_i} (\z_i)\bigr)$ for every $\z_i \in \vertex(\Delta^{d_i})$. It follows readily that $(\phi \mcirc \ell)(\z) = (\phi \mcirc \ell^\sigma)\bigl(L^\sigma(\z)\bigr)$. The hypograph of $\phi \mcirc \ell$ can be obtained as an affine transformation of that of $\phi \mcirc \ell^\sigma$, that is, 
\[
\hypo(\phi \mcirc \ell) = \Bigl\{ (\z, \mu) \Bigm| \z = (L^\sigma)^{-1}(\z'), (\z',\mu) \in \hypo\bigl(\phi \mcirc \ell^\sigma\bigr) \Bigr\}.
\]
The proof is complete by observing that 
\[
\begin{aligned}
\conv\bigl(\hypo(\phi \mcirc f)\bigr) &=  \Bigl\{(\z,\mu) \in \Delta \times \R \Bigm| \z=(L^\sigma)^{-1}(\z'),\ (\z',\mu) \in \conv\bigl(\hypo(\phi \mcirc \ell^\sigma)\bigr) \Bigr\} \\
&= \bigl\{(\z,\mu) \in \Delta \times \R \bigm| \z'=L^\sigma(\z),\ \eqref{eq:staircase} \bigr\},   
\end{aligned}
\]
where the first equality holds since convexification commutes with affine transformation, and the second eqaulity holds since  by Lemma~\ref{lemma:switchtosupermodular},  $\phi \mcirc \ell^\sigma(\z)$ is supermodular, and then invoking Lemma~\ref{lemma:staircase}, we obtain $\conv\bigl(\hypo(\phi \mcirc \ell^\sigma)\bigr) = \bigl\{(\z',\mu) \in \Delta \times \R \bigm|~\eqref{eq:staircase}\bigr\}$. \hfill \Halmos

\section{Proof of Proposition~\ref{prop:bilinear-sub}}
Let $\ell_i(\cdot)$ be the lifted representation of $f_i(\cdot)$. 
By Lemma~\ref{lemma:MIP-discrete}, it suffices to derive the convex hull of the graph of $\ell_1\ell_2$ over $\vertex(\Delta)$, denoted as $\graph(\ell_1\ell_2)$. Moreover, the convex hull of $\graph(\ell_1\ell_2)$ is $ \conv\bigl(\epi(\ell_1\ell_2)\bigr) \cap \conv\bigl(\hypo(\ell_1\ell_2)\bigr)$, where $\epi(\ell_1\ell_2)$ (resp. $\hypo(\ell_1\ell_2)$) is the epigraph (resp. hypograph) of $\ell_1 \ell_2$. Let $\sigma = (\sigma_1,\sigma_2)$ satisfying~\eqref{eq:IO}.  Since the product of two variables is supermodular, by Lemma~\ref{lemma:switchtosupermodular}, $\ell_1^{\sigma_1}\ell_2^{\sigma_2}$ is supermodular over $\vertex(\Delta)$. Following the line of arguments in the proof of Theorem~\ref{them:supermodular-model}, 
\[
\conv\bigl(\hypo(\ell_1\ell_2)\bigr) = \bigl\{(\z,\mu) \bigm| \z' = L^{\sigma}(\z),\ \eqref{eq:staircase}\bigr\}.
\]
For $\varsigma = (\varsigma_1, \varsigma_2)$ satisfying~\eqref{eq:mix-order},  $\ell_1^{\varsigma_1}\ell_2^{\varsigma_2}$ is submodular over $\vertex(\Delta)$. Thus, using  Lemma~\ref{lemma:staircase} and the line of arguments in the proof of Theorem~\ref{them:supermodular-model}, we obtain 
\[
 \conv\bigl(\epi(\ell_1\ell_2)\bigr)   =\biggl\{(\z,\mu) \biggm| \z'' = L^{\varsigma}(\z),\  \mu \geq		\psi^\varsigma(\j^0_\pi) + \sum_{t  \in [D]} \bigl( \psi^\varsigma(\j^t_\pi) - \psi^\varsigma(\j^{t-1}_\pi) \bigr) \cdot z''_{\pi_t}  \for \pi \in \Pi\biggr\}.
\]
This completes the proof. \hfill \Halmos

\section{Proof of Lemma~\ref{lemma:MIP-discrete-log}}
Let $E$ denote the formulation given in Lemma~\ref{lemma:MIP-discrete-log}, that is,
\[
\begin{aligned}
E:=\Bigl\{(\x, \z, \mu, \bdelta) \Bigm| \; (\z,\mu) &\in \conv\bigl(\graph(\phi \mcirc \ell)\bigr),\ \x \in P, \\
& x_i = \pt_{i0} + \sum_{ j \in [d_i]} (\pt_{ij} - \pt_{ij-1}) z_{ij},\ (\z_i, \bdelta_i) \in Q_i,\ \bdelta_i \in \{0,1\}^{r_i}\, \for i \in [n] \Bigr\}.
\end{aligned}
\]
Its validity follows from a similar argument as the one given in the proof of Lemma~\ref{lemma:MIP-discrete}. In the following, we will focus on showing that $E$ is ideal. To establish the ideality of $E$, we invoke a decomposition result, which has been used in the literature~\citep{schrijver1983short,kim2025reciprocity,he2024mip}. 
\begin{lemma}\label{lemma:commonsimplex}
For $i  = 1,2$, let $D_i$ be a subset in the space of $(x_i,y)\in \R^{n_i+m}$,  and let $D :=\{(x,y) \mid (x_i,y) \in D_i \; i = 1,2 \} $. If $\proj_{y}(D_1) = \proj_{y}(D_2)$ is a finite set of affinely independent points then $\conv(D) = \{(x,y) \mid (x_i,y) \in \conv(D_i) \;  i =1,2\}$.
\end{lemma}
\noindent{\bf Proof of Lemma~\ref{lemma:MIP-discrete-log}:} Let $R$ be the LP relaxation of $E$. Consider a recursive definition of $R$  given as follows. Let 
\[
R^0 :=\Bigl\{(\x, \z, \mu) \Bigm| \; (\z,\mu) \in \conv\bigl(\graph(\phi \mcirc \ell)\bigr), 
x_i = \pt_{i0} + \sum_{ j \in [d_i]} (\pt_{ij} - \pt_{ij-1}) z_{ij} \for i \in [n] \Bigr\}, 
\]
and for $t \in [n]$, let $\y^t:=(\bdelta_1, \ldots, \bdelta_t )$  and define 
\[
R^t : = \bigl\{(\x,\z,\mu,\y^t) \bigm| (\z_t,\boldsymbol{\delta}_t) \in Q_t,\ (\x,\z,\mu,\y^{t-1})  \in R^{t-1}  \bigr\}.
\]
Clearly, $R^n = R$. For $t \in 0 \cup [n]$, we  define 
 \[
 E^t := \Bigl\{(\x,\z,\mu,\y^t) \in R^t \Bigm| \z \in \{0,1\}^{\sum_{i\in [n]}d_i},\ \y^t \in \{0,1\}^{\sum_{i \in [t]}r_i}  \Bigr\}.
 \]
Next, we will recursively invoke Lemma~\ref{lemma:commonsimplex} to argue that $R^n = \conv(E^n)$. This implies that $\vertex(R) = \vertex(R^n)  \subseteq E^n$, proving the ideality of $E$. First,
since $R^0$ is a linear transformation of $\conv(\graph(\phi \mcirc \ell))$ and convexification commutes with linear transformation, we have $R^0 = \conv(E^0)$. For $t \in [n]$, suppose that $R^{t-1} = \conv(E^{t-1})$. Then, we obtain
\[
\begin{aligned}
R^t &=  \bigl\{(\x,\z,\mu,\y^t) \bigm|  (\x,\z,\mu,\y^{t-1})  \in R^{t-1} ,\ (\z_t,\bdelta_t) \in Q_t \bigr\} \\
& = \Bigl\{(\x,\z,\mu, \y^t) \Bigm| (\x,\z,\mu, \y^{t-1}) \in \conv(E^{t-1}),\ (\z_t,\bdelta_t) \in \conv(Q_t \cap \{0,1\}^{d_t} \times \{0,1\}^{r_t}) \Bigr\}  \\
& = \conv\Bigl( \bigl\{(\x,\z,\mu, \y^t) \bigm| (\x,\z,\mu, \y^{t-1}) \in E^{t-1},\ (\z_t,\bdelta_t) \in Q_t \cap \{0,1\}^{d_t} \times \{0,1\}^{r_t} \bigr\}\Bigr),
\end{aligned}
\]
where the second equality holds since $R^{t-1} = \conv(E^{t-1})$ and $Q_t \cap (\R^{d_t} \times \{0,1\}^{r_t})$ is an ideal formulation of $\Delta^{d_t} \cap \{0,1\}^{d_t}$, and the third equality follows from Lemma~\ref{lemma:commonsimplex} since 
\[
\proj_{\z_t}\bigl(Q_t \cap \{0,1\}^{d_t} \times \{0,1\}^{r_t})\bigr) = \proj_{\z_t}(E^{t-1}) =  \Delta^{d_t} \cap \{0,1\}^{d_t},
\]
which is the vertices of the simplex $\Delta^{d_t}$, a set of affinely independent points. Therefore, we conclude that $R^n = \conv(E^n)$. \hfill \Halmos

\section{Proof of Corollary~\ref{cor:supermodular-model-log}}
By Lemma~\ref{lemma:MIP-discrete-log} and the proof of Theorem~\ref{them:supermodular-model}, it suffices to argue that~\eqref{eq:log} is an ideal formulation of $\vertex(\Delta^{d_i})$. To do so, we invoke the ideal formulation for the SOS1 constraint from~\cite{vielma2011modeling}. In particular, Theorem 2 of~\cite{vielma2011modeling} provides an ideal formulation for the vertices of the standard simplex $C^{d_i+1}:=\{\blambda_i \in \R^{d_i+1} \mid \blambda_i \geq 0,\ \sum_{j = 0}^{d_i} \lambda_{ij} = 1\}$: 
\begin{equation}\label{eq:SOS1-log}
\begin{aligned}
	\blambda_i \in \Lambda^{d_i+1} \quad \bdelta_{i} \in \{0,1\}^{\lceil \log_2(d_i+1) \rceil} \quad &\sum_{j \in A_{ik}} \lambda_{ij} \leq \delta_{ik}  \text{ and } \\
 &\sum_{j \notin A_{ik}} \lambda_{ij} \leq 1- \delta_{ik}\quad \for  k = 1, \ldots, \lceil\log_2(n+1) \rceil .
	\end{aligned}\tag{\textsc{SOS1-log}}
\end{equation}
Consider the affine transformation $\Lambda_i$ defined in~\eqref{eq:fromztolambda} and its inverse. Replacing $\blambda_i$ in~\eqref{eq:SOS1-log} with $\Lambda_i(\z_i)$ yields~\eqref{eq:log}. Since the inverse of $\Lambda_i$ maps $\vertex(C^{d_i+1})$ to $\vertex(\Delta^{d_i})$, the validity of~\eqref{eq:SOS1-log} yields that of~\eqref{eq:log}. Next, we argue that~\eqref{eq:log} is ideal. Let $M$ be the LP relaxation of~\eqref{eq:SOS1-log}, and $R$ be that of~\eqref{eq:log}. Clearly, $R = \bigl\{(\z_i,\bdelta_i) \bigm| \blambda_i = \Lambda_i(\z_i),\ (\blambda_i,\bdelta_i) \in M \bigr\}$. Let $(\z_i,\bdelta_i)$ be an extreme point of $R$, and $\blambda_i = \Lambda_i(\z_i)$. It follows that $(\blambda_i, \bdelta_i)$ is an extreme point of $M$ since otherwise, via the invertibility of $\Lambda_i$, we can argue that $(\z_i,\bdelta_i)$ is also not an extreme point of $R$. This, together with the ideality of~\eqref{eq:SOS1-log},   shows that $\bdelta_i$ is binary, proving that~\eqref{eq:log} is ideal. \hfill  \Halmos
\section{Proof of Proposition~\ref{prop:base-linear}}\label{ec:base-linear}
We start with showing that a subset of staircase inequalities in Proposition~\ref{prop:bilinear-sub} suffices to provide an exact MIP formulation for the product of two univariate functions. Here, we focus on treating the case where two functions are non-decreasing. Similar to Proposition~\ref{prop:bilinear-sub}, an exact formulation for  the general case can be obtained using the affine transformation~\eqref{eq:variable-switching-general}. Letting $\psi(j_1,j_2) = f_1(p_{1j_1})f_2(p_{2j_2})$ for $(j_1,j_2) \in \prod_{i=1}^2\{p_{i0}, \ldots, p_{id_i}\}$, we use the following inequalities: 
\begin{equation}\label{eq:ecbistairover}
\begin{aligned}
\mu \leq \psi^{\tau}_\leq(\z):= \psi(0,0) &+ \sum_{j = 1}^{\tau} \bigl(\psi(j,0) - \psi(j-1,0)\bigr)z_{1j} + \sum_{j = \tau +1}^{d_1} \bigl(\psi(j,d_2) - \psi(j-1,d_2)\bigr)z_{1j} \\ 
&  + \sum_{j = 1}^{d_2}\bigl(\psi(\tau,j) -\psi(\tau,j-1)\bigr)z_{2j}     \quad \for \tau \in  0 \cup [d_1],
\end{aligned}    
\end{equation}
and
\begin{equation}\label{eq:ecbistairunder}
\begin{aligned}
\mu \geq \psi^\tau_\geq(\z):= \psi(0,d_2) &+ \sum_{j = 1}^{\tau} \bigl(\psi(j,d_2) - \psi(j-1,d_2)\bigr)z_{1j} + \sum_{j = \tau +1}^{d_1} \bigl(\psi(j,0) - \psi(j-1,0)\bigr)z_{1j} \\ 
&  + \sum_{j = 1}^{d_2}\bigl(\psi(\tau,d_2-j) -\psi(\tau,d_2-j+1)\bigr)(1-z_{2, d_2-j+1})    \quad \for \tau \in  0 \cup [d_1]. 
\end{aligned}    
\end{equation}

\begin{lemma}\label{lemma:bilinear-exact}
Assuming that $f_i(\cdot)$ is a non-decreasing univariate function for $i=1,2$. Then, the unary binarization~\eqref{eq:unary}, and staircase inequalities~\eqref{eq:ecbistairover} and ~\eqref{eq:ecbistairunder} yield an exact MIP formulation for 
\begin{equation}\label{eq:bigraph}
    \Bigl\{(\x,\mu) \Bigm| \mu = f_1(x_1)f_2(x_2),\ x_i \in \{p_{i0}, \ldots, p_{id_i}\} \, \for i =1,2 \Bigr\}.
\end{equation}
\end{lemma}
\noindent{\bf Proof of Lemma~\ref{lemma:bilinear-exact}}. 
By Lemma~\ref{lemma:MIP-discrete}, it suffices to show that 
\[
A:=\bigl\{(\z, \mu) \bigm| \z \in \vertex(\Delta),\ \eqref{eq:ecbistairover} \text{ and }  ~\eqref{eq:ecbistairunder} \bigr\} = \bigl\{(\z, \mu) \bigm| \z \in \vertex(\Delta),\ \mu = \ell_1(\z_1)\ell_2(\z_2) \bigr\}=:B,
\]
where recall that $\ell_i(\z_i) := f_i(p_{i0}) + \sum_{j \in [d_i]}(f_i(p_{ij}) - f_i(p_{ij-1}) )z_{ij}$. Inequalities in \eqref{eq:ecbistairover} and~\eqref{eq:ecbistairunder} are valid since they are obtained from Proposition~\ref{prop:bilinear-sub} by using the following $d_1 + 1$ staircases 
\[
(0,0), \cdots, (\tau,0),(\tau,1), \ldots,(\tau,d_2), \ldots, (d_1,d_2) \qquad \for \tau \in 0 \cup [d_1],
\]
and the linear transformation $\z_2 \mapsto (1-z_{2,d_2}, 1-z_{2,d_2-1}, \ldots, 1 - z_{2,1})$. This shows that $B \subseteq A$. To show $A \subseteq B$, we consider a point  $(\z,\mu) \in A$. Then, $\z = (\v_{1j_1},\v_{2j_2})$ for some $j_1$ and $j_2$, where $\v_{ij}$ is defined as in~\eqref{eq:verticesDelta}. Moreover,  \eqref{eq:ecbistairover} implies that 
\[
\mu \leq \psi_\leq^{j_1}(\z) =  \psi(0,0) + \sum_{j = 1}^{j_1} \bigl(\psi(j,0) - \psi(j-1,0)\bigr)+ \sum_{j = 1}^{j_2}\bigl(\psi(j_1,j) -\psi(j_1,j-1)\bigr)  = \psi(j_1,j_2),
\]
and~\eqref{eq:ecbistairunder} implies that 
\[
\mu \geq \psi_\geq^{j_1}(\z) = \psi(0,d_2) + \sum_{j=1}^{j_1}(\psi(j,d_2)-\psi(j-1,d_2)) + \sum_{j=1}^{d_2-j_2}(\psi(j_1,d_2-j) - \psi(j_1,d_2-j+1)) = \psi(j_1,j_2).
\]
In other words, $ \mu = \psi(j_1,j_2) = \ell_1(\z_1)\ell_2(\z_2)$, showing $(\z,\mu) \in B$. Thus, $A \subseteq B$. \hfill \Halmos

\noindent{\bf Proof of Proposition~\ref{prop:base-linear}:}  By disaggregating each term $(x^i_t-c_t^i) \cdot D^i_t(\x)$, and using the conditions $b^i_m\geq 0$, and $\delta^{ik} \geq 0$ for $k \in A_i$ and $\delta^{ik} \leq 0$ for $k\in B_i$, we obtain the following equivalent formulation of POP:
\begin{subequations}\label{eq:ecbase-linear}
\begin{alignat}{3}
\max \quad &  \sum_{i \in [N]} \sum_{t \in [T]} s^i_t + \sum_{m \in [M_i]} b^i_m \cdot u^i_{tm} \notag
 + \sum_{k \in [N] \setminus i}  \delta^{ik} \cdot v^{ik}_t  -  c^i_t\cdot D_t^i(\x) \\
 \text{s.t.} \quad &~\eqref{eq:brulesMIP} \notag \\
 &s^i_t = x^i_t  (a^i_t-b^i_0x^i_t) \quad \for i \in [N] \text{ and }  t \in [T]  \label{eq:ecbase-linear-1} \\
 & u^i_{tm} \leq x^i_t \cdot x^i_{t-m} \quad \for i \in [N], \, t \in [T] \text{ and } m \in [M_i] \label{eq:ecbase-linear-2} \\
 & v^{ik}_t \leq x^{i}_{t} \cdot x^k_t \quad \for i \in [N],\, t \in [T] \text{ and } k \in A_i \label{eq:ecbase-linear-3} \\
 &v^{ik}_t \geq x^{i}_{t} \cdot x^k_t  \quad \for i \in [N],\, t \in [T] \text{ and } k \in B_i.\label{eq:ecbase-linear-4}
\end{alignat}
\end{subequations}
The constraint~\eqref{eq:ecbase-linear-1} involves a univariate function, and thus can be expressed as the affine function in terms of $\z^i_t$ defined as in~\eqref{eq:linear-self}. By Lemma~\ref{lemma:bilinear-exact}, the bilinear constraints~\eqref{eq:ecbase-linear-2},~\eqref{eq:ecbase-linear-3} and~\eqref{eq:ecbase-linear-4}  are exactly represented as linear constraints~\eqref{eq:linear-model-period},\eqref{eq:linear-model-item-over} and \eqref{eq:linear-model-item-under}, respectively.  Hence, we obtain~\eqref{eq:linear-base}  as an exact MIP formulation of POP. \hfill \Halmos

\section{Proof of Proposition~\ref{prop:base-loglog}}\label{ec:base-loglog}
In this proof, we will use a special case of Lemma~\ref{lemma:bilinear-exact} where $f_2(x_2) =x_2$. In this case, variable $\z_2$ can be projected out to obtain an exact formulation. More specifically, an MIP formulation for the hypograph of $f_1(\cdot)f_2(\cdot)$ is given by~\eqref{eq:unary} with $i=1$, and  staircase inequalities
\begin{equation}\label{eq:ecbistair-simple}
    \begin{aligned}
        \mu \leq \psi(0,0)& + \sum_{j = 1}^{\tau} \bigl(\psi(j,0) - \psi(j-1,0)\bigr)z_{1j} + \sum_{j = \tau +1}^{d_1} \bigl(\psi(j,d_2) - \psi(j-1,d_2)\bigr)z_{1j} \\
        & \qquad + \bigl(\psi(\tau,d_2) - \psi(\tau,0)\bigr) \cdot \frac{x_2-p_{20}}{p_{2d_2}-p_{20}} \quad \for \tau \in \{0,1, \ldots, d_1\}.
    \end{aligned}
\end{equation}
This is obtained by using the relation $x_2 -p_{20} = \sum_{j =1}^{d_2}(p_{2j} - p_{2j-1})z_{2j}$.

By recursively expressing a product of multiple functions into products of two functions, the problem POP can be represented as follows:   
\begin{subequations}
\begin{alignat}{3}
		\max \quad & \sum_{i \in [N]} \sum_{t \in [T]} s_t^i
 + \sum_{k \in [N] \setminus \{i\} }  \delta^{ik}( v^{ik}_t - c^
i_tx^i_t) \notag \\
 \text{s.t.} \quad &~\eqref{eq:brulesMIP} \notag \\
 &u^i_{t1} = (x^i_{t-1})^{b^i_1}  \text{ and }  u^i_{tm} \leq  (x^i_{t-m})^{b^i_m} u^i_{tm-1}       \quad \for i \in [N], \, t \in [T] \text{ and } m  =2, \ldots, M_i  \label{eq:ecbase-loglog-1} \\
 & s^i_t \leq   w^i_t(x^i_t)   u^i_{tM_i}  \quad \for i \in [N] \text{ and } t \in [T] \label{eq:ecbase-loglog-2} \\
 & v^{ik}_t \leq x^{i}_{t} \cdot x^k_t \quad \for i \in [N],\, t \in [T] \text{ and } k \in A_i \label{eq:ecbase-loglog-3} \\
 &v^{ik}_t \geq x^{i}_{t} \cdot x^k_t  \quad \for i \in [N],\, t \in [T] \text{ and } k \in B_i,\label{eq:ecbase-loglog-4}
\end{alignat}
\end{subequations}
where $w^i_t(x^i_t):= (x^i_t-c^i_t) \cdot a^i_t \cdot (x^i_t)^{-b^i_0}$. It suffices to use inequalities in constraints~\eqref{eq:ecbase-loglog-1} and~\eqref{eq:ecbase-loglog-2} because this is a maximization problem and each univariate function is non-negative. Now, using the staircase inequalities in~\eqref{eq:ecbistair-simple}, the constraint~\eqref{eq:ecbase-loglog-1} is formulated as~\eqref{eq:log-crossperiod}. To treat the constraint~\eqref{eq:ecbase-loglog-2}, we need a permutation $\sigma^i_t$ to sort the self-effect $w^i_t(\cdot)$ such that~\eqref{eq:permute-self-effects} is satisfied. Using the affine transformation associated with $\sigma^i_t$ defined as in~\eqref{eq:swtich-self-effect} and  the staircase inequalities in~\eqref{eq:ecbistair-simple}, the constraint~\eqref{eq:ecbase-loglog-2} is formulated as~\eqref{eq:self-crossperiod}. The constraints~\eqref{eq:ecbase-loglog-3} and~\eqref{eq:ecbase-loglog-4} are treated as in the proof of Proposition~\ref{prop:base-linear}. Hence, we obtain~\eqref{eq:formulation-multi} as an exact formulation of POP. \hfill \Halmos

\section{Proof of Theorem~\ref{thm:ppo-lp}}
We will need to use the following decomposition lemma, and a similar result has been obtained in~\citep{tawarmalani2010inclusion}. 

\begin{lemma}\label{lemma:decompositionsupermodular}  
For $i \in [r]$, let $g_i: \vertex(\Delta) \to \R$ be a submodular function, and for $j \in [m]$, let $h_j: \vertex(\Delta) \to \R$ be a supermodular function, where $\Delta:=\Delta^{d_1} \times \cdots \times \Delta^{d_n}$. Consider a set $\mathcal{S}$ defined as the intersection of the epigraph of $g_i$ and the hypograph of $h_j$, that is, $\mathcal{S} := \bigl\{(\z, \bmu, \boldsymbol{\nu}) \bigm| \mu_i 
\geq g_i(\z) \for i \in [r] \text{ and } \nu_j \leq h_j(\z) \for  j \in [m],\ \z \in \vertex(\Delta) \bigr\}$.
Then, $\conv(\mathcal{S}) = R$, where
\[
\begin{aligned}
    R :&= \Bigl\{(\z, \bmu, \boldsymbol{\nu}) \Bigm| (\z, \mu_i) \in \conv(\epi(g_i)) 
 \for i \in [r],\ (\z, \nu_j) \in \conv\bigl(\hypo(h_j)\bigr) \for  j \in [m],\   \z \in \Delta \Bigr\} \\
&=    \Bigl\{(\z, \bmu, \boldsymbol{\nu}) \Bigm| \mu_i 
\geq g_i^\pi(\z) \for i \in [r]\; \pi \in \Pi \text{ and } \nu_j \leq h^\pi_j(\z) \for  j \in [m]\;  \pi \in \Pi,\ \z \in \Delta \Bigr\},
\end{aligned}
\]
where for a staircase $\pi$, $g^\pi_i(\cdot)$ and $h^\pi_j(\cdot)$ are defined as in~\eqref{ec:staircaseineq}. 
\end{lemma}
\noindent{\bf Proof of Lemma~\ref{lemma:decompositionsupermodular}}. The second equality holds due to Lemma~\ref{lemma:staircase}. Now, we show that $\conv(\mathcal{S}) = R$. Clearly, $\conv(\mathcal{S}) \subseteq R$ since $\mathcal{S} \subseteq R$ and $R$ is convex. To prove the opposite direction, we consider a point $(\bar{\z}, \bar{\bmu}, \bar{\boldsymbol{\nu}}) \in R$. To show that this point can be expressed a convex combination of points in $\mathcal{S}$, we follow the argument in the proof of Lemma~\ref{lemma:staircase}. There exits a staircase $\omega$ such that $\bar{\z} \in \Upsilon_\omega$ and 
\begin{equation}\label{eq:decompositionsupermodular}
    \boldsymbol{g}^\omega(\v) = \boldsymbol{g}(\v)   \quad \text{ and }\quad \boldsymbol{h}^\omega(\v) = \boldsymbol{h}(\v)  \quad  \for \v \in \vertex(\Upsilon_\omega). 
\end{equation}
Then, there exist a convex multiplier $\{\lambda_{\v}\}$ such that  $\bar{\z} = \sum_{\v \in \vertex(\Upsilon_\omega)} \lambda_{\v} \v$, and 
\[
\begin{aligned}
 \bigl(\boldsymbol{g}^\omega(\bar{\z}),\boldsymbol{h}^\omega(\bar{\z})\bigr) & = \sum_{\v \in \vertex(\Upsilon_\omega)} \lambda_{\v}\bigl(\boldsymbol{g}^\omega(\v),\boldsymbol{h}^\omega(\v) \bigr) = \sum_{\v \in \vertex(\Upsilon_\omega)} \lambda_{\v}\bigl(\boldsymbol{g}(\v),\boldsymbol{h}(\v)\bigr).
        \end{aligned}
\]
Let $\boldsymbol{\epsilon}:= \boldsymbol{g}^\omega(\bar{\z}) - \bar{\boldsymbol{\mu}}$ and $\boldsymbol{\varepsilon}:=\boldsymbol{h}^\omega(\bar{\z}) - \bar{\boldsymbol{\nu}}$.  Then,
\[
(\bar{\z},\bar{\bmu},\bar{\boldsymbol{\nu}}) = \bigl(\bar{\z}, \boldsymbol{g}^\omega(\bar{\z})-\boldsymbol{\epsilon}, \boldsymbol{h}^\omega(\bar{\z})-\bar{\boldsymbol{\nu}}\bigr) = \sum_{\v \in \vertex(\Upsilon_\omega)}\lambda_{\v}\bigl(\v, \boldsymbol{g}(\v)-\boldsymbol{\epsilon}, \boldsymbol{h}(\v) - \boldsymbol{\varepsilon}\bigr) \in \conv(\mathcal{S}),
\]
where the second equality holds due to the convex decomposition in~\eqref{eq:decompositionsupermodular}, and the last relation holds due to $\boldsymbol{\epsilon} \leq 0$ and $\boldsymbol{\varepsilon} \geq 0$. \hfill \Halmos

Now, using Lemma~\ref{eclemma:MIP-discrete} and Lemma~\ref{lemma:decompositionsupermodular}, we obtain a generalization of Theorem~\ref{them:supermodular-model}, yielding an ideal formulation for the hypograph of a vector of supermodular compositions. 
\begin{corollary}\label{claim:simutanuoushull}
    Let $\phi \mcirc f : \X \to \R^m$ be a vector of composite functions defined as 
            \[
    (\phi \mcirc f)(\x) = \Bigl(\phi_1\bigl(f_{11}(x_1), 
    \ldots, f_{1n}(x_n) \bigr), \ldots, \phi_m\bigl(f_{m1}(x_1), 
    \ldots, f_{mn}(x_n) \bigr)\Bigr) \qquad \for \x \in \X, 
    \]
    where we assume that $f_{ki}(\cdot)$ is a non-decreasing function. 
    For $k \in [m]$, let $\psi_k(j_1, \ldots,j_n) = (\phi_k \mcirc f_k)(p_{1j_1}, \ldots, p_{nj_n})$ for $\boldsymbol{j} \in \prod_{i\in [n]}\{0, \ldots, d_i\}$. Then, an ideal formulation for the hypograph of $\phi \mcirc f$, that is, $\{(\x,\bmu) \mid \mu_{k} \leq (\phi_k \mcirc f_k)(\x) \for k\in [m],\ \x \in \X \}$, is given as follows
    \[
   \begin{aligned}
    \Bigl\{(\x, \boldsymbol{\mu}, \z) \Bigm| (x_i,\z_i) & \in\eqref{eq:unary} \, \for i \in [n] \\
       & \mu_k \leq	\psi_k(\j^0_\pi) + \sum_{t  \in [D]} \bigl( \psi_k(\j^t_\pi) - \psi_k(\j^{t-1}_\pi) \bigr) \cdot z_{\pi_t} \, \for k \in [m] \text{ and } \pi \in \Pi \Bigr\}.
    \end{aligned}
    \]
\end{corollary}

\noindent{\bf Proof of Theorem~\ref{thm:ppo-lp}:} First, we consider the case where the demand model is as in~\eqref{eq:demand-linear} and all items are substitutable. Under this case, POP can be expressed as follows: 
\begin{subequations}
\begin{alignat}{3}
\max \quad &  \sum_{i \in [N]} \sum_{t \in [T]} s^i_t + \sum_{m \in [M_i]} b^i_m \cdot u^i_{tm} \notag
 + \sum_{k \in [N] \setminus i}  \delta^{ik} \cdot v^{ik}_t  -  c^i_t\cdot D_t^i(\x) \\
 \text{s.t.} \quad & x^i_t \in \{p^i_0, \ldots, p^i_{d_i}\} \quad  \for i \in [N] \text{ and } t \in [T]   \notag \\
 &s^i_t \leq  x^i_t(a^i_t-b^i_0x^i_t) \quad \for i \in [N] \text{ and }  t \in [T]  \label{eq:ec-ppo-lp-bi-1} \\
 & u^i_{tm} \leq x^i_t \cdot x^i_{t-m} \quad \for i \in [N] \, t \in [T] \text{ and } m \in [M_i] \label{eq:ec-ppo-lp-bi-2} \\
 & v^{ik}_t \leq x^{i}_{t} \cdot x^k_t \quad \for i \in [N]\, t \in [T] \text{ and } k \in [N] \label{eq:ec-ppo-lp-bi-3}. 
\end{alignat}
\end{subequations}
Since functions appearing in~\eqref{eq:ec-ppo-lp-bi-1},~\eqref{eq:ec-ppo-lp-bi-2} and~\eqref{eq:ec-ppo-lp-bi-3} are compositions of supermodular functions and non-decreasing unvariate functions, it follows from Corollary~\ref{claim:simutanuoushull} that individually modeling each constraint using staircase inequalities yields  an ideal formulation, and thus its natural LP relaxation solves POP. 

Second, we consider the case where the demand model is as in~\eqref{eq:demand-loglog} and  all items are substitutable. Under this case, POP can be expressed as follows: 
\begin{subequations}
\begin{alignat}{3}
		\max \quad & \sum_{i \in [N]} \sum_{t \in [T]} s_t^i
 + \sum_{k \in [N] \setminus \{i\} }  \delta^{ik}( v^{ik}_t - c^
i_tx^i_t) \notag \\
 \text{s.t.} \quad & x^i_t \in \{p^i_0, \ldots, p^i_{d_i}\} \quad  \for i \in [N] \text{ and } t \in [T] \notag \\
  & s^i_t \leq   w^i_t(x^i_t) \cdot \prod_{m \in [M_i]}  (x^i_{t-m})^{b^i_m}   \quad \for i \in [N] \text{ and } t \in [T] \label{eq:ec-ppo-lp-loglog-2} \\
 & v^{ik}_t \leq x^{i}_{t} \cdot x^k_t \quad \for i \in [N]\, t \in [T] \text{ and } k \in [N]. \label{eq:ec-ppo-lp-loglog-3} 
 \end{alignat}
\end{subequations}
Since the self-effect $w^i_t(\cdot)$ is non-negative and is assumed to be non-decreasing, the function appearing in~\eqref{eq:ec-ppo-lp-loglog-2} is a composition of a supermodular function and a non-decreasing univariate function. It follows from Corollary~\ref{claim:simutanuoushull} that individually modeling each constraint using staircase inequalities yields  an ideal formulation, and thus its natural LP relaxation solves POP. \hfill \Halmos

\section{Computational Comparison on Root Node Relaxations}\label{app:bounds}
Here, we compare the root node relaxations of the nonlinear model~$\mathsf{NL}$ and the formulations~$\mathsf{Stair}$-K for $K\in \{0,2,5,7\}$. The root node gap quantifies the relative tightness of the relaxation and is computed as:
\[
\frac{v_{\mathsf{model}} - v^*}{v^*} \times 100\% \quad \text{ for } \mathsf{model} \in \{\mathsf{NL},\mathsf{Stair\text{-}0},\mathsf{Stair\text{-}2},\mathsf{Stair\text{-}5},\mathsf{Stair\text{-}7} \},
\]
where $v_{\mathsf{model}}$ is the optimal objective value of the root node relaxation for each $\mathsf{model}$, and $v^*$ denotes the objective value of the best known feasible solution. 

Table~\ref{tab:rootrlx} reports the average root gaps across different parameter configurations. Overall, the root node gaps of our formulations $\mathsf{Stair}$-K are significantly smaller than that for the nonlinear model $\mathsf{NL}$, with all root gaps consistently below $3\%$. This reduced gap demonstrates that our formulations are significantly tighter than the nonlinear model $\mathsf{NL}$. In particular, model $\mathsf{Stair}$-5 achieves the best tradeoff between formulation strength and size, making it the most effective choice for our computational experiments. 

\renewcommand{\arraystretch}{1}
\begin{table}[htbp]
\TABLE 
{{Average root gaps and sizes of alternative formulations}\label{tab:rootrlx}}
{\begin{tabular}{cc rrrrr rrr}
\toprule
\multicolumn{2}{c}{Instance} & \multicolumn{5}{c}{Root gap (\%)} &  \multicolumn{3}{c}{\# cuts} \\
\cmidrule(r){1-2} \cmidrule(lr){3-7} \cmidrule(lr){8-10}

$N$ & $T$ &  \textsf{NL} & \textsf{Stair}-0 & \textsf{Stair}-2 & \textsf{Stair}-5  &  \textsf{Stair}-7 & \textsf{Stair}-2 & \textsf{Stair}-5  &  \textsf{Stair}-7 \\
\midrule
10 & 5 
& 157.544 & 0.897 & 0.638 & 0.633 & 0.633 & 7.9 & 9.3 & 9.3  \\
10 & 10 
& 131.749 & 1.72 & 1.263 &    1.261 &   1.261 & 24.9 & 26.3 & 26.3 \\
\midrule
50 & 5 
& 139.466 & 1.312 & 0.931 & 0.925 & 0.925 & 61.4 & 70.4 & 70.8  \\
50 & 10 
& 136.53 & 1.72 & 1.184 & 1.178 &   1.178 & 170.1 & 193.3 & 193.3 \\
\midrule
75 & 5 
& 137.12 & 1.434 & 1.111 & 1.110 &   1.110 & 107.7 & 115.3 & 115.3 \\
75 & 10 
& 117.97 & 1.563 & 1.124 & 1.117 &   1.117 & 320 & 356.8 & 360.4 \\
\midrule
100 & 5 
& 136.92 & 1.650 & 1.311 & 1.311 &  1.311 & 169.2 & 187.8 & 188.5 \\
100 & 10 
& 123.79 & 2.094 & 1.762 & 1.759 &1.759 & 379.3 & 414.8 & 416.3 \\
\midrule
125 & 5 
& 132.80 & 2.081 & 1.755 & 1.754 &   1.754  & 214.9 & 226.5 & 227.2\\
125 & 10 
& 127.76 & 2.359 & 1.996 & 1.993 &   1.993 & 424.5 & 450.5 & 451.1 \\

\bottomrule
\end{tabular}}
{}
\end{table}
\section{Proof of Theorem~\ref{them:Lnatural}}
We denote by $\conv(\phi)$ the convex envelope of $\phi(\cdot)$ over $H:=\prod_{i \in [n]}[0,d_i]$. We begin with proving the first statement. Suppose that $\phi(\cdot)$ is a submodular function. We observe that the epigraph of the staircase extension $\check{\phi}$ can be expressed as follows:
\[
\begin{aligned}
\epi(\check{\phi}) &= \bigl\{ (\x, \mu) \in H \times \R  \bigm| \mu \geq \check{\phi}(\x) \bigr\} \\
&= \bigl\{(\x, \mu) \in H \times \R \bigm| \exists \z \, \st \, (\mu,\x,\z) \in \Phi  \bigr\} \\
&=  \conv\bigl(\epi(\phi)\bigr),
\end{aligned}
\]
where the first and second equality hold by definition. To see the last equality, we observe that, by Theorem~\ref{them:supermodular-model}, $\bigl\{(\mu,\x,\z) \in  \Phi \bigm| \z \in \{0,1\}^{\sum_{i \in [n]}d_i} \bigr\}$ is an ideal MIP formulation of the epigraph of $\phi(\cdot)$, and thus the projection of its LP formulation $\Phi$ onto the space of $(\x ,\mu)$ variables yields the convex hull of $\epi(\phi)$. Since $\epi(\check{\phi}) = \conv\bigl(\epi(\phi)\bigr)$, by Theorem 5.6 in~\cite{rockafellar2015convex}, we obtain  that $\check{\phi}(\x) = \conv(\phi)(\x)$ for $\x \in H$. To show the reverse direction, we assume that $\phi(\cdot)$ is not submodular. In other words, there exist two distinct points $\p$ and $\boldsymbol{q}$ such that  $\p \vee \boldsymbol{q}$ and $\p \wedge \boldsymbol{q}$ are distinct and  $\phi(\p \vee \boldsymbol{q} ) + \phi(\p \wedge \boldsymbol{q} ) > \phi(\p) + \phi(\boldsymbol{q})$. Let $\pi$ be a staircase passing through $\p \wedge \boldsymbol{q}$, $\p$, and $\p \vee \boldsymbol{q}$, and consider the affine function defined by the staircase
\[
\phi^\pi(\z):= \phi(\j^0_\pi) + \sum_{s \in [D]}\bigl(\phi(\j^s_\pi) - \phi(\j^{s-1}_\pi)\bigr)z_{\pi_t}
\]
 where $\j^0_\pi, \ldots,\j^{D}_\pi$ is the point representation of the staircase $\pi$ defined as in~\eqref{eq:pointrep}. Clearly, $\phi^\pi(\v_{1j_1}, \ldots, \v_{nj_n}) = \phi(\j)$ for every $\j = (j_1, \ldots,j_n) \in \{\j^0_\pi, \ldots, \j^D_\pi\}$. Thus, we can obtain that  
\[
\begin{aligned}
\phi (\boldsymbol{q}) &< \phi(\p \vee \boldsymbol{q} ) + \phi(\p \wedge \boldsymbol{q} ) - \phi(\p) \\
& = \phi^\pi(\v_{1 \max\{p_1,q_1\} }, \ldots, \v_{n \max\{p_n,q_n\} } )	+  \phi^\pi(\v_{1 \min\{p_1,q_1\} }, \ldots, \v_{n \min\{p_n,q_n\} } ) - \phi^\pi(\v_{1p_1}, \ldots, \v_{np_n} ) \\
& = \phi^\pi(\v_{1q_1}, \ldots, \v_{nq_n} ) \\
& \leq \check{\phi}(\boldsymbol{q}), 
\end{aligned}
\]
where the first inequality follows from the assumption that $\phi$ is not submodular, the first equality holds since $\p \wedge \boldsymbol{q}$, $\p$ and $\p \vee \boldsymbol{q}$ belong to $\j^0_\pi, \ldots, \j^D_\pi$,  the second equality holds since $\phi^\pi(\cdot)$ is an affine function and $\v_{i\max\{p_i,q_i\}} + \v_{i\min\{p_i,q_i\}} =\v_{ip_i} + \v_{iq_i}$, and the second inequality holds by the definition of $\check{\phi}$. Therefore, $\check{\phi}(\cdot)$ is not an underestimator of $\phi(\cdot)$, and thus is not the convex envelope of $\phi(\cdot)$. 

Next, we prove the second statement. By definition, $\hat{\phi}(\x) \geq \conv(\phi)(\x)$ for $\x \in H$. If $\phi(\cdot)$ is $\lnatural$-convex then by Lemma~\ref{lemma:murota} the locally Lov\'asz extension $\hat{\phi}(\cdot)$ is convex. This, together with $\hat{\phi}(\x) \leq \phi(\x)$ for $\x \in \prod_{i \in [n]}\{0,1, \ldots d_i\}$, implies that $\hat{\phi}(\x) \leq \conv(\phi)(\x)$ for $\x \in H$. Moreover, since a $\lnatural$-convex function is a submodular function, by the first statement, we obtain that  $\conv(\phi)(\x) = \check{\phi}(\x)$ for $\x \in H$. Therefore, we can conclude that $\check{\phi}(\x) = \hat{\phi}(\x)$ for $x \in H$. Now, suppose that $\check{\phi}(\x) = \hat{\phi}(\x)$ for $x \in H$. Since the staircase extension is always convex, it follows that the locally Lov\'asz extension $\hat{\phi}(\cdot)$ is convex. Thus, by Lemma~\ref{lemma:murota}, $\phi(\cdot)$ is $\lnatural$-convex. \hfill \Halmos

\end{document}